\documentclass[a4paper]{article}
\usepackage{theorem,amsmath,amssymb,amscd}
\usepackage[small]{titlesec}

\theoremstyle{change} \allowdisplaybreaks 

\newcommand{\R}{\mathbb{R}}
\newcommand{\C}{\mathbb{C}}
\newcommand{\Z}{\mathbb{Z}}
\newcommand{\Q}{\mathbb{Q}}
\newcommand{\A}{\mathbb{A}}

\newcommand{\SH}{\mathfrak h}
\newcommand{\HH}{\mathbb{H}}
\newcommand{\St}{\mathrm{St}}
\newcommand{\GL}{\mathrm{GL}}
\newcommand{\GU}{\mathrm{GU}}

\newcommand{\SL}{\mathrm{SL}}
\newcommand{\SO}{\mathrm{SO}}
\newcommand{\SU}{\mathrm{SU}}
\newcommand{\GSp}{\mathrm{GSp}}
\newcommand{\SSp}{\mathrm{Sp}}
\newcommand{\OF}{\mathfrak{o}}
\newcommand{\p}{\mathfrak{p}}
\renewcommand{\P}{\mathfrak{P}}

\newcommand{\tr}{{\rm tr}}
\newcommand{\AI}{\mathcal{AI}}
\newcommand{\mat}[4]{{\setlength{\arraycolsep}{0.5mm}\left[
\begin{array}{cc}#1&#2\\#3&#4\end{array}\right]}}
\newcommand{\qed}{\hspace*{\fill}\rule{1ex}{1ex}}
\def\qdots{\mathinner{\mkern1mu\raise0pt\vbox{\kern7pt\hbox{.}}\mkern2mu
\raise3.4pt\hbox{.}\mkern2mu\raise7pt\hbox{.}\mkern1mu}}

\newcommand{\nl}{

\vspace{2ex}}
\newcommand{\nll}{

\vspace{1ex}}

\newtheorem{thm}{Theorem.}[subsection]
\newtheorem{theorem}{Theorem.}[subsection]
\newtheorem{lem}[thm]{Lemma.}
\newtheorem{lemma}[thm]{Lemma.}

\begin{document}

\begin{center}
{\Large Integral Representation for $L$-functions for $\GSp_4\times \GL_2$}

\vspace{2ex} Ameya Pitale\footnote{Department of Mathematics, University of Oklahoma,
Norman, OK 73019-0315, {\tt ameya@math.ou.edu}},
Ralf Schmidt\footnote{Department of Mathematics, University of Oklahoma,
Norman, OK 73019-0315, {\tt rschmidt@math.ou.edu}}

\vspace{3ex}
\begin{minipage}{80ex}
{\sc Abstract.} Let $\pi$ be a cuspidal, automorphic representation of $\GSp_4$
attached to a Siegel modular form of degree $2$.
We refine the method of Furusawa \cite{Fu} to obtain an integral
representation for the degree-$8$ $L$-function $L(s,\pi\times\tau)$, where
$\tau$ runs through certain cuspidal, automorphic representation of $\GL_2$.
Our calculations include the case of square-free level for the $p$-adic
components of $\tau$, and a wide class of archimedean types including Maa{\ss} forms.
As an application we obtain a special value result for $L(s,\pi\times\tau)$.

\end{minipage}
\end{center}
\section{Introduction}
Let $\pi = \otimes \pi_{\nu}$ and $\tau = \otimes \tau_{\nu}$ be
irreducible, cuspidal, automorphic representations of $\GSp_4(\A)$
and $\GL_2(\A)$, respectively. Here, $\A$ is the ring of adeles of
a number field $F$. We want to investigate the degree eight
twisted $L$-functions $L(s, \pi \times \tau)$ of $\pi$ and $\tau$,
which are important for a number of reasons. For
example, when $\pi$ and $\tau$ are obtained from holomorphic
modular forms, then Deligne \cite{De1} has conjectured that a
finite set of special values of $L(s, \pi \times \tau)$ are
algebraic up to certain period integrals. Another very important
application is the conjectured Langlands functorial transfer of
$\pi$ to an automorphic representation of $\GL_4(\A)$. One
approach to obtain the transfer to $\GL_4(\A)$ is to use the
converse theorem due to Cogdell and Piatetski-Shapiro \cite{CPS},
which requires precise information about the $L$-functions $L(s, \pi
\times \tau)$.

In the special case that $\pi$ is generic, Asgari and Shahidi
\cite{ASh} have been successful in obtaining the above transfer
using the converse theorem. They analyze the twisted $L$-functions
using the Langlands-Shahidi method. In this method, one has to
consider a larger group in which $\GSp_4$ is embedded and then use
the representation $\pi$ to construct an Eisenstein series on the
larger group. Then the $L$-functions are obtained in the constant
and non-constant terms of the Eisenstein series. Unfortunately,
this method only works when $\pi$ is generic. It is known that if
$\pi$ is obtained from a holomorphic Siegel modular form then it
is not generic.

Another method to understand $L$-functions is via integral
representations. For this method one constructs an integral
that is Eulerian, i.e., one that can be written as an infinite
product of local integrals, $Z(s) = \prod_{\nu} Z_{\nu}(s)$.
Then the local integrals are computed to obtain the local
$L$-functions. In many of the constructions, the local
calculations are done only when all the local data is unramified.
This gives information about the partial $L$-functions, which
already leads to remarkable applications. The calculations for the
ramified data are unfortunately often very involved and not
available in the literature. (For more on integral representations
of $L$-functions, see \cite{GPR1}, \cite{GPR2}, \cite{PR}.)

In the $\GSp_4 \times \GL_2$ case, Novodvorsky, Piatetski-Shapiro
and Soudry (see \cite{N}, \cite{PS}, \cite{PS-S}) were the first
ones to construct integral representations for $L(s, \pi \times
\tau)$. Their constructions were for the special case when $\pi$
is either generic or has a special Bessel model. Examples of
Siegel modular forms which do not have a special Bessel model have
been constructed by Schulze-Pillot \cite{S-P}. The first
construction of an integral representation for $L(s, \pi \times
\tau)$ with no restriction on the Bessel model of $\pi$ is the
work of Furusawa \cite{Fu}. In this remarkable paper, Furusawa
embeds $\GSp_4$ in a unitary group $\GU(2,2)$ and constructs an
Eisenstein series on $\GU(2,2)$ using the $\GL_2$ representation
$\tau$. He then integrates the Eisenstein series against a vector
in $\pi$. He shows that this integral is Eulerian and, when the
local data is unramified, he computes the local integral to obtain
the local $L$-function $L(s, \pi_{\nu} \times \tau_{\nu})$ up to a
normalizing factor. He also calculates the archimedean
integral for the case that both $\pi$ and $\tau$ are holomorphic of the same weight.
Thus, Furusawa obtains an integral representation for the
completed $L$-function $L(s,\pi \times \tau)$ in the case when
$\pi$ and $\tau$ are obtained from holomorphic modular forms of
full level and same weight. He uses this to obtain a special value
result, which fits into the context of Deligne's conjectures,
and to prove meromorphic continuation
and functional equation for the $L$-function. The main limitation
of \cite{Fu} is that, if we fix a Siegel modular form, then the
results allow us to obtain information on a very small family of
twists only, namely those coming from elliptic modular forms of
full level and the same weight as the Siegel modular form, which
is a finite dimensional vector space.

For the applications that we discussed above, we need twists of
$\pi$ by all representations $\tau$ of $\GL_2$, i.e., twists by
all $\GL_2$ modular forms, holomorphic or non-holomorphic, of
arbitrary weight and level. For this purpose, one needs to compute
the non-archimedean local integral obtained in \cite{Fu} when the
local representation $\tau_{\nu}$ is ramified. Also, one needs to
extend Furusawa's archimedean calculation to include more general
archimedean representations.

In this paper, we will compute the local non-archimedean integral
from \cite{Fu} in a mildly ramified case, namely
when $\tau_{\nu}$ is an unramified twist of the Steinberg representation. We will also
compute the archimedean integral for a larger family of
archimedean representations $\tau_{\infty}$.

Before we state the results of this paper, let us recall the
integral representation of \cite{Fu} in some more detail. Let $L$ be a quadratic
extension of the number field $F$, and let $\GU(2,2)$ be the
unitary group defined using the field $L$. Let $P$ be the standard
maximal parabolic subgroup of the unitary group $\GU(2,2)$ with a
non-abelian radical. Given an irreducible, admissible
representation $\tau$ of $\GL_2(\A)$ and suitable characters
$\chi$ and $\chi_0$ of $\A_L^\times$, one considers an induced representation
$I(s,\chi,\chi_0,\tau)$ from $P$ to $\GU(2,2)$, where $s$ is a complex parameter.
Let $f(g,s)$ be an analytic family in $I(s,\chi,\chi_0,\tau)$.
Define an Eisenstein series on $\GU(2,2)$ by the formula
$$
 E(g,s) = E(g,s;f) = \sum\limits_{\gamma\in P(F)\backslash
 \GU(2,2)(F)}f(\gamma g, s),\qquad g \in \GU(2,2)(\A).
$$
For an automorphic form $\phi$ in the space of $\pi$, consider the integral
\begin{equation}\label{intro-main-int}
 Z(s)=Z(s,f,\bar{\phi}) = \int\limits_{Z(\A)\GSp_4(F)\backslash
 \GSp_4(\A)}E(h,s;f)\bar{\phi}(h)dh.
\end{equation}
In \cite{Fu}, Furusawa has shown that these integrals have the
following two important properties.
\begin{enumerate}
 \item There is a ``basic identity''
  \begin{equation}\label{intro-basic-identity}
   Z(s)=\int\limits_{R(\A)\backslash \GSp_4(\A)}W_f(\eta h,s)B_{\bar{\phi}}(h) dh,
  \end{equation}
  where $R\subset\GSp(4)$ is a Bessel subgroup of the Siegel parabolic subgroup,
  $\eta$ is a certain fixed element, $B_{\bar{\phi}}$ corresponds to $\bar{\phi}$
  in the Bessel model for $\pi$, and $W_f$ is a function on
  $\GU(2,2)$ obtained from the Whittaker model of $\tau$ and depending
  on the section $f$ used to define the Eisenstein series.
 \item $Z(s)$ is Eulerian, i.e.,
  \begin{equation}\label{intro-eulerian-int}
   Z(s) = \prod\limits_{\nu}Z_{\nu}(s) = \prod\limits_{\nu}\:
   \int\limits_{R(F_{\nu})\backslash \GSp_4(F_{\nu})}W_{\nu}(\eta h, s)B_{\nu}(h)\,dh.
  \end{equation}
\end{enumerate}
In Theorem \ref{localintegralmaintheorem} below we show that if
$\tau_{\nu}$ is mildly ramified then the local integral can be computed
to give $L(3s+\frac 12,\tilde\pi_{\nu} \times \tilde\tau_{\nu})$
up to a normalizing factor.

{\bf Theorem 1.}
 \emph{Let $F_{\nu}$ be a non-archimedean local field with characteristic zero. Let
 $\pi_{\nu}$ be an unramified, irreducible, admissible
 representation of $\GSp_4(F_{\nu})$. Let $\tau_{\nu}$ be an
 unramified twist of the Steinberg representation. Then we can make a choice of
 vectors $W_{\nu}$ and $B_{\nu}$ such that the local integral in
 (\ref{intro-eulerian-int}) is given by
 $$
  Z_{\nu}(s) = \frac{q(q-1)}{(q+1)(q^4-1)}\Big(1-\Big(\frac{L_{\nu}}{\p}\Big)q^{-1}\Big)
  \frac{L(3s+\frac 12,\tilde\pi_{\nu} \times \tilde\tau_{\nu})}
  {L(3s+1, \tau_{\nu} \times \AI(\Lambda_{\nu})\times (\chi_{\nu}|_{F_{\nu}^{\times}}))}.
 $$
 Here, $q$ is the cardinality of the residue class field of $F_{\nu}$, $\Lambda_{\nu}$ is
 the Bessel character on $L_{\nu}^{\times}$ used to define the Bessel model $B_{\nu}$, and
 $\AI(\Lambda_{\nu})$ is the representation of $\GL_2(F_{\nu})$ obtained from
 $\Lambda_{\nu}$ by automorphic induction.
}

We point out that the ramified calculation is not a trivial
generalization of the unramified calculation in \cite{Fu}. There
are two main steps. First is the choice of the vector $W_{\nu}$ --
for several obvious choices $Z_\nu(s)$ evaluates to zero. This
choice depends crucially on the underlying number theory.
Secondly, the actual computation of the local integral is
complicated and depends heavily on the structure theory of the
groups involved. We will explain this in detail in Sect.\
\ref{non-arch-ram-sec}.

In Theorem \ref{archmaintheorem}, we compute the local archimedean
integral in the following cases:
\begin{enumerate}
\item  $\pi_{\infty}$ is the holomorphic discrete series
representation of $\GSp_4(\R)$ with trivial central character and
Harish-Chandra parameter $(l-1,l-2)$.

\item $\tau_{\infty}$ is either a principal series representation
of $\GL_2(\R)$ whose $K$-types have the same parity as $l$ or is a
holomorphic discrete series representation of $\GL_2(\R)$ with
lowest weight $l_2$ satisfying $l_2 \leq l$ and $l_2 \equiv l
\pmod{2}$.
\end{enumerate}
This extends the calculations in \cite{Fu}, where $\tau_{\infty}$
is only allowed to be a holomorphic discrete series representation
with lowest weight $l$.

Putting together the local computations we get the following global result in Theorem \ref{main-global-thm}.

{\bf Theorem 2.}
 \emph{Let $\Phi$ be a cuspidal Siegel  eigenform of weight $l$ with
 respect to $\SSp_4(\Z)$. Let $N$ be a square-free, positive integer. Let $f$
 be a cuspidal Maa{\ss} eigenform of weight $l_1 \in \Z$ with
 respect to $\Gamma_0(N)$. If (the adelic function corresponding to)
 $f$ lies in a holomorphic discrete
 series representation with lowest weight $l_2$, then assume that $l_2 \leq l$.
 Let $\pi_{\Phi}$ and $\tau_f$ be the corresponding cuspidal
 automorphic representations of $\GSp_4(\A_{\Q})$ and
 $\GL_2(\A_{\Q})$, respectively. Then a choice of local vectors can be made
 such that the global integral $Z(s)$ defined in (\ref{intro-main-int}) is given by
 \begin{equation}\label{intro-integral-l-fn-formula}
  Z(s) = \kappa_{\infty} \kappa_N
  \frac{L(3s+\frac 12, \pi_{\Phi} \times \tau_f)}{\zeta(6s+1)
  L(3s+1, \tau_f\times\AI(\Lambda))},
 \end{equation}
 where $\kappa_{\infty}$ and $\kappa_N$ are obtained from the local computations.
}

Note that the above theorem still gives information on the
twisted $L$-functions for a smaller family of representations
$\tau$ than is required for the application of the converse
theorem mentioned earlier. For a general representation $\tau$ the
local calculations are conceptually the same but the calculations
are much more complicated. This is work in progress and will be a
subject of a future paper.

Using (\ref{intro-integral-l-fn-formula}), we get the following special value result
in Theorem \ref{special values thm}.

{\bf Theorem 3.}
 \emph{Let $\Phi$ be a cuspidal Siegel  eigen-form of weight $l$ with
 respect to $\SSp_4(\Z)$. Let $N$ be a square-free, positive integer. Let
 $\Psi$ be a holomorphic, cuspidal Hecke eigen-form of weight $l$
 with respect to $\Gamma_0(N)$. Then
 $$
  \frac{L(\frac l2 - 1, \pi_{\Phi} \times \tau_{\Psi})}{\pi^{5l-8}
  (\Phi,\Phi)_2(\Psi,\Psi)_1} \in \overline{\Q}.
 $$
}

Note that in \cite{BH}, using a completely different method,
special value results, in the spirit of Deligne's conjectures,
were proven under the assumption that $\Psi$ is a cusp form with
respect to $\SL(2,\Z)$ with weight $k \leq 2l-2$, where $l$ is the
weight of the Siegel modular form $\Phi$. Since the results of
\cite{BH} cannot be applied to modular forms with respect to
congruence subgroups, there is no overlap of \cite{BH} with this
paper.

This paper is organized as follows. In Sect.\ \ref{defn-sec} we
make the basic definitions and describe the setup for the local
integrals from (\ref{intro-eulerian-int}) for a non-archimedean local field $F$
of characteristic zero or $F = \R$. We use the fact that the basic local setup is uniform
and can be stated in full generality. The main input of the local
integrals are the choices of the functions $W$ and $B$ from
(\ref{intro-eulerian-int}). In Sects.\ \ref{non-arch-ram-sec} and
\ref{arch-sec} we consider the non-archimedean and archimedean
case, respectively. We make the choice of the appropriate functions
$W$ and $B$ and compute the local integrals. In Sect.\
\ref{modular forms sec}, we consider the global situation
corresponding to modular forms on $\GSp_4$ and $\GL_2$. We use the
local calculations from Sects.\ \ref{non-arch-ram-sec} and
\ref{arch-sec} to obtain an integral representation for the global
$L$-function. Finally, in Sect.\ \ref{special values sec}, we use
the global theorem to obtain a special values result.
\section{General setup}\label{defn-sec}
In this section, we give the basic definitions and set up the data
required to compute the local integrals. Let $F$
be a non-archimedean local field of characteristic zero, or
$F=\R$. We fix three elements $a,b,c\in F$ such that
$d:=b^2-4ac\neq0$. Let
\begin{equation}\label{Ldefeq}\renewcommand{\arraystretch}{1.2}
 L = \left\{
      \begin{array}{l@{\qquad\mbox{if }}l}
        F(\sqrt{d})&d \notin F^{\times2},\\
        F \oplus F&d \in F^{\times2}.
      \end{array}
    \right.
\end{equation}
In case $L = F \oplus F$, we consider $F$ diagonally embedded. If
$L$ is a field, we denote by $\bar x$ the Galois conjugate of
$x\in L$ over $F$. If $L=F\oplus F$, let $\overline{(x,y)}=(y,x)$.
In any case we let $N(x)=x\bar x$ and $\tr(x)=x+\bar x$.
\subsection{The unitary group}\label{unitarygroupsec}
We define the symplectic and unitary similitude groups by
\begin{align*}
 H(F) = \GSp_4(F) &:= \{g \in \GL_4(F) : \, ^tgJg = \mu(g)J,\:\mu(g)\in F^{\times} \}, \\
 G(F) = \GU(2,2;L) &:= \{g \in \GL_4(L) : \, ^t\bar{g}Jg = \mu(g)J,\:\mu(g)\in F^{\times}\},
\end{align*}
where $J = \mat{}{1_2}{-1_2}{}$. Note that $H(F) = G(F)
\cap\GL_4(F)$. As a minimal parabolic subgroup we choose the
subgroup of all matrices that become upper triangular after
switching the last two rows and last two columns. Let $P$ be the
standard maximal parabolic subgroup of $G(F)$ with a non-abelian
unipotent radical. Let $P = MN$ be the Levi decomposition of $P$.
We have $M = M^{(1)}M^{(2)}$, where
\begin{align}
 M^{(1)}(F)&=\{\begin{bmatrix}\zeta\\&1\\&&\bar
  \zeta^{-1}\\&&&1\end{bmatrix}:\:\zeta\in L^\times\} \label{M1defn},\\
 M^{(2)}(F)&=\{\begin{bmatrix}1\\&\alpha&&\beta\\&&\mu\\&\gamma&&\delta\end{bmatrix}
  \in G(F)\} \label{M2defn},\\
 N(F) &=\{\begin{bmatrix}
                 1 & z &  &  \\
                  & 1 &  &  \\
                  &  & 1 &  \\
                  &  & -\overline{z} & 1 \\
               \end{bmatrix}
               \begin{bmatrix}
                 1 &  & w & y \\
                  & 1 & \overline{y} &  \\
                  &  & 1 &  \\
                  &  &  & 1 \\
               \end{bmatrix}
              : w\in F,\;y,z \in L \}\label{Ndefn}.
\end{align}
For a matrix in $M^{(2)}(F)$ as the one above, the unitary
conditions are equivalent to $\mu=\bar\mu$ (i.e., $\mu\in
F^\times$), $\mu=\bar\alpha\delta-\beta\bar\gamma$,
$\bar\alpha\gamma=\bar\gamma\alpha$ and $\bar\delta\beta=\bar\beta\delta$.
In addition, we have $\bar\alpha\beta=\bar\beta\alpha$,
$\bar\delta\gamma=\bar\gamma\delta$, $\bar\alpha\delta=\bar\delta\alpha$,
$\bar\gamma\beta=\bar\beta\gamma$. Hence the following holds.

\begin{lemma}\label{M2structurelemma}
 Let
 $$
  \begin{bmatrix}1\\&\alpha&&\beta\\&&\mu\\&\gamma&&\delta\end{bmatrix}
 $$
 be an element of $M^{(2)}(F)$, as above. Then the quotient of any two entries
 of the matrix $\mat{\alpha}{\beta}{\gamma}{\delta}$, if defined, lies in $F$.
 Hence, if $\lambda$ is any invertible entry of $\mat{\alpha}{\beta}{\gamma}{\delta}$, then
 $$
  \mat{\alpha}{\beta}{\gamma}{\delta}=\lambda
  \underbrace{\mat{\alpha/\lambda\;}{\beta/\lambda}
  {\gamma/\lambda\;}{\delta/\lambda}}_{\in\GL_2(F)}.
 $$
 Consequently, the map
 \begin{align}\label{M2structurelemmaeq1}
  L^\times\times\GL_2(F)&\longrightarrow M^{(2)}(F),\\
  (\lambda,\mat{\alpha}{\beta}{\gamma}{\delta})&\longmapsto
  \begin{bmatrix}1\\&\lambda\alpha&&\lambda\beta\\&&N(\lambda)(\alpha\delta-\beta\gamma)\\
  &\lambda\gamma&&\lambda\delta\end{bmatrix},\nonumber
 \end{align}
 is surjective with kernel $\{(\lambda,\lambda^{-1}):\:\lambda\in F^\times\}$.
\end{lemma}
The modular factor of the parabolic $P$ is given by
\begin{equation}\label{deltaPformulaeq}
 \delta_P(\begin{bmatrix}\zeta\\&1\\&&\bar\zeta^{-1}\\&&&1\end{bmatrix}
 \begin{bmatrix}1\\&\alpha&&\beta\\&&\mu\\&\gamma&&\delta\end{bmatrix})
 =|N(\zeta)\mu^{-1}|^3\qquad(\mu=\bar\alpha\delta-\beta\bar\gamma),
\end{equation}
where $|\cdot|$ is the normalized absolute value on $F$.
\subsection{The Bessel subgroup}\label{besselsubgroupsec}
Recall that we fixed three elements $a,b,c\in F$ such that $d=b^2-4ac\neq0$. Let
$$
 S=\mat{a}{\frac b2}{\frac b2}{c},\qquad
 \xi=\mat{\frac b2}{c}{-a}{\frac{-b}2}.
$$
Then $F(\xi)=F+F\xi$ is a two-dimensional $F$-algebra isomorphic to $L$.
If $L=F(\sqrt{d})$ is a field, then an isomorphism is given by
$x+y\xi\mapsto x+y\frac{\sqrt{d}}2$. If $L=F\oplus F$, then an
isomorphism is given by $x+y\xi\mapsto(x+y\frac{\sqrt{d}}2,x-y\frac{\sqrt{d}}2)$.
The determinant map on $F(\xi)$ corresponds to the norm map on $L$. Let
$$
 T(F)=\{g\in\GL_2(F):\:^tgSg=\det(g)S\}.
$$
One can check that $T(F)=F(\xi)^\times$. Note that $T(F)\cong
L^\times$ via the isomorphism $F(\xi)\cong L$. We consider $T(F)$
a subgroup of $H(F)=\GSp_4(F)$ via
$$
 T(F)\ni g\longmapsto\mat{g}{}{}{\det(g)\,^tg^{-1}}\in H(F).
$$
Let
$$
 U(F)=\{\mat{1_2}{X}{}{1_2}\in\GSp_4(F):\:^tX=X\}
$$
and $R(F)=T(F)U(F)$. We call $R(F)$ the \emph{Bessel subgroup} of
$\GSp_4(F)$ (with respect to the given data $a,b,c$). Let $\psi$
be any non-trivial character $F\rightarrow\C^\times$. Let
$\theta:\:U(F)\rightarrow\C^\times$ be the character given by
\begin{equation}\label{thetadefeq}
 \theta(\mat{1}{X}{}{1})=\psi(\tr(SX)).
\end{equation}
Explicitly,
\begin{equation}\label{thetadef2eq}
 \theta(\begin{bmatrix}1&&x&y\\&1&y&z\\&&1\\&&&1\end{bmatrix})=\psi(ax+by+cz).
\end{equation}
We have $\theta(t^{-1}ut)=\theta(u)$ for all $u\in U(F)$ and $t\in T(F)$. Hence,
if $\Lambda$ is any character of $T(F)$, then the map
$tu\mapsto\Lambda(t)\theta(u)$ defines a character of $R(F)$. We denote
this character by $\Lambda\otimes\theta$.
\subsection{Parabolic induction from $P(F)$ to $G(F)$}\label{parabolicinductionsec}
Let $(\tau,V_\tau)$ be an irreducible, admissible representation
of $\GL_2(F)$, and let $\chi_0$ be a character of $L^\times$ such
that $\chi_0\big|_{F^\times}$ coincides with $\omega_{\tau}$, the
central character of $\tau$. Then the representation
$(\lambda,g)\mapsto\chi_0(\lambda)\tau(g)$ of
$L^\times\times\GL_2(F)$ factors through
$\{(\lambda,\lambda^{-1}):\:\lambda\in F^\times\}$, and
consequently, by Lemma \ref{M2structurelemma}, defines a
representation of $M^{(2)}(F)$ on the same space $V_\tau$. Let us
denote this representation by $\chi_0\times\tau$. Every
irreducible, admissible representation of $M^{(2)}(F)$ is of this
form. If $V_\tau$ is a space of functions on $\GL_2(F)$ on which
$\GL_2(F)$ acts by right translation, then $\chi_0\times\tau$ can
be realized as a space of functions on $M^{(2)}(F)$ on which
$M^{(2)}(F)$ acts by right translation. This is accomplished by
extending every $W\in V_\tau$ to a function on $M^{(2)}(F)$ via
\begin{equation}\label{extendedWformulaeq}
 W(\lambda g)=\chi_0(\lambda)W(g),\qquad\lambda\in L^\times,\:g\in\GL_2(F).
\end{equation}
If $V_\tau$ is the Whittaker model of $\tau$ with respect to the character $\psi$,
then the extended functions $W$ satisfy the transformation property
\begin{equation}\label{extendedWformulaWhittakereq}
 W(\begin{bmatrix}1\\&1&&x\\&&1\\&&&1\end{bmatrix}g)
 =\psi(x)W(g),\qquad x\in F,\:\:g\in M^{(2)}(F).
\end{equation}
If $s$ is a complex parameter,
$\chi$ is any character of $L^\times$, and $\chi_0\times\tau$ is
a representation of $M^{(2)}(F)$ as above, we denote by $I(s,\chi,\chi_0,\tau)$
the representation of $G(F)$ obtained by parabolic induction from the representation
of $P(F)=M(F)N(F)$ given on the Levi part by
$$
 \begin{bmatrix}\zeta\\&1\\&&\bar\zeta^{-1}\\&&&1\end{bmatrix}
 \begin{bmatrix}1\\&\lambda\alpha&&\lambda\beta\\
 &&N(\lambda)(\alpha\delta-\beta\gamma)\\&\lambda\gamma&&\lambda\delta\end{bmatrix}
 \longmapsto\big|N(\zeta\lambda^{-1})(\alpha\delta-\beta\gamma)^{-1}\big|^{3s}
 \chi(\zeta)\chi_0(\lambda)\tau(\mat{\alpha}{\beta}{\gamma}{\delta}).
$$
Explicitly, the space of $I(s,\chi,\chi_0,\tau)$ consists of functions
$f:\:G(F)\rightarrow V_\tau$ with the transformation property
\begin{align}\label{Isfctnspropeq}
 \nonumber &f(\begin{bmatrix}\zeta\\&1\\&&\bar\zeta^{-1}\\&&&1\end{bmatrix}
 \begin{bmatrix}1\\&\lambda\alpha&&\lambda\beta\\
 &&N(\lambda)(\alpha\delta-\beta\gamma)\\&\lambda\gamma&&\lambda\delta\end{bmatrix}g)\\
 &\hspace{20ex}=\big|N(\zeta\lambda^{-1})(\alpha\delta-\beta\gamma)^{-1}\big|^{3(s+\frac12)}
 \chi(\zeta)\chi_0(\lambda)\tau(\mat{\alpha}{\beta}{\gamma}{\delta})f(g).
\end{align}
Now assume that $V_\tau$ is the Whittaker model of $\tau$ with respect to the
character $\psi$ of $F$. If we associate to each $f$ as above the function
on $G(F)$ given by $W^\#_f(g)=f(g)(1)$, then we obtain another model of
$I(s,\chi,\chi_0,\tau)$ consisting of functions $W^\#:\:G(F)\rightarrow\C$.
These functions satisfy
\begin{equation}\label{Wsharpproperty1eq}
 W^\#(\begin{bmatrix}\zeta\\&1\\&&\bar\zeta^{-1}\\&&&1\end{bmatrix}
 \begin{bmatrix}1\\&\lambda\\&&N(\lambda)\\&&&\lambda\end{bmatrix}g)
 =|N(\zeta\lambda^{-1})|^{3(s+\frac12)}
 \chi(\zeta)\chi_0(\lambda)W^\#(g),\qquad \zeta,\lambda\in L^\times,
\end{equation}
and
\begin{equation}\label{Wsharpproperty2eq}
 W^\#(\begin{bmatrix}1 & z\\& 1\\&  & 1\\&  & -\overline{z} & 1 \\\end{bmatrix}
  \begin{bmatrix}1&&w&y\\& 1 & \overline{y} &x\\&  & 1\\&  &  & 1 \\\end{bmatrix}g)
 =\psi(x)W^\#(g),\qquad w,x\in F,\;y,z\in L.
\end{equation}
The following lemma gives a transformation property of $W^\#$
under the action of the elements of the Bessel subgroup $R(F)$.
\begin{lemma}\label{WsharpinvBessellemma}
 Let $(\tau,V_\tau)$ be a generic, irreducible, admissible representation of
 $\GL_2(F)$. We assume that $V_\tau$ is the Whittaker model of $\tau$ with
 respect to the non-trivial character $\psi^{-c}(x)=\psi(-cx)$ of $F$.
 Let $\chi$ and $\chi_0$ be characters of $L^\times$ such that
 $\chi_0\big|_{F^\times}=\omega_\tau$.
 Let $W^\#(\,\cdot\,,s):\:G(F)\rightarrow\C$ be a function in the above model
 of the induced representation $I(s,\chi,\chi_0,\tau)$, where $s$ is
 a complex parameter.
 Let $\theta$ be the character of $U(F)$ defined in (\ref{thetadefeq}).
 Let $\Lambda$ be the character of $L^\times\cong T(F)$ given by
 \begin{equation}\label{WsharpinvBessellemmaeq1}
  \Lambda(\zeta)=\chi(\bar\zeta)^{-1}\chi_0(\zeta)^{-1}.
 \end{equation}
 Let
\begin{equation}\label{etadefeq}
 \eta =\begin{bmatrix}1&0&&\\\alpha&1&&\\&&1&-\bar{\alpha}\\&&0&1\end{bmatrix},\qquad
 \text{where}\quad
 \renewcommand{\arraystretch}{2.1}\alpha:=\left\{\begin{array}{l@{\qquad\text{if }L}l}
 \displaystyle\frac{b+\sqrt{d}}{2c}&\text{ is a field},\\
 \displaystyle\Big(\frac{b+\sqrt{d}}{2c},\frac{b-\sqrt{d}}{2c}\Big)&=F\oplus F.
 \end{array}\right.
\end{equation}
 Then
 \begin{equation}\label{WsharpinvBessellemmaeq2}
  W^\#(\eta tuh,s)=\Lambda(t)^{-1}\theta(u)^{-1}W^\#(\eta h,s)
 \end{equation}
 for $t\in T(F)$, $u\in U(F)$ and $h\in G(F)$.
\end{lemma}
{\bf Proof.} If $L$ is a field, then the proof is word for word the same as on
p.\ 197/198 of \cite{Fu}. The case $L=F\oplus F$ requires the only modification
that the element $\zeta=x+\frac y2\sqrt{d}$ is to be replaced by
$\zeta=x+\frac y2(\sqrt{d},-\sqrt{d})$.\qed
\subsection{The local integral}\label{localintegralsec}
Let $(\pi,V_\pi)$ be an irreducible, admissible representation of
$H(F)=\GSp_4(F)$. Let the Bessel subgroup $R(F)$ be as defined  in
Section \ref{besselsubgroupsec}; it depends on the given data
$a,b,c\in F$. We assume that $V_\pi$ is a Bessel model for $\pi$
with respect to the character $\Lambda\otimes\theta$ of $R(F)$.
Hence, $V_\pi$ consists of functions $B:\:H(F)\rightarrow\C$
satisfying the Bessel transformation property
$$
 B(tuh)=\Lambda(t)\theta(u)B(h)\qquad
 \text{for }t\in T(F),\:u\in U(F),\:h\in H(F).
$$
Let $(\tau,V_\tau)$ be a generic, irreducible, admissible
representation of $\GL_2(F)$ such that $V_\tau$ is the
$\psi^{-c}$--Whittaker model of $\tau$ (we assume $c\neq0$). Let
$\chi_0$ be a character of $L^\times$ such that
$\chi_0\big|_{F^\times}=\omega_\tau$. Let $\chi$ be the character
of $L^\times$ for which (\ref{WsharpinvBessellemmaeq1}) holds. Let
$W^\#(\,\cdot\,,s)$ be an element of $I(s,\chi,\chi_0,\tau)$ for
which the restriction of $W^\#(\,\cdot\,,s)$ to the standard
maximal compact subgroup of $G(F)$ (see below for more details) is
independent of $s$, i.e., $W^\#(\,\cdot\,,s)$ is a ``section'' of
the family of induced representations $I(s,\chi,\chi_0,\tau)$. By
Lemma \ref{WsharpinvBessellemma} it is meaningful to consider the
integral
\begin{equation}\label{localZseq}
 Z(s)=\int\limits_{R(F)\backslash H(F)}W^\#(\eta h,s)B(h)\,dh.
\end{equation}
In the following we shall compute these integrals for certain
choices of $W^\#$ and $B$. We shall only consider $\GSp_4(F)$
representations $\pi$ that are relevant for the global application
to Siegel modular forms we have in mind. In the real case we shall
assume that $\pi$ is a holomorphic discrete series representation
and that $B$ corresponds to the lowest weight vector. In the
$p$-adic case we shall assume that $\pi$ is an unramified
representation and that $B$ corresponds to the spherical vector.

The generic $\GL_2(F)$ representation $\tau$, however, will be
only mildly restricted in the real case, and, in the $p$-adic
case, will be a Steinberg representation twisted by an unramified
character. In the real case, the function $W^\#$ will be
constructed from a certain vector of the ``correct'' weight in
$V_\tau$. In the $p$-adic case, the function $W^\#$ will be
constructed from the local newform in $V_\tau$.

In each case our calculations will show that the integral $(\ref{localZseq})$
converges absolutely for ${\rm Re}(s)$ large enough and has meromorphic
continuation to all of $\C$. Our choice of $W^\#$ will be such that $Z(s)$
is closely related to the local $L$-factor $L(s,\pi\times\tau)$.

Note that the integral (\ref{localZseq}) has been calculated in
\cite{Fu} for $\pi$ and $\tau$ both holomorphic discrete series
representations with related lowest weights in the real case and
$\pi$ and $\tau$ both unramified representations in the $p$-adic
case.
\section{Local non-archimedean ramified theory}\label{non-arch-ram-sec}
In this section, we evaluate (\ref{localZseq}) in the
non-archimedean setting. The key steps are the choices of the
vector $W^\#$ and the actual computation of the integral $Z(s)$.
\subsection{Setup}
Let $F$ be a non-archimedean local field of characteristic zero.
Let $\OF$, $\p$, $\varpi$, $q$ be the ring of integers, prime ideal, uniformizer
and cardinality of the residue class field $\OF/\p$, respectively. Recall that we fix three
elements $a,b,c\in F$ such that $d := b^2-4ac \neq 0$. Let $L$ be as in (\ref{Ldefeq}).
We shall make the following {\bf assumptions}:
\begin{description}
 \item[(A1)] $a,b\in\OF$ and $c\in\OF^\times$.
 \item[(A2)] If $d \not\in F^{\times2}$, then $d$ is the generator of the discriminant of $L/F$.
  If $d \in F^{\times2}$, then $d \in \OF^{\times}$.
\end{description}

{\bf Remark:} In \cite[p.\ 198]{Fu}, Furusawa makes a stronger
assumption on $a,b,c$, namely, $\mat{a}{b/2}{b/2}{c} \in
M_2(\OF)$. However, it is necessary to make the weaker
assumption $a,b,c \in \OF$ for the global integral calculation
$(4.5)$ in \cite[p.\ 210]{Fu} to be valid for $D \equiv 3
\pmod{4}$. (This is because the matrix $S(-D)$ on p.\ 208 is not in
$M_2(\OF_2)$ for $D \equiv 3 \pmod{4}$.) One can check that the
non-archimedean unramified calculation in \cite{Fu} is valid with
the weaker assumption $a,b,c \in \OF$. Hence, the global result of
\cite{Fu} is still valid but the assumptions (A1) and (A2)
above are the correct ones.

We set the Legendre symbol as follows,
\begin{equation}\label{legendresymboldefeq}
 \Big(\frac L{\p}\Big) := \left\{
                  \begin{array}{l@{\qquad\text{if }}l@{\qquad}l}
                    -1, &d \not\in F^{\times2},\:d \not\in \p&\mbox{(the inert case)}, \\
                    0, &d\not\in F^{\times2},\:d\in\p&\mbox{(the ramified case)}, \\
                    1, &d \in F^{\times2}&\mbox{(the split case)}.
                  \end{array}\right.
\end{equation}
If $L$ is a field, then let $\OF_L$ be its ring of integers.
If $L = F \oplus F$, then let $\OF_L = \OF \oplus \OF$.
Note that $x\in \OF_L$ if and only if $N(x), \tr(x) \in \OF$.
If $L$ is a field then we have $x \in \OF_L^{\times}$ if and only if
$N(x) \in \OF^{\times}$. If $L$ is not a field then
$x \in \OF_L, N(x) \in \OF^{\times}$ implies that
$x \in \OF_L^{\times} = \OF^{\times} \oplus \OF^{\times}$.
Let $\varpi_L$ be the uniformizer of $\OF_L$ if $L$ is a field and set
$\varpi_L = (\varpi,1)$ if $L$ is not a field.
Note that, if $(\frac{L}{\p}) \neq -1$, then $N(\varpi_L) \in \varpi \OF^{\times}$. Let
\begin{equation}\label{xi0defeq}\renewcommand{\arraystretch}{2.1}
 \xi_0:=\left\{\begin{array}{l@{\qquad\text{if }L}l}
 \displaystyle\frac{-b+\sqrt{d}}2&\text{ is a field},\\
 \displaystyle\Big(\frac{-b+\sqrt{d}}2,\frac{-b-\sqrt{d}}2\Big)&=F\oplus F.
 \end{array}\right.
\end{equation}
and
\begin{equation}\label{alphadefeq}\renewcommand{\arraystretch}{2.1}
 \alpha:=\left\{\begin{array}{l@{\qquad\text{if }L}l}
 \displaystyle\frac{b+\sqrt{d}}{2c}&\text{ is a field},\\
 \displaystyle\Big(\frac{b+\sqrt{d}}{2c},\frac{b-\sqrt{d}}{2c}\Big)&=F\oplus F.
 \end{array}\right.
\end{equation}
We fix the following ideal in $\OF_L$,
\begin{equation}\label{ideal defn}\renewcommand{\arraystretch}{1.3}
 \P := \p\OF_L = \left\{
                  \begin{array}{l@{\qquad\text{if }}l}
                    \p_L & \big(\frac L{\p}\big) = -1,\\
                    \p_L^2 & \big(\frac L{\p}\big) = 0,\\
                    \p \oplus \p & \big(\frac L{\p}\big) = 1.
                  \end{array}
                \right.
\end{equation}
Here, $\p_L$ is the maximal ideal of $\OF_L$ when $L$ is a field
extension. Note that $\P$ is prime only if $\big(\frac
L\p\big)=-1$. We have $\P^n\cap\OF=\p^n$ for all $n\geq0$. We now
state a number-theoretic lemma which will be crucial in Section
\ref{double coset decomposition sec}.

\begin{lemma}\label{quadraticextdisclemma}
 Let notations be as above.
 \begin{enumerate}
  \item The elements $1$ and $\xi_0$ constitute an integral basis of $L/F$ (i.e.,
   a basis of the free $\OF$-module $\OF_L$). The
   elements $1$ and $\alpha$ also constitute an integral basis of $L/F$.
  \item There exists no $x\in\OF$ such that $\alpha+x\in\P$.
 \end{enumerate}
\end{lemma}
{\bf Proof.} i) Since $c\in\OF^\times$ and $b\in\OF$, the second
assertion of i) follows from the first one. To prove the first
assertion, first note that $\xi_0$ satisfies
$\xi_0^2+\xi_0b+ac=0$, and therefore belongs to $\OF_L$. Since the
claim is easily verified if $L=F\oplus F$, we will assume that $L$
is a field. Let $A,B\in F$ be such that $1$ and
$\xi_1:=A+B\sqrt{d}$ is an integral basis of $L/F$. Then
$$
 \det(\mat{1}{\xi_1}{1}{\bar\xi_1})^2=4B^2d
$$
generates the discriminant of $L/F$. Since $d$ also generates the
discriminant by assumption $(A2)$, it follows that
$2B\in\OF_F^\times$. Dividing $\xi_1$ by this unit, we may assume
$\xi_1=A+\frac12\sqrt{d}$ for some $A\in F$. Now let us represent
$\xi_0$ in this integral basis,
$$
 \xi_0=x+y\xi_1,\qquad x,y\in\OF_F,
$$
i.e.,
$$
 \frac{-b+\sqrt{d}}2=x+y(A+\frac12\sqrt{d}).
$$
Comparing coefficients, we get $y=1$ and $A=-\frac b2-x$. We may modify $\xi_1$
by adding the integral element $x$ and still obtain an integral basis. But
$\xi_1+x=\xi_0$, and the assertion follows.
\nll
ii) Let $X\subset\OF_L/\P$ be the image of the injection
$$
 \OF/\p\longrightarrow\OF_L/\P.
$$
Note that the field on the left hand side has $q$ elements, and the ring on the right hand
side has $q^2$ elements, for any value of $\big(\frac L\p\big)$. Our claim is equivalent
to the statement that $\bar\alpha$, the image of $\alpha$ in $\OF_L/\P$, does
not lie in the subring $X$ of $\OF_L/\P$.
Assume that $\bar\alpha\in X$. By i), any element $z\in\OF_L$ can be
(uniquely) written as
$$
 z=x\alpha+y, \qquad x,y\in\OF.
$$
Applying the projection to $\OF_L/\P$, it follows that $\bar
z=\bar x\bar\alpha+\bar y\in X$. This is a contradiction, since
$\bar z$ runs through all elements of $\OF_L/\P$, but $X$ is a
proper subset.\qed

\vspace{2ex}
Note that, via the identification
$T(F)=L^\times$ described in Sect.\ \ref{besselsubgroupsec}, the
element $\xi_0$ corresponds to the matrix $\mat{0}{c}{-a}{-b}$.
Therefore, by Lemma \ref{quadraticextdisclemma} i),
\begin{equation}\label{OLexpliciteq}
 \OF_L=\OF\oplus\OF\xi_0=\{\mat{x}{yc}{-ya}{\:x-yb}:\:x,y\in\OF\}.
\end{equation}
Since $c$ is assumed to be a unit, it follows that
\begin{equation}\label{OLGL2OFeq}
 \OF_L=T(F)\cap M_2(\OF)\qquad\text{and}\qquad
 \OF_L^\times=T(F)\cap\GL_2(\OF).
\end{equation}
\subsection{The spherical Bessel function}\label{sphericalbesselfunctionsec}
Let $(\pi,V_\pi)$ be an unramified, irreducible, admissible
representation of $\GSp_4(F)$. Then $\pi$ can be realized as the
unramified constituent of an induced representation of the form
$\chi_1\times\chi_2\rtimes\sigma$, where $\chi_1$, $\chi_2$ and
$\sigma$ are unramified characters of $F^\times$; here, we used
the notation of \cite{ST} for parabolic induction. Let
$$
 \gamma^{(1)}=\chi_1\chi_2\sigma,\qquad\gamma^{(2)}=\chi_1\sigma,\qquad
 \gamma^{(3)}=\sigma,\qquad\gamma^{(4)}=\chi_2\sigma.
$$
Then $\gamma^{(1)}\gamma^{(3)}=\gamma^{(2)}\gamma^{(4)}$ is the central character
of $\pi$. The numbers $\gamma^{(1)}(\varpi),\ldots,\gamma^{(4)}(\varpi)$ are the
Satake parameters of $\pi$. The degree-$4$ $L$-factor of $\pi$ is given by
$\prod_{i=1}^4(1-\gamma^{(i)}(\varpi)q^{-s})^{-1}$.

Let $\Lambda$ be any character of $T(F)\cong L^\times$. We assume
that $V_\pi$ is the Bessel model with respect to the character
$\Lambda\otimes\theta$ of $R(F)$; see Sect.\
\ref{besselsubgroupsec}. Let $B\in V_\pi$ be a spherical vector.
By \cite{Su}, Proposition 2-5, we have $B(1)\neq0$. It follows
from $B(1)\neq0$ and (\ref{OLGL2OFeq}) that necessarily
$\Lambda\big|_{\OF_L^\times}=1$. For $l,m\in\Z$ let
\begin{equation}\label{hlmdefeq}
 h(l,m)=\begin{bmatrix}\varpi^{2m+l}\\&\varpi^{m+l}\\&&1\\&&&\varpi^m\end{bmatrix}.
\end{equation}
Then, as in (3.4.2) of \cite{Fu},
\begin{equation}\label{RFKHrepresentativeseq}
 H(F)=\bigsqcup_{l\in\Z}\bigsqcup_{m\geq0}R(F)h(l,m)K^H,\qquad
 K^H=\GSp_4(\OF).
\end{equation}
The double cosets on the right hand side are pairwise disjoint.
Since $B$ transforms on the left under $R(F)$ by the character
$\Lambda \otimes \theta$ and is right $K^H$-invariant, it follows
that $B$ is determined by the values $B(h(l,m))$. By Lemma (3.4.4)
of \cite{Fu} we have $B(h(l,m))=0$ for $l<0$, so that $B$ is
determined by the values $B(h(l,m))$ for $l,m\geq0$.

In \cite{Su}, 2-4, Sugano has given a formula for $B(h(l,m))$ in
terms of a generating function. It turns out that for our purposes
we only require the values $B(h(l,0))$. In this special case
Sugano's formula reads
\begin{equation}\label{suganox0eq}
 \sum_{l\geq0}B(h(l,0))y^l= \frac{1-A_5y-A_2A_4y^2}{Q(y)},
\end{equation}
where
\begin{align}\label{suganoQeq}
 Q(y)&=\prod_{i=1}^4\big(1-\gamma^{(i)}(\varpi_F)q^{-3/2}y\big),
\end{align}
and where $A_2,A_4,A_5$ are given in the following table. Set $H(y) = 1-A_5y-A_2A_4y^2$.
$$\renewcommand{\arraystretch}{1.4}
 \begin{array}{|c|c|c|c|}
 \hline
 &\big(\frac L\p\big)=-1&\big(\frac L\p\big)=0&\big(\frac L\p\big)=1\\\hline \hline
 A_2&q^{-2}\Lambda(\varpi_F)&q^{-2}\Lambda(\varpi_F)&q^{-2}\Lambda(\varpi_F)\\\hline
 A_4&q^{-2}&0&-q^{-2}\\\hline
 A_5&0&q^{-2}\Lambda(\varpi_L)
  &q^{-2}\big(\Lambda(\varpi_L)+\Lambda(\varpi_F\varpi_L^{-1})\big)\\\hline
  H(y)&1-q^{-4}\Lambda(\varpi_F)y^2&1-q^{-2}\Lambda(\varpi_L)y&
   \begin{minipage}{36ex}
   $\begin{array}{c}1-q^{-2}\big(\Lambda(\varpi_L)+\Lambda(\varpi_F\varpi_L^{-1})\big)y\\
   +q^{-4}\Lambda(\varpi_F)y^2\end{array}$\end{minipage}\\\hline
 \end{array}
$$
\subsection{The local compact subgroup}
We define congruence subgroups of $\GL_2(F)$, as follows. For
$n=0$ let $K^{(1)}(\p^0)=\GL_2(\OF)$. For $n>0$ let
\begin{equation}\label{K1defeq}
 K^{(1)}(\p^n)=\GL_2(F)\cap\mat{\OF^\times}{\OF}{\p^n}{1+\p^n}.
\end{equation}
The following result is well known (see \cite{Cas}, \cite{De}).
\begin{theorem}\label{GL2newformtheorem}
 Let $(\tau,V)$ be a generic, irreducible, admissible representation of $\GL_2(F)$.
 Then the spaces
 $$
  V(n)=\{v\in V:\:\tau(g)v=v\text{ for all }g\in K^{(1)}(\p^n)\}
 $$
 are non-zero for $n$ large enough. If $n$ is minimal with $V(n)\neq0$, then
 $\dim(V(n))=1$.
\end{theorem}
If $n$ is minimal such that $V(n)\neq0$, then $\p^n$ is called the
\emph{conductor} of $\tau$.
In this section we shall define a family $K^\#(\P^n)$, $n\geq0$, of compact-open subgroups
of $G(F)$, the relevance of which is as follows.
Recall that our goal is to evaluate integrals of the form
\begin{equation}\label{Zs2eq}
 Z(s)=\int\limits_{R(F)\backslash H(F)}W^\#(\eta h,s)B(h)\,dh,
\end{equation}
where $W^\#(\,\cdot\,,s)$ is a section in a family of induced
representations $I(s,\chi,\chi_0,\tau)$. The choice of the
function $W^\#(\,\cdot\,,s)$ is crucial for our purposes. We will
define it in such a way that $W^\#(\,\cdot\,,s)$ is supported on
$M(F)N(F)K^\#(\P^n)$, where $\p^n$ is the conductor of the
$\GL_2(F)$ representation $\tau$.

Recall that $\P = \p\OF_L$. For $n=0$ we let
$K^\#(\P^0)=G(F)\cap\GL_4(\OF_L)$ be the standard maximal compact subgroup of $G(F)$.
For $n>0$ we define
\begin{equation}\label{K21defeq}
 K^\#_1(\P^n):=\{g\in G(F):\:\mu(g)=1\}\cap
 \begin{bmatrix}1+\P^n&\P^n&\OF_L&\OF_L\\\OF_L&1+\P^n&\OF_L&\OF_L\\
 \P^n&\P^n&1+\P^n&\OF_L\\\P^n&\P^n&\P^n&1+\P^n\end{bmatrix}
\end{equation}
and
\begin{equation}\label{K2defeq}
 K^\#(\P^n)=K^\#_1(\P^n)\rtimes \{{\rm diag}(1,\mu,\mu,1):\:\mu\in\OF^\times\}.
\end{equation}
The $\GL_2$ congruence subgroup $K^{(1)}(\p^n)$ defined above can
be embedded into $K^\#(\P^n)$ in the following way,
\begin{equation}\label{K1K2embeddingeq}
 \mat{\alpha}{\beta}{\gamma}{\delta}\longmapsto
 \begin{bmatrix}1\\&\alpha/\mu&&\beta\\&&1\\&\gamma/\mu&&\delta\end{bmatrix}
 \begin{bmatrix}1\\&\mu&\\&&\mu\\&&&1\end{bmatrix},\qquad
 \text{where }\mu=\alpha\delta-\beta\gamma.
\end{equation}
Important for us will be the intersection
\begin{align}\label{K2capHeq}
 K^\#(\p^n)&:=H(F)\cap K^\#(\P^n)\nonumber\\
 &=\begin{bmatrix}1+\p^n&\p^n&\OF&\OF\\\OF&1+\p^n&\OF&\OF\\
 \p^n&\p^n&1+\p^n&\OF\\\p^n&\p^n&\p^n&1+\p^n\end{bmatrix}
 \begin{bmatrix}1\\&\mu&\\&&\mu\\&&&1\end{bmatrix},
 \quad\mu\in\OF^\times.
\end{align}
Note that $H(F)\cap K^\#(\P^n)=K^H\cap K^\#(\P^n)$, where $K^H =
\GSp_4(\OF)$ is the maximal compact subgroup of $H(F)$. It follows
from Lemma \ref{M2structurelemma} that the map
\begin{align}\label{M2structurelemmaeq2}
  \OF_L^\times\times\GL_2(\OF)&\longrightarrow M^{(2)}(F)\cap\GL_4(\OF_L),\\
  (\lambda,\mat{\alpha}{\beta}{\gamma}{\delta})&\longmapsto
  \begin{bmatrix}1\\&\lambda\alpha&&\lambda\beta\\&&N(c)(\alpha\delta-\beta\gamma)\\
  &\lambda\gamma&&\lambda\delta\end{bmatrix},\nonumber
\end{align}
is surjective with kernel $\{(\lambda,\lambda^{-1}):\:\lambda\in\OF_F^\times\}$.
For $n>0$ this map induces a surjection
\begin{equation}\label{M2structurelemmaeq3}
  (1+\P^n)\times K^{(1)}(\p^n)\longrightarrow M^{(2)}(F)\cap K^\#(\P^n)
\end{equation}
with kernel $\{(\lambda,\lambda^{-1}):\:\lambda\in1+\P^n\}$.
\subsection{The function $W^\#$}\label{W sharp function}
We shall now define the specific function $W^\#(\,\cdot\,,s)$ for
which we shall evaluate the integral (\ref{Zs2eq}). Let
$(\tau,V_\tau)$ be a generic, irreducible, admissible
representation of $\GL_2(F)$. We assume that $V_\tau$ is the
Whittaker model of $\tau$ with respect to the character of $F$
given by $\psi^{-c}(x)=\psi(-cx)$. Let $\p^n$ be the conductor of
$\tau$. Let $W^{(0)}\in V(n)$ be the local newform, i.e., the
essentially unique non-zero $K^{(1)}(\p^n)$ invariant vector in
$V_\tau$. We can make it unique by requiring that $W^{(0)}(1)=1$.

We choose any character $\chi_0$ of $L^\times$ such that
\begin{equation}\label{chi0defpropertieseq}
 \chi_0\big|_{F^\times}=\omega_{\tau}\qquad\text{and}\qquad
 \chi_0\big|_{1+\P^n}=1.
\end{equation}
This can be accomplished by extending $\omega_\tau\big|_{\OF^\times}$ to
$\OF_L^\times$ using the injection $\OF^\times/(1+\p^n)\hookrightarrow
\OF_L^\times/(1+\P^n)$, and defining $\chi_0$ suitably on prime elements.
We extend $W^{(0)}$ to a function on $M^{(2)}(F)$ via
\begin{equation}\label{W1extensioneq}
 W^{(0)}(ag)=\chi_0(a)W^{(0)}(g),\qquad a\in L^\times,\:g\in\GL_2(F)
\end{equation}
(see (\ref{M2structurelemmaeq1})). It follows from (\ref{M2structurelemmaeq3}) that
\begin{equation}\label{W1extensioninvarianceeq}
 W^{(0)}(g\kappa)=W^{(0)}(g),\qquad\text{for }g\in M^{(2)}(F)\text{ and }
 \kappa\in M^{(2)}(F)\cap K^\#(\P^n).
\end{equation}
As in Sect.\ \ref{sphericalbesselfunctionsec}, let $(\pi,V_\pi)$
be an unramified, irreducible, admissible representation of
$\GSp_4(F)$, where $V_\pi$ is the Bessel model for $\pi$ with
respect to the character $\Lambda\otimes\theta$ of
$R(F)=T(F)U(F)$. As was pointed out in Sect.\
\ref{sphericalbesselfunctionsec}, the character $\Lambda$ is
necessarily unramified. Let $\chi$ be the character of $L^\times$
given by
\begin{equation}\label{chidefeq}
 \chi(\zeta)=\Lambda(\bar \zeta)^{-1}\chi_0(\bar\zeta)^{-1},
\end{equation}
so that (\ref{WsharpinvBessellemmaeq1}) holds.

Given a complex number $s$, there exists a unique function
$W^\#(\,\cdot\,,s):\:G(F)\rightarrow\C$ with the following properties.
\begin{enumerate}
 \item If $g\notin M(F)N(F)K^\#(\P^n)$, then $W^\#(g,s)=0$.
 \item If $g=mnk$ with $m\in M(F)$, $n\in N(F)$, $k\in K^\#(\P^n)$, then
  $W^\#(g,s)=W^\#(m,s)$.
 \item For $\zeta\in L^\times$ and $\mat{\alpha}{\beta}{\gamma}{\delta}\in M^{(2)}(F)$,
  \begin{equation}\label{Wsharpformulaeq}
   W^\#(\begin{bmatrix}\zeta\\&1\\&&\bar{\zeta}^{-1}\\&&&1\end{bmatrix}
   \begin{bmatrix}1\\&\alpha&&\beta\\&&\mu\\&\gamma&&\delta\end{bmatrix},s)
   =|N(\zeta)\cdot\mu^{-1}|^{3(s+1/2)}\chi(\zeta)\,
   W^{(0)}(\mat{\alpha}{\beta}{\gamma}{\delta}).
  \end{equation}
  Here $\mu=\bar\alpha\delta-\beta\bar\gamma$.
\end{enumerate}
To verify that such a function exists, use (\ref{W1extensioninvarianceeq}) and
$$
 (M(F)N(F))\cap K^\#(\P^n)=\big(M(F)\cap K^\#(\P^n)\big)
 \big(N(F)\cap K^\#(\P^n)\big).
$$
Also, one has to use the fact that $\chi\big|_{1+\P^n}=1$. Note that
$W^\#(\,\cdot\,,s)$ is an element of the induced representation
$I(s,\chi,\chi_0,\tau)$ discussed in Section \ref{parabolicinductionsec}.
In particular, Lemma \ref{WsharpinvBessellemma} applies. Note that if $n=0$,
i.e., if $\tau$ is unramified, then $W^\#(\,\cdot\,,s)$ coincides with the
function $W_v(\,\cdot\,,s)$ defined on p.\ 200 of \cite{Fu}.
\subsection{Basic local integral computation}\label{basic integral}
Let $W^\#(\,\cdot\,,s)$ be the element of $I(s,\chi,\chi_0,\tau)$
defined in the previous section. Let $B$ be the spherical vector
in the $\Lambda\otimes\theta$ Bessel model of the unramified
representation $\pi$ of $\GSp_4(F)$, as in Sect.\ \ref{sphericalbesselfunctionsec}.
We shall compute the integral
\begin{equation}\label{Zvdefeq}
 Z(s)=\int\limits_{R(F)\backslash H(F)}W^\#(\eta h,s)B(h)\,dh.
\end{equation}
By Lemma \ref{WsharpinvBessellemma}, the integral (\ref{Zvdefeq})
is well-defined. By (\ref{RFKHrepresentativeseq}) and the fact
that $B(h(l,m))=0$ for $l<0$ (\cite{Fu} Lemma 3.4.4), we have
\begin{align}\label{Zvcalceq1}
 Z(s)&=\sum_{l,m\geq0}\;\int\limits_{R(F)\backslash R(F)h(l,m)K^H}
  W^\#(\eta h,s)B(h)\,dh\nonumber\\
 &=\sum_{l,m\geq0}\;\int\limits_{h(l,m)^{-1}R(F)h(l,m)\cap K^H\backslash K^H}
  W^\#(\eta h(l,m)h,s)B(h(l,m)h)\,dh\nonumber\\
 &=\sum_{l,m\geq0}\;B(h(l,m))\int\limits_{h(l,m)^{-1}R(F)h(l,m)\cap K^H\backslash K^H}
  W^\#(\eta h(l,m)h,s)\,dh.
\end{align}
The function $W^\#$ is only invariant under $K^\#(\P^n)$. Since
our integral (\ref{Zvcalceq1}) is over elements of $H(F)$, all
that is relevant is that $W^\#$ is invariant under the group
$K^\#(\p^n)$ defined in (\ref{K2capHeq}). Let us abbreviate
$K_{l,m}:=h(l,m)^{-1}R(F)h(l,m)\cap K^H$. Suppose we had a system
of representatives $\{s_i\}$ for the double coset space
$K_{l,m}\backslash K^H/K^\#(\p^n)$ (it will depend on $l$ and $m$,
of course). Then, from (\ref{Zvcalceq1}),
\begin{align}\label{Zvcalceq3}
 Z(s)&=\sum_{l,m\geq0}\;\sum_i\;B(h(l,m))
  \int\limits_{K_{l,m}\backslash K_{l,m}s_iK^\#(\p^n)}
  W^\#(\eta h(l,m)h,s)\,dh\nonumber\\
 &=\sum_{l,m\geq0}\;\sum_i\;B(h(l,m))\,W^\#(\eta h(l,m)s_i,s)
  \int\limits_{K_{l,m}\backslash K_{l,m}s_iK^\#(\p^n)}\,dh.
\end{align}
In practice it will be difficult to obtain the system $\{s_i\}$.
However, we can save some work by exploiting the fact that $W^\#$
is supported on the small subset $M(F)N(F)K^\#(\P^n)$ of $G(F)$.
Hence, we shall proceed as follows.
\begin{description}
 \item[Step $1$:] First we determine a \emph{preliminary} decomposition
  \begin{equation}\label{prelimdecompeq}
   K^H=\bigcup_jK_{l,m}\,s'_j\,K^\#(\p^n),
  \end{equation}
  which is not necessarily disjoint. We may assume that the $s'_j$ are
  taken from the system of representatives for $K^H/K^\#(\p^n)$ to be
  determined in the next section (but some of these will be absorbed
  in $K_{l,m}$, so that we get an initial reduction).
 \item[Step $2$:] Then we consider the values $ W^\#(\eta h(l,m)s_j',s)$.
  If $\eta h(l,m)s'_j\notin M(F)N(F)K^\#(\P^n)$, then $s'_j$ makes no contribution to the
  integral (\ref{Zvcalceq1}). Therefore, all that is relevant is the
  subset $\{s''_j\}\subset\{s'_j\}$ of representatives for which
  $\eta h(l,m)s''_j\in M(F)N(F)K^\#(\P^n)$. Hence we consider the set
  $$
   S:=\bigcup_j K_{l,m}s''_jK^\#(\p^n).
  $$
 \item[Step $3$:] Now, from this much smaller set of representatives $\{s''_j\}$ we
  determine a subset $\{s'''_j\}$ such that this union becomes disjoint:
  $$
   S=\bigsqcup_j K_{l,m}s'''_jK^\#(\p^n).
  $$
\end{description}
The integral (\ref{Zvcalceq1}) is then given by
\begin{align}\label{Zvcalceq4}
 Z(s)&=\sum_{l,m\geq0}\;\sum_j\;B(h(l,m))\,W^\#(\eta h(l,m)s_j''',s)
  \int\limits_{K_{l,m}\backslash K_{l,m}s_j'''K^\#(\p^n)}\,dh.
\end{align}
Finally, we have to compute the volumes, evaluate $W^\#$, and
carry out the summations with the help of Sugano's formula
(\ref{suganox0eq}).
\subsection{The cosets $K^\#(\p^0)/K^\#(\p)$}\label{Ksharpcosetssec}
From this point on we will assume that the conductor $\p^n$ of the given $\GL(2)$
representation $\tau$ satisfies \underline{$n=1$}. We need to
determine representatives for the coset space
\begin{equation}\label{newK2KHcosetseq1}
 K^\#(\p^0)/K^\#(\p),\qquad\text{where }K^\#(\p^0)=K^H=\GSp_4(\OF).
\end{equation}
Note that this coset space is isomorphic to
$K^\#_1(\p^0)/K^\#_1(\p)$, where $K^\#_1(\p)=K^\#(\p)\cap\{g\in H(F):\mu(g)=1\}$.
Let
\begin{equation}\label{Weylgroupelementsdefeq}
 \qquad s_1=\begin{bmatrix}&1\\1\\&&&1\\&&1\end{bmatrix},
 \qquad s_2=\begin{bmatrix}&&1\\&1\\-1\\&&&1\end{bmatrix}.
\end{equation}
It follows from the Bruhat decomposition for $\SSp_4(\OF/\p)$ that
\begin{align}
 \label{newK2KHcosetseq41}&K^\#(\p^0)=\bigsqcup_{a_1,a_2\in\OF^\times/(1+\p)}
  \begin{bmatrix}a_1\\&a_2\\&&a_1^{-1}\\&&&a_2^{-1}\end{bmatrix}K^\#(\p)\\
 \label{newK2KHcosetseq42}
 &\sqcup\bigsqcup_{a_1,a_2\in\OF^\times/(1+\p)}\:\bigsqcup_{x\in\OF/\p}
  \begin{bmatrix}a_1\\&a_2\\&&a_1^{-1}\\&&&a_2^{-1}\end{bmatrix}
  \begin{bmatrix}1\\x&1\\&&1&-x\\&&&1\end{bmatrix}s_1K^\#(\p)\\
 \label{newK2KHcosetseq43}
 &\sqcup\bigsqcup_{a_1,a_2\in\OF^\times/(1+\p)}\:\bigsqcup_{x\in\OF/\p}
  \begin{bmatrix}a_1\\&a_2\\&&a_1^{-1}\\&&&a_2^{-1}\end{bmatrix}
  \begin{bmatrix}1&&x\\&1\\&&1&\\&&&1\end{bmatrix}s_2K^\#(\p)\\
 \label{newK2KHcosetseq44}
 &\sqcup\bigsqcup_{a_1,a_2\in\OF^\times/(1+\p)}\:\bigsqcup_{x,y\in\OF/\p}
  \begin{bmatrix}a_1\\&a_2\\&&a_1^{-1}\\&&&a_2^{-1}\end{bmatrix}
  \begin{bmatrix}1\\x&1&&y\\&&1&-x\\&&&1\end{bmatrix}s_1s_2K^\#(\p)\\
 \label{newK2KHcosetseq45}
 &\sqcup\bigsqcup_{a_1,a_2\in\OF^\times/(1+\p)}\:\bigsqcup_{x,y\in\OF/\p}
  \begin{bmatrix}a_1\\&a_2\\&&a_1^{-1}\\&&&a_2^{-1}\end{bmatrix}
  \begin{bmatrix}1&&x&y\\&1&y\\&&1\\&&&1\end{bmatrix}s_2s_1K^\#(\p)\\
 \label{newK2KHcosetseq46}
 &\sqcup\bigsqcup_{a_1,a_2\in\OF^\times/(1+\p)}\:\bigsqcup_{x,y,z\in\OF/\p}
  \begin{bmatrix}a_1\\&a_2\\&&a_1^{-1}\\&&&a_2^{-1}\end{bmatrix}
  \begin{bmatrix}1&&&y\\x&1&y&xy+z\\&&1&-x\\&&&1\end{bmatrix}s_1s_2s_1K^\#(\p)\\
 \label{newK2KHcosetseq47}
 &\sqcup\bigsqcup_{a_1,a_2\in\OF^\times/(1+\p)}\:\bigsqcup_{x,y,z\in\OF/\p}
  \begin{bmatrix}a_1\\&a_2\\&&a_1^{-1}\\&&&a_2^{-1}\end{bmatrix}
  \begin{bmatrix}1&&x&y\\&1&y&z\\&&1\\&&&1\end{bmatrix}s_2s_1s_2K^\#(\p)\\
 \label{newK2KHcosetseq48}
 &\sqcup\bigsqcup_{a_1,a_2\in\OF^\times/(1+\p)}\:\bigsqcup_{w,x,y,z\in\OF/\p}
  \begin{bmatrix}a_1\\&a_2\\&&a_1^{-1}\\&&&a_2^{-1}\end{bmatrix}
  \begin{bmatrix}1&&x&y\\w&1&wx+y&wy+z\\&&1&-w\\&&&1\end{bmatrix}s_1s_2s_1s_2K^\#(\p).
\end{align}
\subsection{Double coset decomposition}\label{double coset
decomposition sec} Recall that we are interested in the double
cosets $K_{l,m}\backslash K^H/K^\#(\p)$, where
$K_{l,m}=h(l,m)^{-1}R(F)h(l,m)\cap K^H$.
\subsubsection{Step $1$: Preliminary decomposition}
Observe that $K_{l,m}$ contains all elements
$$
 \begin{bmatrix}u\\&u\\&&u\\&&&u\end{bmatrix},\qquad u\in\OF^\times,
 \qquad\text{and}\qquad
 \begin{bmatrix}1&&\OF&\OF\\&1&\OF&\OF\\&&1\\&&&1\end{bmatrix},
$$
and that $K^\#(\p)$ contains all elements of the form
${\rm diag}(1,\mu,\mu,1)$, $\mu\in\OF^\times$. From
(\ref{newK2KHcosetseq41}) -- (\ref{newK2KHcosetseq48}) we therefore obtain the following
preliminary decomposition, which is not disjoint.
\begin{align}
 \label{prelimnewK2KHcosetseq41} K^H&=\bigcup
  \limits_{u\in\OF^\times/(1+\p)}
  K_{l,m}\begin{bmatrix}1\\&u\\&&1\\&&&u^{-1}\end{bmatrix} K^\#(\p)\\
 \label{prelimnewK2KHcosetseq42} &\cup\bigcup\limits_{\substack{u\in\OF^\times/(1+\p)\\
  w\in\OF/\p}}K_{l,m}\begin{bmatrix}1\\&u\\&&1\\&&&u^{-1}\end{bmatrix}
  \begin{bmatrix}1\\w&1\\&&1&-w\\&&&1\end{bmatrix}s_1K^\#(\p)\\
 \label{prelimnewK2KHcosetseq43}
  &\cup\bigcup\limits_{u\in\OF^\times/(1+\p)}
  K_{l,m}\begin{bmatrix}1\\&u\\&&1\\&&&u^{-1}\end{bmatrix}
  s_2K^\#(\p)\\
 \label{prelimnewK2KHcosetseq44}
  &\cup\bigcup\limits_{\substack{u\in\OF^\times/(1+\p)\\w\in\OF/\p}}
  K_{l,m}\begin{bmatrix}1\\&u\\&&1\\&&&u^{-1}\end{bmatrix}
  \begin{bmatrix}1\\w&1\\&&1&-w\\&&&1\end{bmatrix}
  s_1s_2K^\#(\p)\\
 \label{prelimnewK2KHcosetseq45}
  &\cup\bigcup\limits_{u\in\OF^\times/(1+\p)}
  K_{l,m}\begin{bmatrix}1\\&u\\&&1\\&&&u^{-1}\end{bmatrix}
 s_2s_1K^\#(\p)\\
 \label{prelimnewK2KHcosetseq46}
  &\cup\bigcup\limits_{\substack{u\in\OF^\times/(1+\p)\\
  w\in\OF/\p}}K_{l,m}\begin{bmatrix}1\\&u\\&&1\\&&&u^{-1}\end{bmatrix}
  \begin{bmatrix}1\\w&1\\&&1&-w\\&&&1\end{bmatrix}s_1s_2s_1K^\#(\p)\\
 \label{prelimnewK2KHcosetseq47}
  &\cup\bigcup\limits_{u\in\OF^\times/(1+\p)}
  K_{l,m}\begin{bmatrix}1\\&u\\&&1\\&&&u^{-1}\end{bmatrix}s_2s_1s_2K^\#(\p)\\
 \label{prelimnewK2KHcosetseq48}
  &\cup\bigcup\limits_{\substack{u\in\OF^\times/(1+\p)\\w\in\OF/\p}}
  K_{l,m}\begin{bmatrix}1\\&u\\&&1\\&&&u^{-1}\end{bmatrix}
  \begin{bmatrix}1\\w&1\\&&1&-w\\&&&1\end{bmatrix}s_1s_2s_1s_2K^\#(\p).
\end{align}
\subsubsection{Step $2$: Support of $W^\#$}\label{support of W}
Recall
$$
\eta =\begin{bmatrix}1&0&&\\\alpha&1&&\\&&1&-\bar{\alpha}\\&&0&1\end{bmatrix},\qquad
 \text{where}\quad
 \renewcommand{\arraystretch}{2.1}\alpha:=\left\{\begin{array}{l@{\qquad\text{if }L}l}
 \displaystyle\frac{b+\sqrt{d}}{2c}&\text{ is a field},\\
 \displaystyle\Big(\frac{b+\sqrt{d}}{2c},\frac{b-\sqrt{d}}{2c}\Big)&=F\oplus F.
 \end{array}\right.
$$
We have assumed that $c\in\OF^\times$, so that $\alpha\in\OF_L$.
We have $\eta h(l,m) = h(l,m)\eta_m$ where for $m\geq0$ we define
\begin{equation}\label{etamdefeq}
 \eta_m=\begin{bmatrix}1\\\alpha\varpi^m&1\\&&1&-\bar\alpha\varpi^m\\&&&1\end{bmatrix}.
\end{equation}
Fix $l,m\geq0$, and let $r$ run through the representatives for
$K_{l,m}\backslash K^H/K^\#(\p)$ from
(\ref{prelimnewK2KHcosetseq41}) -- (\ref{prelimnewK2KHcosetseq48}).
We want to find out for which $r$ is $\eta h(l,m)r \in M(F)N(F)
K^\#(\P)$, since this set is the support of $W^\#$. Since
$h(l,m)\in M(F)$, this is equivalent to $\eta_m r \in M(F)N(F)
K^\#(\P)$. Hence, this condition depends only on $m\geq0$ and not
on the integer $l$. Recall that
$$
 K^\#(\P) =  \begin{bmatrix}1+\P&\P&\OF_L&\OF_L\\\OF_L&1+\P&\OF_L&\OF_L\\
 \P&\P&1+\P&\OF_L\\\P&\P&\P&1+\P\end{bmatrix}
 \begin{bmatrix}1\\&\mu&\\&&\mu\\&&&1\end{bmatrix},
 \qquad\mu\in\OF^\times.\
$$
\begin{enumerate}
 \item Let $r=\begin{bmatrix}1\\&u\\&&1\\&&&u^{-1}\end{bmatrix}$
  with $u \in\OF^\times/(1+\p)$. Then $\eta_m r \in M(F)N(F) K^\#(\P)$ for all $u, m$. More
  precisely,
  \begin{equation}\label{nsunramn1matrixidentity1eq}
   \eta_mr=\begin{bmatrix}1\\&u\\&&1\\&&&u^{-1}\end{bmatrix}
    \begin{bmatrix}1\\u^{-1}\varpi^m\alpha&1\\
     &&1&-u^{-1}\varpi^m\bar\alpha\\&&&1\\\end{bmatrix}\in M(F)K^\#(\P).
  \end{equation}
 \item
  Let $r=\begin{bmatrix}1\\&u\\&&1\\&&&u^{-1}\end{bmatrix}
  \begin{bmatrix}1\\w&1\\&&1&-w\\&&&1\end{bmatrix}s_1$ with
  $u\in\OF^\times/(1+\p)$ and $w\in\OF/\p$.
  If $\beta = \varpi^m \alpha + uw \in \OF_L^{\times}$, then
  \begin{align}\label{nsunramn1matrixidentity2eq}
   \eta_mr&=
   \begin{bmatrix}-\beta^{-1}u\\&\beta\\&&-\bar \beta u^{-1}\\&&&\bar \beta^{-1}\end{bmatrix}
   \begin{bmatrix}1&-\beta u^{-1}\\&1\\&&1\\&&\bar \beta u^{-1}&1\end{bmatrix}\nonumber\\
   &\hspace{10ex}\times
    \begin{bmatrix}1\\u\beta^{-1}&1\\&&1&-u \bar \beta^{-1}\\&&&1\end{bmatrix}
   \in M(F)N(F)K^\#(\P).
  \end{align}
  But if $\beta\notin\OF_L^\times$, then $\eta_mr\notin
  M(F)N(F)K^\#(\P)$ since the $(3,3)$-coefficient of any matrix product of the form
  $\tilde n^{-1}\tilde m^{-1}\eta_m r$, $\tilde m\in M(F)$, $\tilde n\in N(F)$,
  is always in $\beta \OF_L^{\times}$.
 \item Let $r=\begin{bmatrix}1\\&u\\&&1\\&&&u^{-1}\end{bmatrix} s_2$ with
  $u\in\OF^\times/(1+\p)$. Then $\eta_m r\notin M(F)N(F)K^\#(\P)$,
  since the $(3,3)$-coefficient of any matrix product of the form
  $\tilde n^{-1}\tilde m^{-1}\eta_m r$, $\tilde m\in M(F)$, $\tilde n\in N(F)$,
  is always zero.
 \item Let $r=\begin{bmatrix}1\\&u\\&&1\\&&&u^{-1}\end{bmatrix}
  \begin{bmatrix}1\\w&1&&\\&&1&-w\\&&&1\end{bmatrix}s_1s_2$ with
  $u\in\OF^\times/(1+\p)$ and $w\in\OF/\p$. Then $\eta_m r\notin M(F)N(F)K^\#(\P)$,
  since the $(3,3)$-coefficient of any matrix product of the form
  $\tilde n^{-1}\tilde m^{-1}\eta_m r$, $\tilde m\in M(F)$, $\tilde n\in N(F)$,
  is always zero.
 \item Let $r=\begin{bmatrix}1\\&u\\&&1\\&&&u^{-1}\end{bmatrix}s_2s_1$
  with  $u\in\OF^\times/(1+\p)$. Then $\eta_m r\notin M(F)N(F)K^\#(\P)$,
  since the $(3,3)$-coefficient of any matrix product of the form
  $\tilde n^{-1}\tilde m^{-1}\eta_m r$, $\tilde m\in M(F)$, $\tilde n\in N(F)$,
  is $\varpi^m\bar{\alpha}u^{-1}$ times the
  $(3,2)$-coefficient.
 \item Let $r=\begin{bmatrix}1\\&u\\&&1\\&&&u^{-1}\end{bmatrix}
  \begin{bmatrix}1&&&\\w&1&&\\&&1&-w\\&&&1\end{bmatrix}s_1s_2s_1$ with
  $u\in\OF^\times/(1+\p)$ and $w\in\OF/\p$.
  If $\beta=\varpi^m \alpha+ uw\in\P$, then
  \begin{equation}\label{nsunramn1matrixidentity3eq}
   \eta_mr=  \begin{bmatrix}1\\&&&u\\&&1\\&-1/u\end{bmatrix}
   \begin{bmatrix}1\\&1\\&\bar\beta u^{-1}&1\\\beta u^{-1}&&&1\end{bmatrix}\in M(F)K^\#(\P).
  \end{equation}
  But if $\beta \not\in\P$, then $\eta_mr\notin M(F)N(F)K^\#(\P)$
  since the $(3,2)$-coefficient of any matrix product of the form
  $\tilde n^{-1}\tilde m^{-1}\eta_m r$, $\tilde m\in M(F)$, $\tilde n\in N(F)$,
  is always in $\beta \OF_L^{\times}$.
 \item Let $r=\begin{bmatrix}1\\&u\\&&1\\&&&u^{-1}\end{bmatrix}s_2s_1s_2$
  with $u\in\OF^\times/(1+\p)$. Then $\eta_m r\notin M(F)N(F)K^\#(\P)$,
  since the $(3,3)$-coefficient of any matrix product of the form
  $\tilde n^{-1}\tilde m^{-1}\eta_m r$, $\tilde m\in M(F)$, $\tilde n\in N(F)$,
  is always zero.
 \item Let $r=\begin{bmatrix}1\\&u\\&&1\\&&&u^{-1}\end{bmatrix}
  \begin{bmatrix}1&&&\\w&1&&\\&&1&-w\\&&&1\end{bmatrix}s_1s_2s_1s_2$ with
  $u\in\OF^\times/(1+\p)$ and $w\in\OF/\p$. Then $\eta_m r\notin M(F)N(F)K^\#(\P)$,
  since the $(3,3)$-coefficient of any matrix product of the form
  $\tilde n^{-1}\tilde m^{-1}\eta_m r$, $\tilde m\in M(F)$, $\tilde n\in N(F)$,
  is always zero.
\end{enumerate}
Hence, for every $l \geq 0$ and $m \geq 0$, the double cosets that
contribute to the computation of the integral (\ref{Zvcalceq1}) are
\begin{align}
 &\bigcup_{u\in\OF^\times/(1+\p)}K_{l,m}
  \begin{bmatrix}1\\&u\\&&1\\&&&u^{-1}\end{bmatrix}K^\#(\p) \label{Wsuppcoset1}\\
 \cup&\bigcup_{u\in\OF^\times/(1+\p)}\:\bigcup_{\substack{w\in\OF/\p\\
  \varpi^m \alpha + uw \in \OF_L^{\times}}} K_{l,m}
  \begin{bmatrix}1\\&u\\&&1\\&&&u^{-1}\end{bmatrix}
  \begin{bmatrix}1\\w&1\\&&1&-w\\&&&1\end{bmatrix}s_1K^\#(\p) \label{Wsuppcoset2} \\
 \cup&\bigcup_{u\in\OF^\times/(1+\p)}\:\bigcup_{\substack{
  w\in\OF/\p \\ \varpi^m \alpha + uw \in \P}} K_{l,m}
  \begin{bmatrix}1\\&u\\&&1\\&&&u^{-1}\end{bmatrix}
  \begin{bmatrix}1&&&\\w&1&&\\&&1&-w\\&&&1\end{bmatrix}s_1s_2s_1K^\#(\p).\label{Wsuppcoset6}
\end{align}
By ii) of Lemma \ref{quadraticextdisclemma}, the condition
$\varpi^m \alpha + uw \in \P$ cannot be satisfied if $m=0$. Hence,
for $\underline{m=0}$ the double cosets that contribute to the
computation of the integral (\ref{Zvcalceq1}) are
\begin{align}
 S &=\bigcup_{u\in\OF^\times/(1+\p)}K_{l,0}
  \begin{bmatrix}1\\&u\\&&1\\&&&u^{-1}\end{bmatrix}K^\#(\p) \label{m0Wsuppcoset1}\\
 &\cup\bigcup_{u\in\OF^\times/(1+\p)}\:\bigcup_{\substack{w\in\OF/\p\\
  \alpha + uw \in \OF_L^{\times}}} K_{l,0}
  \begin{bmatrix}1\\&u\\&&1\\&&&u^{-1}\end{bmatrix}
  \begin{bmatrix}1\\w&1\\&&1&-w\\&&&1\end{bmatrix}s_1K^\#(\p).\label{m0Wsuppcoset2}
\end{align}
For $m>0$ the condition $\varpi^m \alpha + uw \in \OF_L^\times$
(resp.\ $\varpi^m \alpha + uw \in \P$) is satisfied if and only if
$w\in\OF^\times$ (resp.\ $w\in\p$). Hence, for $\underline{m>0}$
the double cosets that contribute to the computation of the
integral (\ref{Zvcalceq1}) are
\begin{align}
 S &=\bigcup_{u\in\OF^\times/(1+\p)}K_{l,m}
  \begin{bmatrix}1\\&u\\&&1\\&&&u^{-1}\end{bmatrix}K^\#(\p) \label{mposWsuppcoset1}\\
 &\cup\bigcup_{u\in\OF^\times/(1+\p)}\:\bigcup_{w\in(\OF/\p)^\times}
  K_{l,m}\begin{bmatrix}1\\&u\\&&1\\&&&u^{-1}\end{bmatrix}
  \begin{bmatrix}1\\w&1\\&&1&-w\\&&&1\end{bmatrix}s_1K^\#(\p) \label{mposWsuppcoset2} \\
 & \cup\bigcup_{u\in\OF^\times/(1+\p)} K_{l,m}
  \begin{bmatrix}1\\&u\\&&1\\&&&u^{-1}\end{bmatrix}s_1s_2s_1K^\#(\p).\label{mposWsuppcoset6}
\end{align}

\subsubsection{Step $3$: Disjointness of double cosets}
We will now investigate possible overlaps between the double cosets given in
(\ref{m0Wsuppcoset1}) and (\ref{m0Wsuppcoset2}) (for $m=0$) resp.\
(\ref{mposWsuppcoset1}), (\ref{mposWsuppcoset2}) and (\ref{mposWsuppcoset6}) (for $m>0$). Recall that $K_{l,m}=h(l,m)^{-1}R(F)h(l,m)\cap K^H$.
\nl
\underline{\bf The case $\mathbf{m=0}$}
\nl
We will now assume $m=0$ and find all equivalences between the double cosets in
(\ref{m0Wsuppcoset1}) and (\ref{m0Wsuppcoset2}).
Let a double coset from (\ref{m0Wsuppcoset2}) be given.
Set $v = a + b(uw) + c(uw)^2$. We claim that the condition $\alpha+uw \in \OF_L^{\times}$
in (\ref{m0Wsuppcoset2}) forces $v \in \OF^{\times}$.
First observe that we have the following identity,
$$
 a+b(uw)+c(uw)^2=-c(\alpha+uw)(\alpha-(uw+bc^{-1})).
$$
Hence, if $v\in\p$, then it would follow that
$\alpha-(uw+bc^{-1})\in\p\OF_L=\P$. By Lemma \ref{quadraticextdisclemma} ii),
this is impossible. It follows that indeed $v\in\OF^\times$. Now let
$y=-u/v$ and $x=-(u/v)(cwu+b/2)$. Let $g=\begin{bmatrix}x+yb/2&yc\\-ya&x-yb/2\end{bmatrix}$.
Then
\begin{align}\label{matrixidentity1eq}
 \nonumber&h(l,0)^{-1}\mat{g}{}{}{\det(g)\,^tg^{-1}} h(l,0)
  \begin{bmatrix}1\\&-u\\&&1\\&&&-u^{-1}\end{bmatrix}\\
 &\qquad = \begin{bmatrix}1\\&u\\&&1\\&&&u^{-1}\end{bmatrix}
  \begin{bmatrix}1\\w&1\\&&1&-w\\&&&1\end{bmatrix}s_1
  \begin{bmatrix}1&0&&\\-\frac{u(b+cuw)}{v}&\frac{cu^2}v&&\\
  &&\frac{cu^2}v&\frac{u(b+cuw)}{v}\\&&0&1\end{bmatrix}.
\end{align}
Since $v \in \OF^{\times}$, the rightmost matrix is in $K^\#(\p)$.
This identity shows that \emph{all} double cosets in (\ref{m0Wsuppcoset2}) are equivalent
to double cosets in (\ref{m0Wsuppcoset1}).
 So far we have shown that the set $S$ in (\ref{m0Wsuppcoset1}) and
(\ref{m0Wsuppcoset2}) reduces to
$$
 S=\bigcup\limits_{u\in\OF^\times/(1+\p)}
  K_{l,0}\begin{bmatrix}1\\&u\\&&1\\&&&u^{-1}\end{bmatrix}K^\#(\p).
$$
We will now show that this is a disjoint union. Let
$$
   h_1=\begin{bmatrix}1\\&u_1\\&&1\\&&&u_1^{-1}\end{bmatrix},\qquad
   h_2=\begin{bmatrix}1\\&u_2\\&&1\\&&&u_2^{-1}\end{bmatrix},
$$
and assume that $K_{l,0}h_1K^\#(\p)=K_{l,0}h_2K^\#(\p)$.
Then there exists $r\in R(F)$ such that
$$
 A = h_2^{-1}h(l,0)^{-1}rh(l,0)h_1\in K^\#(\p).
$$
Dividing the $(1,1)$ coefficient of $A$ by the $(4,4)$
coefficient, we get $u_1/u_2 \in 1+\p$. Hence $h_1$ and $h_2$
define the same double coset if and only if $u_1=u_2$. It follows
that
$$
 S=\bigsqcup\limits_{u\in\OF^\times/(1+\p)}
  K_{l,0}\begin{bmatrix}1\\&u\\&&1\\&&&u^{-1}\end{bmatrix}K^\#(\p),
$$
as claimed.
\nl
\underline{\bf The case $\mathbf{m>0}$}
\nl
We will now assume $m>0$ and find all equivalences between the double cosets in
(\ref{mposWsuppcoset1}), (\ref{mposWsuppcoset2}) and (\ref{mposWsuppcoset6}).
An argument similar to the one above shows that all the double cosets in (\ref{mposWsuppcoset1}) are disjoint.

\begin{description}
\item[Equivalence of double cosets from (\ref{mposWsuppcoset1})
and (\ref{mposWsuppcoset2}):] Let a double coset from
(\ref{mposWsuppcoset2}) be given. Let
$x=\frac{b\varpi^m}{2cw^2u}-\frac{1}w$ and
$y=-\frac{\varpi^{m}}{cw^2u}$. Let
$g=\begin{bmatrix}x+yb/2&yc\\-ya&x-yb/2\end{bmatrix}$. Then
\begin{align}\label{matrixidentity2eq}
 \nonumber&h(l,m)^{-1}\mat{g}{}{}{\det(g)\,^tg^{-1}}h(l,m)
 \begin{bmatrix}1\\&-u\\&&1\\&&&-u^{-1}\end{bmatrix}\\
 &\qquad=\begin{bmatrix}1\\&u\\&&1\\&&&u^{-1}\end{bmatrix}
 \begin{bmatrix}1\\w&1\\&&1&-w\\&&&1\end{bmatrix}s_1
 \begin{bmatrix}1+\frac{\varpi^{2m}a}{cw^2u^2}
 &-\frac{b\varpi^m}{cw^2u}&&\\-\frac 1w&\frac{1}{w^2}\\&&\frac{1}{w^2}&\frac 1w\\
 &&\frac{b\varpi^m}{cw^2u}&1+\frac{\varpi^{2m}a}{cw^2u^2}
 \end{bmatrix}.
\end{align}
The rightmost matrix lies in $K^\#(\p)$ since $w \in \OF^{\times}$. Hence cosets of
(\ref{mposWsuppcoset2}) all coincide with cosets from
(\ref{mposWsuppcoset1}).

\item[Equivalence of double cosets from (\ref{mposWsuppcoset1})
and (\ref{mposWsuppcoset6}):]

Let $h_1$ be a double coset representative obtained in
(\ref{mposWsuppcoset1})  and let $h_2$ be a double coset
representative obtained in (\ref{mposWsuppcoset6}). Then the
double cosets are not equivalent, since, for every $k \in
K_{l,m}$, the $(2,2)$ entry of the matrix $A = h_2^{-1} k h_1$ is
zero.
\item[Equivalence amongst double coset from (\ref{mposWsuppcoset6}):]
We will show that the double cosets in (\ref{mposWsuppcoset6}) are disjoint. Let
$$
 h_1 =\begin{bmatrix}1\\&u_1\\&&1\\&&&u_1^{-1}\end{bmatrix}s_1 s_2
  s_1, \qquad h_2 =
  \begin{bmatrix}1\\&u_2\\&&1\\&&&u_2^{-1}\end{bmatrix}s_1 s_2 s_1.
$$
Assume that $K_{l,m}h_1K^\#(\p)=K_{l,m}h_2K^\#(\p)$.
Then there exists an $r\in R(F)$ such that
$$
 A=h_2^{-1}h(l,m)^{-1}r h(l,m)h_1\in K^\#(\p).
$$
Let $r=gY$, with $g\in T(F)$ and $Y\in U(F)$. We write an element
of $T(F)$ as $g = \begin{bmatrix}x+yb/2&yc\\-ya&x-yb/2\end{bmatrix}$ with $x,y\in F$.
Looking at the $(2,3)$ coefficient of $A$, we see that $y \in \p$
(since $m > 0$). The $(1,1)$ coefficient gives us that $x+\frac
b2y \in 1 + \p$ and hence $x-\frac b2y \in 1+\p$. Looking at the
$(4,4)$ coefficient, we get $u_1/u_2 \in 1+\p$, which implies
$u_1=u_2$. Hence, the double cosets in (\ref{mposWsuppcoset6}) are
disjoint.
\end{description}
We summarize the results of this section in the following lemma.
\begin{lemma}\label{relevant disjoint double cosets}
 The following are the disjoint double cosets in
 $\{K_{l,m} k K^\#(\p): k \in K^H,\: l,m \geq 0 \}$ that have a non-trivial
 intersection with the support of $W^\#$.
 \begin{align*}
  &\bigsqcup_{\substack{l\geq0 \\m\geq0}}\:\bigsqcup_{u\in\OF^\times/(1+\p)}K_{l,m}
  \begin{bmatrix}1\\&u\\&&1\\&&&u^{-1}\end{bmatrix}K^\#(\p)\\
  \sqcup&\bigsqcup_{\substack{l\geq0 \\m>0}}\:\bigsqcup_{u\in\OF^\times/(1+\p)}K_{l,m}
  \begin{bmatrix}1\\&u\\&&1\\&&&u^{-1}\end{bmatrix}s_1 s_2s_1K^\#(\p).
 \end{align*}
\end{lemma}
\subsection{Volume computations}
In Lemma \ref{relevant disjoint double cosets} we obtained the double coset
representatives $\{s_j'''\}$ needed to evaluate the integral (\ref{Zvcalceq4}).
In this section we will compute the corresponding volumes. More precisely,
we have to compute
\begin{alignat}{2}
 V_1^{l,m} &:= \int_{K_{l,m}\backslash K_{l,m}
  \left[\begin{smallmatrix}1\\&u\\&&1\\&&&u^{-1}\end{smallmatrix}\right]
  K^\#(\p)}dh &&\text{for }l \geq 0,\:m \geq 0,\label{k1voleq1}\\
 V_2^{l,m} &:= \int_{K_{l,m}\backslash K_{l,m}
  \left[\begin{smallmatrix}1\\&u\\&&1\\&&&u^{-1}\end{smallmatrix}\right]
  s_1s_2s_1K^\#(\p)}dh \qquad &&\text{for }l \geq 0,\:m > 0.\label{k1voleq2}
\end{alignat}
In the notations above we have suppressed the dependence on $u$
since it will turn out that these volumes are independent of $u$.
We will first show that the calculation reduces to the calculation
of volumes of certain compact subgroups of $\GL_2(F)$. First we
make some general remarks that apply to both cases. The volumes
(\ref{k1voleq1}), (\ref{k1voleq2}) are of the form
$$
 \int\limits_{K_{l,m}\backslash K_{l,m}AK^\#(\p)}dh,
$$
with some $A\in K^H$. Let $\chi_1:\:K_{l,m}\backslash K^H\rightarrow\C$ be the
characteristic function of $K_{l,m}\backslash K_{l,m}AK^\#(\p)$, and let
$\delta_1:\:K^H\rightarrow\C$ be the characteristic function of $AK^\#(\p)$.
\begin{lemma}\label{generalvolumecalclemma}
 For all $g\in K^H$ we have
 \begin{equation}\label{generalvolumecalclemmaeq1}
  \int\limits_{K_{l,m}}\delta_1(tg)\,dt=\chi_1(\dot g)\int\limits_{K_{l,m}\cap
  \big(AK^\#(\p)A^{-1}\big)}dt,
 \end{equation}
 where $\dot g$ denotes the image of $g$ in $K_{l,m}\backslash K^H$.
\end{lemma}
{\bf Proof.} First assume that $g\notin K_{l,m}AK^\#(\p)$. Then
$tg\notin AK^\#(\p)$ for all $t\in K_{l,m}$, and hence the left side is zero.
The right side is also zero by definition of $\chi_1$. Thus the equality holds
under our assumption. Now assume that $g\in K_{l,m}AK^\#(\p)$.
In this case $\chi_1(\dot g)=1$. Write
$g=kA\kappa$ with $k\in K_{l,m}$ and $\kappa\in K^\#(\p)$. We have
\begin{alignat*}{2}
 tg\in AK^\#(\p)\;&\Longleftrightarrow\; tkA\kappa\in AK^\#(\p)
 &&\Longleftrightarrow\;tkA\in AK^\#(\p)\\
 &\Longleftrightarrow\;tk\in AK^\#(\p)A^{-1}
 &\;&\Longleftrightarrow\;t\in \big(AK^\#(\p)A^{-1}\big)k^{-1}.
\end{alignat*}
Hence the left side equals
$$
 \int\limits_{K_{l,m}\cap \big(AK^\#(\p)A^{-1}\big)k^{-1}}dt.
$$
But since $k\in K_{l,m}$, this integral equals
$\int_{K_{l,m}\cap(AK^\#(\p)A^{-1})}dt$.
This proves the lemma.\qed

\vspace{2ex}
Integrating both sides of (\ref{generalvolumecalclemmaeq1}) over
$K_{l,m}\backslash K^H$, we obtain
\begin{equation}\label{generalvolumecalclemmaeq2}
 \int\limits_{K^H}\delta_1(g)\,dg
 =\big(\int\limits_{K_{l,m}\backslash K_{l,m}AK^\#(\p)}dh\big)
 \big(\int\limits_{K_{l,m}\cap\big(AK^\#(\p)A^{-1}\big)}dt\big),
\end{equation}
so that
\begin{equation}\label{generalvolumecalclemmaeq3}
 \int\limits_{K_{l,m}\backslash K_{l,m}AK^\#(\p)}dh
 ={\rm vol}(K^\#(\p))
 \Big(\int\limits_{K_{l,m}\cap\big(AK^\#(\p)A^{-1}\big)}dt\Big)^{-1}.
\end{equation}
Note that
\begin{equation}\label{Ksharppvolumeeq}
 {\rm vol}(K^\#(\p))=\frac1{(q-1)^2(1+2q+2q^2+2q^3+q^4)}
 =\frac1{(q^2-1)(q^4-1)}
\end{equation}
from (\ref{newK2KHcosetseq41}) -- (\ref{newK2KHcosetseq48})
and the fact that ${\rm vol}(K^H)=1$. Hence we are reduced to computing
\begin{equation}\label{generalvolumecalclemmaeq4}
 V(l,m,A):=\int\limits_{K_{l,m}\cap\big(AK^\#(\p)A^{-1}\big)}dt.
\end{equation}
In both the volumes
(\ref{k1voleq1}), (\ref{k1voleq2}) we have $A=Bw$, where $w$ is a Weyl group
element and
$$ B=\begin{bmatrix}1\\&u\\&&1\\&&&u^{-1}\end{bmatrix}.$$
We have
\begin{equation}\label{generalvolumecalclemmaeq5}
 V(l,m,A)=\int\limits_{\big(B^{-1}K_{l,m}B\big)\cap\big(wK^\#(\p)w^{-1}\big)}dt = \int\limits_{\big(B^{-1}h(l,m)^{-1}R(F)h(l,m)B\big)
  \cap\big(wK^\#(\p)w^{-1}\big)}dt.
\end{equation}
The relevant Weyl group elements are $w=1$ and $w=s_1s_2s_1$. We have
\begin{align}\label{K2capHeq2}
 K^\#(\p)&=\begin{bmatrix}1+\p&\p&\OF&\OF\\\OF&\OF^\times&\OF&\OF\\
 \p&\p&\OF^\times&\OF\\\p&\p&\p&1+\p\end{bmatrix},\\
 s_1s_2s_1K^\#(\p)s_1s_2s_1
  &=\begin{bmatrix}1+\p&\OF&\OF&\p\\\p&1+\p&\p&\p\\
 \p&\OF&\OF^\times&\p\\\OF&\OF&\OF&\OF^\times\end{bmatrix}.
\end{align}
 We have to find the intersections of these
compact groups with $B^{-1}h(l,m)^{-1}R(F)h(l,m)B$. Note that, mod
$\p$, the upper left block of any matrix in any of these compact
subgroups, lies in $\GL_2(\OF)$. Let $L_{K,w}$ be the subgroup of
$\GL_2(\OF)$ occurring in the upper left block of
$wK^\#(\p)w^{-1}$, and let $N_{K,w}$ be the compact subgroup of
$F^3$ occurring in the upper right block. Write a given element of
$R(F)$ as $tn$ with $t\in T(F)$ and $n\in U(F)$. Then
$$
 B^{-1}h(l,m)^{-1}(tn)h(l,m)B
 =\big(B^{-1}h(l,m)^{-1}t\,h(l,m)B\big)\big(B^{-1}h(l,m)^{-1}n\,h(l,m)B\big).
$$
A direct computation shows that this element lies in $wK^\#(\p)w^{-1}$ if and only if
\begin{equation}\label{volcomp1cond1eq}
 \text{the upper left block of }B^{-1}h(l,m)^{-1}t\,h(l,m)B\;\text{ lies in }\;L_{K,w}
\end{equation}
and
\begin{equation}\label{volcomp1cond2eq}
 \text{the upper right block of }B^{-1}h(l,m)^{-1}n\,h(l,m)B\;\text{ lies in }\;N_{K,w}.
\end{equation}
We have $B=\mat{g}{}{}{^tg^{-1}}$ with $g=\mat{1}{}{}{u}
\in\GL_2(\OF)$. The conditions (\ref{volcomp1cond1eq}) and
(\ref{volcomp1cond2eq}) become
\begin{equation}\label{volcomp1cond1aeq}
 \mat{1}{}{}{u^{-1}}\mat{\varpi^{-m}}{}{}{1}t\mat{\varpi^m}{}{}{1}\mat{1}{}{}{u}\in L_{K,w}
\end{equation}
and
\begin{equation}\label{volcomp1cond2aeq}
 \mat{1}{}{}{u^{-1}}\mat{\varpi^{-2m-l}}{}{}{\varpi^{-m-l}}X\mat{1}{}{}{\varpi^m}
 \mat{1}{}{}{u^{-1}}\in N_{K,w},\qquad\text{where }n=\mat{1}{X}{}{1}.
\end{equation}
It follows that
\begin{align*}
 &{\rm vol}\big(\{X\in F^3:\:\mat{1}{}{}{u^{-1}}
  \mat{\varpi^{-2m-l}}{}{}{\varpi^{-m-l}}X\mat{1}{}{}{\varpi^m}\,
  \mat{1}{}{}{u^{-1}}\in N_{K,w}\}\big)\\
 &\qquad={\rm vol}\big(\{X\in F^3:\:X\in
  \mat{\varpi^{2m+l}}{}{}{\varpi^{m+l}}\mat{1}{}{}{u}N_{K,w}
  \mat{1}{}{}{u}\mat{1}{}{}{\varpi^{-m}}\}\big)\\
 &\qquad={\rm vol}
  \big(\mat{\varpi^{2m+l}}{}{}{\varpi^{m+l}}\mat{1}{}{}{u}N_{K,w}\mat{1}{}{}{u}
  \mat{1}{}{}{\varpi^{-m}}\}\big)\\
 &\qquad={\rm vol}\big(\mat{1}{}{}{u}\mat{\varpi^{2m+l}}{}{}{\varpi^{m+l}}N_{K,w}
  \mat{1}{}{}{\varpi^{-m}}\mat{1}{}{}{u}\}\big)\\
 &\qquad={\rm vol}\big(\mat{\varpi^{2m+l}}{}{}{\varpi^{m+l}}N_{K,w}
  \mat{1}{}{}{\varpi^{-m}}\}\big)\\
 &\qquad=q^{-3m-3l}{\rm vol}(N_{K,w}).
\end{align*}
The volume of $N_{K,w}$ is $1$ if $w=1$ and $q^{-2}$ if
$w=s_1s_2s_1$. Let us set
 $$
  L_{K,w}=\left\{\begin{array}{l@{\qquad\mbox{if }}l}
  \mat{1+\p}{\p}{\OF}{\OF^\times}&w=1,\\[2ex]
  \mat{1+\p}{\OF}{\p}{1+\p}&w=s_1s_2s_1,
  \end{array}\right.
 $$
as above. Let
\begin{align}\label{generalvolumecalclemma2eq2}
  T_{m,w}&=\{t\in T(F):\:\mat{1}{}{}{u^{-1}}
  \mat{\varpi^{-m}}{}{}{1}t\mat{\varpi^m}{}{}{1}\mat{1}{}{}{u}\in L_{K,w}\}\nonumber\\
  &=T(F)\cap\mat{\varpi^m}{}{}{1}\mat{1}{}{}{u}L_{K,w}\mat{1}{}{}{u^{-1}}
  \mat{\varpi^{-m}}{}{}{1}\nonumber \\
  &=T(F)\cap\mat{\varpi^m}{}{}{1}L_{K,w}\mat{\varpi^{-m}}{}{}{1}.
\end{align}
We summarize the above considerations.
\begin{lemma}\label{generalvolumecalclemma2}
 Let $l$ and $m$ be non-negative integers. Let $w\in\{1,s_1s_2s_1\}$,
 and set $\delta=0$ if $w=1$, and $\delta=2$ if $w=s_1s_2s_1$. Let
 $A={\rm diag}(1,u,1,u^{-1})w$, where $u\in\OF^\times$. Then
 $$
  \int\limits_{K_{l,m}\backslash K_{l,m}AwK^\#(\p)}dh={\rm vol}(K^\#(\p))\,
  {\rm vol}(T_{m,w})^{-1}\,q^{3(m+l)+\delta}.
 $$
 Here, ${\rm vol}(K^\#(\p))=\big((q^2-1)(q^4-1)\big)^{-1}$.\qed
\end{lemma}
Thus we are reduced to computing the volumes of the groups $T_{m,1}$ for all $m \geq 0$
and $T_{m,s_1s_2s_1}$ for all $m > 0$.
\begin{lemma}\label{GL2volumecalclemma1}
 For any $m\geq0$ we have
 $$
  {\rm vol}(T_{m,1})^{-1}=(q-1)\Big(1-\Big(\frac L\p\Big)q^{-1}\Big)q^{m+1}.
 $$
\end{lemma}
{\bf Proof.} By definition,
$$
 T_{m,1}=T(F)\cap\mat{\varpi^m}{}{}{1}\mat{1+\p}{\p}{\OF}{\:\OF^\times}
 \mat{\varpi^{-m}}{}{}{1}.
$$
Since
$$
 \mat{\OF^\times}{\p}{\OF\:}{\:\OF^\times}=\bigsqcup\limits_{u\in\OF^\times/(1+\p)}
 \mat{u}{}{}{u}\mat{1+\p}{\p}{\OF}{\:\OF^\times}
$$
and $\mat{u}{}{}{u}\in T(F)$, we have
\begin{align*}
 \int\limits_{T(F)\cap
  \left[\begin{smallmatrix}\varpi^m\\&1\end{smallmatrix}\right]
  \left[\begin{smallmatrix}\OF^\times&\p\\\OF&\OF^{\times}\end{smallmatrix}\right]
  \left[\begin{smallmatrix}\varpi^{-m}\\&1\end{smallmatrix}\right]}dt
 &=\sum_{u\in\OF^\times/(1+\p)}\;(\int\limits_{T(F)\cap
  \left[\begin{smallmatrix}\varpi^m\\&1\end{smallmatrix}\right]
  \left[\begin{smallmatrix}u\\&u\end{smallmatrix}\right]
  \left[\begin{smallmatrix}1+\p&\p\\\OF&\OF^{\times}\end{smallmatrix}\right]
  \left[\begin{smallmatrix}\varpi^{-m}\\&1\end{smallmatrix}\right]}dt)\\
 &=\sum_{u\in\OF^\times/(1+\p)}\;(\int\limits_{T(F)\cap
  \left[\begin{smallmatrix}\varpi^m\\&1\end{smallmatrix}\right]
  \left[\begin{smallmatrix}1+\p&\p\\\OF&\OF^{\times}\end{smallmatrix}\right]
  \left[\begin{smallmatrix}\varpi^{-m}\\&1\end{smallmatrix}\right]}dt)\\
 &=(q-1)\;(\int\limits_{T(F)\cap
  \left[\begin{smallmatrix}\varpi^m\\&1\end{smallmatrix}\right]
  \left[\begin{smallmatrix}1+\p&\p\\\OF&\OF^{\times}\end{smallmatrix}\right]
  \left[\begin{smallmatrix}\varpi^{-m}\\&1\end{smallmatrix}\right]}dt).
\end{align*}
Therefore,
$$
 \Big(\int\limits_{T(F)\cap
 \left[\begin{smallmatrix}\varpi^m\\&1\end{smallmatrix}\right]
 \left[\begin{smallmatrix}1+\p&\p\\\OF&\OF^{\times}\end{smallmatrix}\right]
 \left[\begin{smallmatrix}\varpi^{-m}\\&1\end{smallmatrix}\right]}dt\Big)^{-1}
 =(q-1)\Big(\int\limits_{T(F)\cap
 \left[\begin{smallmatrix}\varpi^m\\&1\end{smallmatrix}\right]
 \left[\begin{smallmatrix}\OF^\times&\p\\\OF&\OF^{\times}\end{smallmatrix}\right]
 \left[\begin{smallmatrix}\varpi^{-m}\\&1\end{smallmatrix}\right]}dt\Big)^{-1}.
$$
Note that the group
$T(F)\cap\mat{\varpi^m}{}{}{1}\GL_2(\OF)\mat{\varpi^{-m}}{}{}{1}$
lies in $\OF_L^\times$, since the determinants of these matrices
lie in $\OF^\times$ and the trace lies in $\OF$. As in \cite{Fu},
p.\ 202, we define a subring $\OF_m$ of $\OF_L$ by
$$
 \OF_m:=\OF_L\cap\mat{\varpi^m}{}{}{1}M_2(\OF)\mat{\varpi^{-m}}{}{}{1}
$$
In addition, we define a smaller subring
$$
 \OF'_m:=\OF_L\cap\mat{\varpi^m}{}{}{1}\mat{\OF}{\p}{\OF}{\OF}\mat{\varpi^{-m}}{}{}{1}.
$$
We normalize the measure so that ${\rm vol}(\OF_L^\times)=1$. Hence, we have
$$
 \Big(\int\limits_{T(F)\cap
 \left[\begin{smallmatrix}\varpi^m\\&1\end{smallmatrix}\right]
 \left[\begin{smallmatrix}\OF^\times&\p\\\OF&\OF^{\times}\end{smallmatrix}\right]
 \left[\begin{smallmatrix}\varpi^{-m}\\&1\end{smallmatrix}\right]}dt\Big)^{-1}
 =(\OF_L^\times:(\OF'_m)^\times).
$$
From (\ref{OLexpliciteq}), we have the integral basis
$$
 \OF_L=\OF+\OF\xi_0=\{\mat{x}{yc}{-ya\:}{\;x-yb}\:x,y\in\OF\},
 \qquad\xi_0=\mat{0}{c}{-a}{-b}.
$$
Such an element lies in
$\mat{\varpi^m}{}{}{1}M_2(\OF)\mat{\varpi^{-m}}{}{}{1}$ if and
only if $y\in\p^m$. Similarly, such an element lies in
$\mat{\varpi^m}{}{}{1}\mat{\OF}{\p}{\OF}{\OF}\mat{\varpi^{-m}}{}{}{1}$
if and only if $y\in\p^{m+1}$. Therefore,
$$
 \OF_m=\{x+\varpi^my\xi_0:\:x,y\in\OF\}
$$
and
$$
 \OF'_m=\{x+\varpi^{m+1}y\xi_0:\:x,y\in\OF\}.
$$
Hence $\OF'_m=\OF_{m+1}$, so that $(\OF_L^\times:(\OF'_m)^\times)
=(\OF_L^\times:(\OF_{m+1})^\times)$. By Lemma (3.5.3) of
\cite{Fu},
$$
 (\OF_L^\times:(\OF'_m)^\times)=\Big(1-\Big(\frac L\p\Big)q^{-1}\Big)q^{m+1}.
$$
This concludes the proof.\qed
\begin{lemma}\label{GL2volumecalclemma2}
 For any $m>0$ we have
 $$
  {\rm vol}(T_{m,s_1s_2s_1})^{-1}=(q-1)\Big(1-\Big(\frac L\p\Big)q^{-1}\Big)q^m.
 $$
\end{lemma}
{\bf Proof.} By definition
$$
 T_{m,s_1s_2s_1}=T(F)\cap\mat{\varpi^m}{}{}{1}
 \mat{1+\p}{\OF}{\p}{1+\p}\mat{\varpi^{-m}}{}{}{1}.
$$
We claim that
$$
 T_{m,s_1s_2s_1}=T(F)\cap\mat{\varpi^m}{}{}{1}
 \mat{1+\p}{\OF}{\p}{\OF^\times}\mat{\varpi^{-m}}{}{}{1}.
$$
Indeed, assume that
$$
 \mat{x+yb/2}{yc}{-ya}{x-yb/2}\in\mat{\varpi^m}{}{}{1}
 \mat{1+\p}{\OF}{\p}{\OF^\times}\mat{\varpi^{-m}}{}{}{1}.
$$
Then in particular $yc\in\p^m\subset\p$, since $m>0$. Since $c$ is
a unit, we get $y\in\p$. Thus, $x+yb/2\in1+\p$ implies
$x-yb/2\in1+\p$, as claimed. Since
$$
 \mat{\OF^\times}{\OF}{\p\:}{\:\OF^\times}=\bigsqcup\limits_{u\in\OF^\times/(1+\p)}
 \mat{u}{}{}{u}\mat{1+\p}{\OF}{\p}{\:\OF^\times}
$$
and $\mat{u}{}{}{u}\in T(F)$, we have
\begin{align*}
 \int\limits_{T(F)\cap
  \left[\begin{smallmatrix}\varpi^m\\&1\end{smallmatrix}\right]
  \left[\begin{smallmatrix}\OF^\times&\OF\\\p&\OF^{\times}\end{smallmatrix}\right]
  \left[\begin{smallmatrix}\varpi^{-m}\\&1\end{smallmatrix}\right]}dt
 &=\sum_{u\in\OF^\times/(1+\p)}\;(\int\limits_{T(F)\cap
  \left[\begin{smallmatrix}\varpi^m\\&1\end{smallmatrix}\right]
  \left[\begin{smallmatrix}u\\&u\end{smallmatrix}\right]
  \left[\begin{smallmatrix}1+\p&\OF\\\p&\OF^{\times}\end{smallmatrix}\right]
  \left[\begin{smallmatrix}\varpi^{-m}\\&1\end{smallmatrix}\right]}dt)\\
 &=\sum_{u\in\OF^\times/(1+\p)}\;(\int\limits_{T(F)\cap
  \left[\begin{smallmatrix}\varpi^m\\&1\end{smallmatrix}\right]
  \left[\begin{smallmatrix}1+\p&\OF\\\p&\OF^{\times}\end{smallmatrix}\right]
  \left[\begin{smallmatrix}\varpi^{-m}\\&1\end{smallmatrix}\right]}dt)\\
 &=(q-1)\;(\int\limits_{T(F)\cap
  \left[\begin{smallmatrix}\varpi^m\\&1\end{smallmatrix}\right]
  \left[\begin{smallmatrix}1+\p&\OF\\\p&\OF^{\times}\end{smallmatrix}\right]
  \left[\begin{smallmatrix}\varpi^{-m}\\&1\end{smallmatrix}\right]}dt).
\end{align*}
Therefore,
$$
 \Big(\int\limits_{T(F)\cap
 \left[\begin{smallmatrix}\varpi^m\\&1\end{smallmatrix}\right]
 \left[\begin{smallmatrix}1+\p&\OF\\\p&\OF^{\times}\end{smallmatrix}\right]
 \left[\begin{smallmatrix}\varpi^{-m}\\&1\end{smallmatrix}\right]}dt\Big)^{-1}
 =(q-1)\Big(\int\limits_{T(F)\cap
 \left[\begin{smallmatrix}\varpi^m\\&1\end{smallmatrix}\right]
 \left[\begin{smallmatrix}\OF^\times&\OF\\\p&\OF^{\times}\end{smallmatrix}\right]
 \left[\begin{smallmatrix}\varpi^{-m}\\&1\end{smallmatrix}\right]}dt\Big)^{-1}.
$$
Let $\OF_m$ be the subring of $\OF_L$ as defined in the proof of Lemma \ref{GL2volumecalclemma1}. In addition, we define another subring
$$
 \OF''_m:=\OF_L\cap\mat{\varpi^m}{}{}{1}\mat{\OF}{\OF}{\p}{\OF}\mat{\varpi^{-m}}{}{}{1}.
$$
Since ${\rm vol}(\OF_L^\times)=1$, we have
$$
 \Big(\int\limits_{T(F)\cap
 \left[\begin{smallmatrix}\varpi^m\\&1\end{smallmatrix}\right]
 \left[\begin{smallmatrix}\OF^\times&\OF\\\p&\OF^{\times}\end{smallmatrix}\right]
 \left[\begin{smallmatrix}\varpi^{-m}\\&1\end{smallmatrix}\right]}dt\Big)^{-1}
 =(\OF_L^\times:(\OF''_m)^\times).
$$
As above we have the integral basis
$\OF_L=\OF+\OF\xi_0=\{\mat{x}{yc}{-ya\:}{\;x-yb}\:x,y\in\OF\}$.
Such an element lies in
$\mat{\varpi^m}{}{}{1}M_2(\OF)\mat{\varpi^{-m}}{}{}{1}$ if and
only if $y\in\p^m$. Similarly, such an element lies in
$\mat{\varpi^m}{}{}{1}\mat{\OF}{\OF}{\p}{\OF}\mat{\varpi^{-m}}{}{}{1}$
if and only if $y\in\p^m$. Therefore,
$$
 \OF_m=\{x+\varpi_F^my\xi_0:\:x,y\in\OF\}
$$
and
$$
 \OF''_m=\{x+\varpi_F^my\xi_0:\:x,y\in\OF\},
$$
so that actually $\OF_m=\OF''_m$. Hence
$(\OF_L^\times:(\OF''_m)^\times) =(\OF_L^\times:(\OF_m)^\times)$.
By Lemma (3.5.3) of \cite{Fu},
$$
 (\OF_L^\times:(\OF''_m)^\times)=\Big(1-\Big(\frac L\p\Big)q^{-1}\Big)q^m.
$$
This concludes the proof.\qed

\vspace{2ex}
Let us summarize the volume computations.
\begin{lemma}\label{k1volumesprop}
 Let $V^{l,m}_1$ and $V^{l,m}_2$
 be the volumes defined in (\ref{k1voleq1}), (\ref{k1voleq2}).
 \begin{enumerate}
  \item For any $l,m\geq0$ we have
   $$
    V^{l,m}_1=\frac1{(q+1)(q^4-1)}\Big(1-\Big(\frac L\p\Big)q^{-1}\Big)q^{4m+3l+1}.
   $$
  \item For any $l\geq0$ and $m>0$ we have
   $$
    V^{l,m}_2=\frac1{(q+1)(q^4-1)}\Big(1-\Big(\frac L\p\Big)q^{-1}\Big)q^{4m+3l+2}.
   $$
 \end{enumerate}
\end{lemma}
{\bf Proof.} This follows from Lemmas \ref{generalvolumecalclemma2},
\ref{GL2volumecalclemma1} and  \ref{GL2volumecalclemma2}. \qed
\subsection{Main local theorem}
In Sections \ref{W sharp function} and \ref{basic integral}, we
have defined the functions $W^\#$ and the integral $Z$ for any
ramified representation $\tau$ of $\GL_2(F)$. We have computed the
relevant double cosets and their corresponding volumes under the
assumption that $\tau$ has conductor $\p$. We will now assume that
$\tau=\Omega\St_{\GL(2)}$, where $\Omega$ is an unramified
character of $F^\times$, and $\St_{\GL(2)}$ is the Steinberg
representation of $\GL(2,F)$. Then $\tau$ has conductor $\p$, and
the central character of $\tau$ is $\omega_{\tau} = \Omega^2$. We
work in the $\psi^{-c}$ Whittaker model for $\tau$. In this model,
the newform $W^{(0)}$ has the properties
\begin{equation}\label{unramifiedsteinbergkirillov}
 W^{(0)}(\begin{bmatrix}a\\&1\end{bmatrix}) = \left\{
 \begin{array}{ll}
    |a|\Omega(a)& \hbox{ if } a \in \OF, \\
    0 & \hbox{ otherwise,} \\
 \end{array}\right.
\end{equation}
and
\begin{equation}\label{unramifiedsteinbergAL}
 W^{(0)}(g\begin{bmatrix}&1\\
 \varpi\end{bmatrix})=-\Omega(\varpi)W^{(0)}(g)\qquad\text{ for all }g\in\GL_2(F).
\end{equation}
We refer to \cite{Sc} for details.
Using Lemma \ref{relevant disjoint double cosets}, we have
\begin{align}\label{Zcosetapp}
 Z(s)&=\sum_{\substack{l\geq0\\m>0}}\;
  B(h(l,m))\sum_{u \in \OF^{\times}/(1+\p)}\Big(W^\#(\eta
  h(l,m)\begin{bmatrix}1\\&u\\&&1\\&&&u^{-1}\end{bmatrix},s)V_1^{l,m} \nonumber \\
 &\hspace{30ex}+W^\#(\eta h(l,m)\begin{bmatrix}1\\&u\\&&1\\&&&u^{-1}\end{bmatrix}
  s_1s_2s_1,s)V_2^{l,m}\Big) \nonumber \\
 &+\sum_{l \geq 0}\;B(h(l,0))\sum_{u \in
  \OF^{\times}/(1+\p)}W^\#(\eta
  h(l,0)\begin{bmatrix}1\\&u\\&&1\\&&&u^{-1}\end{bmatrix},s)V_1^{l,0}.
\end{align}
Recall formula (\ref{Wsharpformulaeq}) for the function $W^\#(\,\cdot\,,s)$.
Note that for $\zeta\in F^\times$ we
have $\chi(\zeta)=\omega_\pi(\zeta)^{-1}\omega_{\tau}(\zeta)^{-1}
=\omega_\pi(\zeta)^{-1}\Omega(\zeta)^{-2}$. Using
(\ref{nsunramn1matrixidentity1eq}) and
(\ref{nsunramn1matrixidentity3eq}), we write the argument of
$W^\#$ as an element of $M(F)N(F)K^\#(\P)$. Then, from (\ref{Wsharpformulaeq}),
\begin{align}
 Z(s)&=\sum_{\substack{l\geq0\\m>0}}\;B(h(l,m))\sum_{u \in
  \OF^{\times}/(1+\p)}|\varpi^{2m+l}|^{3(s+\frac12)}\omega_\pi(\varpi^{2m+l})^{-1}
   \Omega(\varpi^{2m+l})^{-2}\Omega^2(\varpi^m)\nonumber\\
 &\hspace{20ex}\times\Big(W^{(0)}(\begin{bmatrix}\varpi^{l}\\&1\end{bmatrix})V_1^{l,m}
  +W^{(0)}(\begin{bmatrix}&\varpi^l\\1\end{bmatrix})V_2^{l,m}\Big)\nonumber \\
 &\qquad+\sum_{l \geq 0}\;B(h(l,0))\sum_{u \in\OF^{\times}/(1+\p)}
  |\varpi^{l}|^{3(s+\frac12)}\omega_\pi(\varpi^l)^{-1}
  \Omega(\varpi^{l})^{-2}W^{(0)}(\begin{bmatrix}\varpi^l\\&1\end{bmatrix})V_1^{l,0}.
\end{align}
Here, we have used the fact that $\Omega^2$ is the central
character of $\tau$ and $W^{(0)}$ is right invariant under
$\begin{bmatrix}\OF^{\times}&\OF\\\p&\OF^{\times}\end{bmatrix}$.
It follows from (\ref{unramifiedsteinbergkirillov}) and (\ref{unramifiedsteinbergAL}) that
$W^{(0)}(\begin{bmatrix}&\varpi^l\\1\end{bmatrix})=-\Omega(\varpi^{l})|\varpi|^{l+1}$
for all $l \geq 0$. Hence, we get
\begin{align}
 Z(s) &= (q-1)\sum_{\substack{l\geq0\\m>0}}\;B(h(l,m))
  |\varpi^{2m+l}|^{3(s+\frac12)}\omega_\pi(\varpi^{2m+l})^{-1}
  \Omega(\varpi^{m+l})^{-2}\nonumber\\
 &\hspace{20ex}\times \Big(|\varpi|^l\Omega(\varpi^l)V_1^{l,m}
  - |\varpi|^{l+1}\Omega(\varpi^{l})V_2^{l,m}\Big) \nonumber \\
 & \qquad+ (q-1)\sum\limits_{l \geq 0} B(h(l,0))|\varpi^{l}|^{3(s+\frac12)}
  \omega_\pi(\varpi^l)^{-1}
  \Omega(\varpi^l)^{-2}|\varpi|^l\Omega(\varpi^l)V_1^{l,0}\nonumber\\
 &=(q-1)\sum_{\substack{l\geq0\\m>0}}\;B(h(l,m))
  |\varpi^{2m+l}|^{3(s+\frac12)}|\varpi|^l(\omega_\pi\Omega)(\varpi)^{-2m-l}
  \Big(V_1^{l,m}- q^{-1}V_2^{l,m}\Big) \nonumber \\
 &\qquad + (q-1)\sum\limits_{l \geq 0} B(h(l,0))|\varpi|^{l(3s+5/2)}
  (\omega_\pi\Omega)(\varpi)^{-l}V_1^{l,0}.
\end{align}
By Lemma \ref{k1volumesprop} i) and ii) the first sum is zero. Hence
\begin{align}\label{k1Zresulteq}
 Z(s)  &=(q-1)\sum\limits_{l \geq 0} B(h(l,0))q^{-l(3s+5/2)}
  (\omega_\pi\Omega)(\varpi)^{-l}V_1^{l,0}\nonumber\\
 &=\frac{q-1}{(q+1)(q^4-1)}\sum\limits_{l \geq 0} B(h(l,0))q^{-l(3s+5/2)}
  (\omega_\pi\Omega)(\varpi)^{-l}\Big(1-\Big(\frac L\p\Big)q^{-1}\Big)q^{3l+1}\nonumber\\
 &=\frac{q(q-1)}{(q+1)(q^4-1)}\Big(1-\Big(\frac L\p\Big)q^{-1}\Big)
  \sum\limits_{l \geq 0} B(h(l,0))\big(q^{-3s+1/2}
  (\omega_\pi\Omega)(\varpi)^{-1}\big)^l.
\end{align}
Let $\pi=\chi_1\times\chi_2\rtimes\sigma$ be an unramified
principal series representation of $\GSp_4(F)$; in case
$\chi_1\times\chi_2\rtimes\sigma$ is not irreducible, take its
unramified constituent. Recall the characters $\gamma^{(1)},\ldots,\gamma^{(4)}$ defined
in Sect.\ \ref{sphericalbesselfunctionsec}.
Let $\nu$ be the absolute value in $F$ normalized by $\nu(\varpi)= q^{-1}$. Set
\begin{equation}\label{pi-tau-lfn}
 L(s,\tilde\pi\times\tilde\tau) =
 \prod_{i=1}^4\Big(1-((\gamma^{(i)})^{-1}\Omega^{-1}\nu^{1/2})(\varpi_F)q^{-s}\Big)^{-1}.
\end{equation}
Then $L(s,\tilde\pi\times\tilde\tau)$ is the standard $L$-factor attached to
the representation $\tilde\pi\times\tilde\tau$ of $\GSp_4(F)\times\GL_2(F)$
by the local Langlands correspondence. Here,
$\tilde\pi$ (resp.\ $\tilde\tau$) denotes the contragredient representation
of $\pi$ (resp.\ $\tau$). Denote by $\AI(\Lambda)$ the irreducible, admissible
representation of $\GL_2(F)$ obtained by automorphic induction
from the character $\Lambda$ of $L^\times$. Set
\begin{equation}\label{tau-twisted-lfn}\renewcommand{\arraystretch}{1.3}
 L(s,\tau\times \AI(\Lambda) \times \chi|_{F^{\times}}) = \left\{\begin{array}{ll}
 (1-\chi(\varpi)q^{-1}q^{-2s})^{-1}, & \hbox{ if } \big(\frac{L}{\p}\big) = -1,\\
 (1-\Lambda(\varpi_L)(\chi\Omega)(\varpi)q^{-1/2}q^{-3s-1})^{-1}, & \hbox{ if }
  \big(\frac{L}{\p}\big) = 0,\\
  (1-\Lambda(\varpi_L)(\chi\Omega)(\varpi)q^{-1/2}q^{-3s-1})^{-1} \\
   \times(1-\Lambda(\varpi \varpi_L^{-1})(\chi\Omega)(\varpi)q^{-1/2}q^{-3s-1})^{-1}, &
   \hbox{ if } \big(\frac{L}{\p}\big) = 1.
   \end{array}\right.
\end{equation}
Then $L(s,\tau\times \AI(\Lambda) \times \chi|_{F^{\times}})$ is the standard
$L$-factor attached to the representation
$\tau\times \AI(\Lambda) \times \chi|_{F^{\times}}$ of
$\GL_2(F)\times\GL_2(F)\times\GL_1(F)$ by the local Langlands correspondence.
We now state the main theorem of the local non-archimedean theory.
\begin{theorem}\label{localintegralmaintheorem}
 Let $\pi$ be an unramified, irreducible, admissible representation of $\GSp_4(F)$
 (not necessarily with trivial central character), and let
 $\tau=\Omega\St_{\GL(2)}$ with an unramified character $\Omega$ of $F^\times$.
 Let $Z(s)$ be the integral (\ref{Zvdefeq}), where $W^\#$ is the function
 defined in Sect.\ \ref{W sharp function}, and $B$ is the spherical Bessel
 function defined in Sect.\ \ref{sphericalbesselfunctionsec}. Then
 \begin{equation}\label{localintegralmaintheoremeq}
  Z(s)=\frac{q(q-1)}{(q+1)(q^4-1)}\Big(1-\Big(\frac L\p\Big)q^{-1}\Big)
  \frac{L(3s+\frac 12,\tilde\pi\times\tilde\tau)}
  {L(3s+1, \tau \times \AI(\Lambda)\times\chi|_{F^{\times}})}.
 \end{equation}
\end{theorem}
\nll
{\bf Proof.} By (\ref{suganox0eq}) and (\ref{k1Zresulteq}),
\begin{equation}\label{k1localintegralmaintheoremeq2}
 Z(s)=\frac{q(q-1)}{(q+1)(q^4-1)}\Big(1-\Big(\frac L\p\Big)q^{-1}\Big)\frac{H(q^{-3s+1/2}
 (\omega_\pi\Omega)(\varpi_F)^{-1})}{Q(q^{-3s+1/2}(\omega_\pi\Omega)(\varpi_F)^{-1})}.
\end{equation}
By (\ref{suganoQeq}),
\begin{align}
 Q(q^{-3s+1/2}(\omega_\pi\Omega)(\varpi_F)^{-1}))
 &=\prod_{i=1}^4\big(1-\gamma^{(i)}(\varpi_F)q^{-3s-1}
  (\omega_\pi\Omega)(\varpi_F)^{-1})\big)\nonumber\\
 &=\prod_{i=1}^4\big(1-(\gamma^{(i)}(\omega_\pi\Omega)^{-1}\nu^{1/2})(\varpi_F)
  q^{-3s-1/2}\big)\nonumber\\
 &=\prod_{i=1}^4\big(1-((\gamma^{(i)})^{-1}\Omega^{-1}\nu^{1/2})(\varpi_F)
  q^{-3s-1/2}\big)\nonumber\\
 &\stackrel{(\ref{pi-tau-lfn})}{=}L(3s+1/2,\tilde\pi\times\tilde\tau)^{-1}.\nonumber
\end{align}
To compute the numerator of (\ref{k1localintegralmaintheoremeq2}), we distinguish
cases. If $\big(\frac L\p\big)=-1$, then $H(y)=1-q^{-4}\Lambda(\varpi_F)y^2$, and
hence
\begin{align*}
 H(q^{-3s+1/2}(\omega_\pi\Omega)(\varpi_F)^{-1})
 &=1-q^{-4}\Lambda(\varpi_F)\big(q^{-3s+1/2}(\omega_\pi\Omega)(\varpi_F)^{-1}\big)^2\\
 &=1-\big(\Lambda\omega_\pi^{-2}\Omega^{-2}\big)(\varpi_F)q^{-6s-3}\\
 &=1-\big(\omega_\pi^{-1}\omega_\tau^{-1}\big)(\varpi_F)q^{-6s-3}\\
 &=1-\chi(\varpi_F)q^{-1}q^{-6s-2}\\
 &\stackrel{(\ref{tau-twisted-lfn})}{=}L(3s+1, \tau \times \AI(\Lambda)\times \chi|_{F^{\times}})^{-1}.
\end{align*}
If $\big(\frac L\p\big)=0$, then $H(y)=1-q^{-2}\Lambda(\varpi_L)y$, and hence
\begin{align*}
 H(q^{-3s+1/2}(\omega_\pi\Omega)(\varpi_F)^{-1})
 &=1-q^{-2}\Lambda(\varpi_L)q^{-3s+1/2}(\omega_\pi\Omega)(\varpi_F)^{-1}\\
 &=1-\Lambda(\varpi_L)(\omega_\pi\omega_\tau\Omega^{-1})(\varpi_F)^{-1}q^{-3s-3/2}\\
 &=1-\Lambda(\varpi_L)(\chi\Omega)(\varpi_F)q^{-1/2}q^{-3s-1}\\
 &\stackrel{(\ref{tau-twisted-lfn})}{=}
  L(3s+1, \tau \times \AI(\Lambda)\times \chi|_{F^{\times}})^{-1}.
\end{align*}
If $\big(\frac L\p\big)=1$, then
$H(y)=(1-q^{-2}\Lambda(\varpi_L)y)(1-q^{-2}\Lambda(\varpi_F\varpi_L^{-1})y)$, and hence
\begin{align*}
 H(q^{-3s+1/2}(\omega_\pi\Omega)(\varpi_F)^{-1})
 &=(1-q^{-2}\Lambda(\varpi_L)q^{-3s+1/2}(\omega_\pi\Omega)(\varpi_F)^{-1})\\
  &\hspace{10ex}(1-q^{-2}\Lambda(\varpi_F\varpi_L^{-1})
   q^{-3s+1/2}(\omega_\pi\Omega)(\varpi_F)^{-1})\\
 &=(1-\Lambda(\varpi_L)(\omega_\pi\omega_\tau\Omega^{-1})(\varpi_F)^{-1}q^{-3s-3/2})\\
  &\hspace{10ex}(1-\Lambda(\varpi_F\varpi_L^{-1})
   (\omega_\pi\omega_\tau\Omega^{-1})(\varpi_F)^{-1}q^{-3s-3/2})\\
 &=(1-\Lambda(\varpi_L)(\chi\Omega)(\varpi_F)q^{-1/2}q^{-3s-1})\\
  &\hspace{10ex}(1-\Lambda(\varpi_F\varpi_L^{-1})
   (\chi\Omega)(\varpi_F)q^{-1/2}q^{-3s-1})\\
 &\stackrel{(\ref{tau-twisted-lfn})}{=}
  L(3s+1, \tau \times \AI(\Lambda)\times \chi|_{F^{\times}})^{-1}.
\end{align*}
Hence $H(q^{-3s+1/2}(\omega_\pi\Omega)(\varpi_F)^{-1}) =L(3s+1,
\tau \times \AI(\Lambda)\times \chi|_{F^{\times}})^{-1}$ in all
cases. This concludes the proof of the theorem. \qed
\section{Local archimedean theory}\label{arch-sec}
In this section we compute the local archimedean integral. As in
Sect.\ \ref{non-arch-ram-sec}, the key step is the choice of
vectors $B$ and $W^\#$.
\subsection{Real groups}\label{real groups}
Consider the symmetric domains $\HH_2 := \{Z \in M_2(\C):\;
i(\,^{t}\!\bar{Z}-Z) \mbox{ is positive definite}\}$ and
$\SH_2 := \{Z \in \HH_2:\;^{t}Z = Z \}$. The group
$G^{+}(\R) := \{g \in G(\R) : \mu_2(g) > 0 \}$ acts on $\HH_2$
via $(g,Z) \mapsto g\langle Z \rangle$, where
$$
 g\langle Z \rangle = (AZ+B)(CZ+D)^{-1}, \mbox{ for } g =
 \mat{A}{B}{C}{D} \in G^{+}(\R), Z \in \HH_2.
$$
Under this action, $\SH_2$ is stable by $H^+(\R) = \GSp_4^+(\R)$.
The group $K_{\infty} = \{g\in G^{+}(\R) : \mu_2(g) = 1, g\langle I\rangle = I \}$
is a maximal compact subgroup of $G^{+}(\R)$. Here, $I=\mat{i}{}{}{i} \in\HH_2$.
Explicitly,
$$
 K_\infty=\{\mat{A}{B}{-B}{A}:\;A,B\in\GL(2,\C),\;^t\!\bar AB=\,^t\bar BA,\;
 ^t\!\bar AA+\,^t\bar BB=1\}.
$$
By the Iwasawa decomposition
\begin{equation}\label{realiwasawaeq}
 G(\R)=M^{(1)}(\R)M^{(2)}(\R)N(\R)K_\infty,
\end{equation}
where $M^{(1)}(\R)$, $M^{(2)}(\R)$ and $N(\R)$ are as defined in (\ref{M1defn}),
(\ref{M2defn}), (\ref{Ndefn}). A calculation shows that
\begin{align}\label{realKintersect1eq}
 &M^{(1)}(\R)M^{(2)}(\R)N(\R)\cap K_\infty\nonumber\\
 &\qquad=\{\begin{bmatrix}\zeta\\&\alpha&&\beta\\&&\zeta\\&-\beta&&\alpha\end{bmatrix}:
 \:\zeta,\alpha,\beta\in\C,\:|\zeta|=1,\:|\alpha|^2+|\beta|^2=1,\:\alpha\bar\beta=\beta\bar\alpha\}.
\end{align}
Note also that
\begin{equation}\label{realKintersect2eq}
 M^{(2)}(\R)\cap K_\infty
 =\{\begin{bmatrix}1\\&\alpha&&\beta\\&&1\\&-\beta&&\alpha\end{bmatrix}:
 \:\alpha,\beta\in\C,\:|\alpha|^2+|\beta|^2=1,\:\alpha\bar\beta=\beta\bar\alpha\},
 \end{equation}
and that there is an isomorphism
\begin{align}\label{realKintersect3eq}
 (S^1\times\SO(2))/\{(\lambda,\mat{\lambda}{}{}{\lambda}):\:\lambda=\pm1\}&\stackrel{\sim}{\longrightarrow}
  M^{(2)}(\R)\cap K_\infty,\nonumber\\
 (\lambda,\mat{\alpha}{\beta}{-\beta\:}{\alpha})&\longmapsto
  \begin{bmatrix}1\\&\lambda\alpha&&\lambda\beta\\&&1\\&-\lambda\beta&&\lambda\alpha\end{bmatrix}.
\end{align}
For $g\in G^+(\R)$ and $Z\in\mathbb{H}_2$, let $J(g,Z)=CZ+D$ be
the automorphy factor. Then, for any integer $l$, the map
\begin{equation}\label{realcompactcharactereq}
 k\longmapsto \det(J(k,I))^l
\end{equation}
defines a character $K_\infty\rightarrow\C^\times$. If
$k\in M^{(2)}(\R)\cap K_\infty$ is written in the form (\ref{realKintersect3eq}), then
$\det(J(k,I))^l=\lambda^le^{-il\theta}$, where $\alpha=\cos(\theta)$, $\beta=\sin(\theta)$.
Let $K^H_\infty = K_{\infty} \cap H^+(\R)$. Then $K^H_\infty$ is a maximal compact subgroup,
explicitly given by
$$
 K^H_\infty=\{\mat{A}{B}{-B}{A}:\;^t\!AB=\,^tBA,\;^t\!AA+\,^tBB=1\}.
$$
Sending $\mat{A}{B}{-B}{A}$ to $A-iB$ gives an isomorphism $K^H_\infty\cong{\rm U}(2)$.
Recall that we have chosen $a,b,c\in\R$ such that $d=b^2-4ac\neq0$. In the archimedean
case we shall assume that $d<0$ and let $D=-d$. Then $\R(\sqrt{-D})=\C$.
As in Sect.\ \ref{besselsubgroupsec} we have
\begin{equation}\label{TReq}
 T(\R)=\{\mat{x+yb/2}{yc}{-ya}{x-yb/2}:\:x,y\in\R,\:x^2+y^2D/4>0\}.
\end{equation}
Let
\begin{equation}\label{Tinfty1eq}
 T^1_\infty=T(\R)\cap\SL(2,\R)
 =\{\mat{x+yb/2}{yc}{-ya}{x-yb/2}:\:x,y\in\R,\:x^2+y^2D/4=1\}.
\end{equation}
We have $T(\R)\cong\C^\times$ via $\mat{x+yb/2}{yc}{-ya}{x-yb/2}\mapsto x+y\sqrt{-D}/2$.
Under this isomorphism $T^1_\infty$ corresponds to the unit circle. We have
\begin{equation}\label{TRTinfty1eq}
 T(\R)=T^1_\infty\cdot\{\mat{\zeta}{}{}{\zeta}:\:\zeta>0\}.
\end{equation}
As in \cite{Fu}, p.\ 211, let $t_0\in\GL_2(\R)^+$ be such that
$T_\infty^1=t_0\SO(2)t_0^{-1}$. We will make a specific choice of
$t_0$ when we choose the matrix $S=\mat{a}{b/2}{b/2}{c}$ below. By
the Cartan decomposition,
\begin{equation}\label{cartan1eq}
 \GL_2^+(\R)=\SO(2)\cdot\{\mat{\zeta_1}{}{}{\zeta_2}:\:\zeta_1,\zeta_2>0,\:\zeta_1\geq \zeta_2\}\cdot\SO(2).
\end{equation}
Therefore,
\begin{align}\label{GL2RplusBesseldecompeq}
 \GL_2^+(\R)&=t_0\SO(2)\cdot\{\mat{\zeta_1}{}{}{\zeta_2}:\:\zeta_1,\zeta_2>0,\:\zeta_1\geq \zeta_2\}\cdot\SO(2)\nonumber\\
 &=T^1_\infty t_0\cdot\{\mat{\sqrt{\zeta_1\zeta_2}}{}{}{\sqrt{\zeta_1\zeta_2}}\mat{\sqrt{\zeta_1/\zeta_2}}{}{}{\sqrt{\zeta_2/\zeta_1}}
  :\:\zeta_1,\zeta_2>0,\:\zeta_1\geq \zeta_2\}\cdot\SO(2)\nonumber\\
 &=T(\R)t_0\cdot\{\mat{\zeta}{}{}{\zeta^{-1}}:\:\zeta\geq1\}\cdot\SO(2).
\end{align}
Using this, it is not hard to see that
\begin{equation}\label{HRBesseldecompeq}
 H(\R)=R(\R)\cdot\big\{\begin{bmatrix}\lambda t_0\mat{\zeta}{}{}{\zeta^{-1}}&\\
  &^tt_0^{-1}\mat{\zeta^{-1}}{}{}{\zeta}\end{bmatrix}:\:\lambda\in\R^\times,\,\zeta\geq1\big\}
  \cdot K^H_\infty.
\end{equation}
Here, $R(\R)=T(\R)U(\R)$ is the Bessel subgroup defined in Sect.\ \ref{besselsubgroupsec}.
One can check that all the double cosets in (\ref{HRBesseldecompeq}) are disjoint.
\subsection{The Bessel function}\label{realbesselfnsec}
Recall that we have chosen three elements $a,b,c\in\R$ such that $d=b^2-4ac\neq0$.
We will now make the stronger assumption that
$S=\mat{a}{b/2}{b/2}{c}\in M_2(\R)$ is a positive definite matrix.
Set $D=4ac-b^2>0$, as above. Given a positive integer $l\geq2$,
consider the function $B:\:H(\R)\rightarrow\C$ defined by
\begin{equation}\label{archBesselformula2eq}\renewcommand{\arraystretch}{1.2}
 B(h) := \left\{\begin{array}{ll}
 \mu_2(h)^l\,\overline{\det(J(h, I))^{-l}}\,e^{-2\pi i\,
 {\rm tr}(S\overline{h\langle I\rangle})}& \hbox{ if } h \in H^+(\R),\\
 0& \hbox{ if } h \notin H^+(\R),\\
\end{array}\right.
\end{equation}
where $I=\mat{i}{}{}{i}$. Note that the function $B$ only depends
on the choice of $S$ and $l$. Recall the character $\theta$ of
$U(\R)$ defined in (\ref{thetadef2eq}). It depends on the choice
of additive character $\psi$, and throughout we choose
$\psi(x)=e^{-2\pi ix}$. Then the function $B$ satisfies
\begin{equation}\label{realBproperty1eq}
 B(tuh)=\theta(u)B(h)\qquad\text{for }
 h\in H(\R),\;t\in T(\R),\;u\in U(\R),
\end{equation}
and
\begin{equation}\label{realBproperty2eq}
 B(hk)=\det(J(k,I))^lB(h)\qquad\text{for }
 h\in H(\R),\;k\in K^H_\infty.
\end{equation}
Property (\ref{realBproperty1eq}) means that
$B$ satisfies the Bessel transformation property with the
character $\Lambda \otimes \theta$ of $R(\R)$, where $\Lambda$ is
trivial. In fact, by the considerations in \cite{Su} 1-3, $B$ is
the lowest weight vector in a holomorphic discrete series
representation of ${\rm PGSp}(4,\R)$ corresponding to Siegel
modular forms of degree $2$ and weight $l$. By
(\ref{realBproperty1eq}) and (\ref{realBproperty2eq}), the
function $B$ is determined by its values on a set of
representatives for $R(\R)\backslash H(\R)/K^H_\infty$. Such a set
is given in (\ref{HRBesseldecompeq}).
\subsection{The function $W^\#$}\label{realW sharp function}
Let $(\tau,V_\tau)$ be a  generic, irreducible, admissible
representation of $\GL_2(\R)$ with central character $\omega_{\tau}$. We assume that
$V_\tau=\mathcal{W}(\tau,\psi_{-c})$ is the Whittaker model of
$\tau$ with respect to the non-trivial additive character
$x\mapsto\psi(-cx)$. Note that $S$ positive definite implies $c>0$. Let
$W^{(0)}\in V_\tau$ have weight $l_1$. Then $W^{(0)}$ has the
following properties.
\begin{enumerate}
 \item
  $$
   W^{(0)}(gr(\theta))=e^{il_1\theta}W^{(0)}(g)\qquad\text{for }g\in \GL_2(\R),\;
    r(\theta)=\mat{\cos(\theta)}{\sin(\theta)}{-\sin(\theta)}{\cos(\theta)}\in\SO(2).
  $$
 \item
  $$
   W^{(0)}(\mat{1}{x}{}{1}g)=\psi(-cx)W^{(0)}(g)\qquad\text{for }g\in \GL_2(\R),\;x\in\R.
  $$
\end{enumerate}
Let $\chi_0$ be the character of $\C^\times$ with the properties
\begin{equation}\label{realchi0defeq}
 \chi_0\big|_{\R^\times}=\omega_\tau,\qquad\chi_0(\zeta)=\zeta^{-l_1}\quad\text{for }
 \zeta\in\C^\times,\:|\zeta|=1.
\end{equation}
Such a character exists since $\omega_\tau(-1)=(-1)^{l_1}$.
We extend $W^{(0)}$ to a function on $M^{(2)}(\R)$ via
\begin{equation}\label{realW1extensioneq}
 W^{(0)}(\zeta g)=\chi_0(\zeta)W^{(0)}(g),\qquad \zeta\in\C^\times,\:g\in\GL_2(\R)
\end{equation}
(see Lemma \ref{M2structurelemma}).
Then it is easy to check that
\begin{equation}\label{realW1extensioninvarianceeq}
 W^{(0)}(gk)=\det(J(k,I))^{-l_1}W^{(0)}(g)\qquad\text{for }g\in M^{(2)}(\R)\text{ and }
 k\in M^{(2)}(\R)\cap K_\infty.
\end{equation}
We will need values of $W^{(0)}$ evaluated at $\mat{t}{}{}{1}$ for $t \neq 0$. For this
we look at the Lie algebra $\mathfrak{g}=\mathfrak{gl}(2,\R)$ and consider the elements
$$
 R=\mat{0}{1}{0}{0},\qquad L=\mat{0}{0}{1}{0},\qquad
 H=\mat{1}{0}{0}{-1},\qquad Z=\mat{1}{0}{0}{1}.
$$
In the universal enveloping algebra $U(\mathfrak{g})$ let
\begin{equation}\label{Deltadefeq}
 \Delta=\frac14(H^2+2RL+2LR).
\end{equation}
Then $\Delta$ lies in the center of $U(\mathfrak{g})$ and acts on $V_\tau$ by
a scalar, which we write in the form $-(\frac14+(\frac r2)^2)$ with $r\in\C$. In particular,
\begin{equation}\label{W1Deltapropertyeq}
 \Delta W^{(0)}=-\Big(\frac14+\Big(\frac r2\Big)^2\Big)W^{(0)}.
\end{equation}
If one restricts the function $W^{(0)}$ to $\mat{t^{1/2}}{}{}{t^{-1/2}}$,
$t>0$, then (\ref{W1Deltapropertyeq}) reduces to the differential equation
satisfied by the classical Whittaker functions. Hence,
there exist constants $a^+,a^-\in\C$ such that
\begin{equation}\label{W1Whittakerformulaeq}
 W^{(0)}(\mat{t}{0}{0}{1}) = \renewcommand{\arraystretch}{1.3}
  \left\{\begin{array}{ll}
 a^+\omega_\tau((4\pi ct)^{1/2})W_{\frac{l_1}2, \frac{ir}2}(4\pi ct) & \hbox{ if } t > 0, \\
 a^-\omega_\tau((-4\pi ct)^{1/2})
  W_{-\frac{l_1}2, \frac{ir}2}(-4\pi ct)& \hbox{ if } t <0 .\end{array}\right.
\end{equation}
Here,
$W_{\pm\frac{l_1}2, \frac{ir}2}$ denotes a classical
Whittaker function; see \cite[p. 244]{Bu}, \cite{MOS}.
Let $\chi$ be the character of $\C^\times$ given by
\begin{equation}\label{archchidefeq}
 \chi(\zeta)=\chi_0(\bar{\zeta})^{-1}.
\end{equation}
We interpret $\chi$ as a character of $M^{(1)}(\R)$; see (\ref{M1defn}).
Given a complex number $s$, we define
a function $W^\#(\,\cdot\,,s):\:G(\R)\rightarrow\C$ as follows.
Given $g\in G(\R)$, write $g=m_1m_2nk$ according to (\ref{realiwasawaeq}).
Then set
\begin{equation}\label{realWsharpdefeq}
 W^\#(g,s)=\delta_P^{s+1/2}(m_1m_2)\det(J(k,I))^{-{l_1}}\chi(m_1)W^{(0)}(m_2).
\end{equation}
Property (\ref{realW1extensioninvarianceeq}) shows that this is well-defined.
Explicitly, for $\zeta\in\C^\times$ and $\mat{\alpha}{\beta}{\gamma}{\delta}\in M^{(2)}(\R)$,
\begin{equation}\label{realWsharpproperty1eq}
   W^\#(\begin{bmatrix}\zeta\\&1\\&&\bar{\zeta}^{-1}\\&&&1\end{bmatrix}
   \begin{bmatrix}1\\&\alpha&&\beta\\&&\mu\\&\gamma&&\delta\end{bmatrix},s)
   =||\zeta|^2\cdot\mu^{-1}|^{3(s+1/2)}\chi(\zeta)\,
   W^{(0)}(\mat{\alpha}{\beta}{\gamma}{\delta}).
\end{equation}
Here $\mu=\bar\alpha\delta-\beta\bar\gamma$. It is clear that $W^\#(\,\cdot\,,s)$
satisfies
\begin{equation}\label{realWsharpproperty2eq}
 W^\#(gk,s)=\det(J(k,I))^{-l_1}W^\#(g,s)\qquad\text{for }
 g\in G(\R),\;k\in K_\infty.
\end{equation}
By Lemma \ref{WsharpinvBessellemma}, we have
\begin{equation}\label{realWsharpeq}
  W^\#(\eta tuh,s)=\theta(u)^{-1}W^\#(\eta h,s)
\end{equation}
for $t\in T(\R)$, $u\in U(\R)$, $h\in G(\R)$ and
$$
  \eta=\begin{bmatrix}1\\\alpha&1\\&&1&-\bar\alpha\\&&&1\end{bmatrix},
  \qquad \alpha=\frac{b+\sqrt{d}}{2c},\;d=b^2-4ac.
$$
\subsection{The local archimedean integral}
Let $B$ and $W^\#$ be as defined in Sections \ref{realbesselfnsec} and \ref{realW sharp function}. By (\ref{realBproperty1eq}) and (\ref{realWsharpeq}),
it makes sense to consider the integral
\begin{equation}\label{realZinftyeq}
 Z_{\infty}(s) = \int\limits_{R(\R)\backslash H(\R)}W^\#(\eta h, s) B(h)dh.
\end{equation}
Our goal in the following is to evaluate this integral. It follows from
(\ref{realBproperty2eq}) and (\ref{realWsharpproperty2eq}) that this
integral is zero if $l_1\neq l$. We shall therefore
assume that \underline{$l_1=l$}. Then the function $W^\#(\eta h, s) B(h)$
is right invariant under $K^H_\infty$.
From the disjoint double coset decomposition (\ref{HRBesseldecompeq}) and the fact
that $W^\#(\eta h, s) B(h)$ is right invariant under $K^H_\infty$ we obtain
\begin{align}\label{archintegral1eq}
 Z_{\infty}(s) &= \pi\int\limits_{\R^{\times}}\int\limits_1^{\infty}
  W^\#\Big(\eta \begin{bmatrix}\lambda t_0\mat{\zeta}{}{}{\zeta^{-1}}&\\
  &^tt_0^{-1}\mat{\zeta^{-1}}{}{}{\zeta}\end{bmatrix},s\Big) \nonumber \\
 &\hspace{10ex}B\Big(\begin{bmatrix}\lambda t_0\mat{\zeta}{}{}{\zeta^{-1}}&\\
  &^tt_0^{-1}\mat{\zeta^{-1}}{}{}{\zeta}\,\end{bmatrix}\Big)(\zeta-\zeta^{-3})\lambda^{-4}\,d\zeta\,d\lambda;
\end{align}
see \cite{Fu} (4.6) for the relevant integration formulas.
The above calculations are valid for any choice of $a,b,c$ as long
as $S=\mat{a}{b/2}{b/2}{c}$ is positive definite. To compute
(\ref{archintegral1eq}), we will fix $D = 4ac-b^2$ and make
special choices for $a,b,c$. First assume that \underline{$D\equiv 0 \pmod{4}$}.
In this case let $S(-D) := \mat{\frac D4}{0}{0}{1}$.
Then $\eta =\begin{bmatrix}1&&&\\\frac{\sqrt{-D}}{2}&1&&\\&&1&
\frac{\sqrt{-D}}{2}\\ &&&1\end{bmatrix}$, and we can choose
$t_0 =\mat{2^{1/2}D^{-1/4}}{}{}{2^{-1/2}D^{1/4}}$. From
(\ref{archBesselformula2eq}) we have
\begin{equation}\label{archbesselformula3eq}
 B\Big(\begin{bmatrix}\lambda t_0\mat{\zeta}{}{}{\zeta^{-1}}&\\
  &^tt_0^{-1}\mat{\zeta^{-1}}{}{}{\zeta}\,\end{bmatrix}\Big) =
  \renewcommand{\arraystretch}{1.3}\left\{\begin{array}{ll}\lambda^le^{-2\pi \lambda
  D^{1/2}\frac{\zeta^2+\zeta^{-2}}2}& \hbox{ if } \lambda > 0,\\
    0& \hbox{ if } \lambda < 0.
 \end{array}\right.
\end{equation}
Next we rewrite the argument of $W^\#$ as an element of $MNK_{\infty}$,
\begin{align*}
 &\eta \begin{bmatrix}\lambda t_0\mat{\zeta}{}{}{\zeta^{-1}}&\\
  &^tt_0^{-1}\mat{\zeta^{-1}}{}{}{\zeta}\end{bmatrix}\\
 &\qquad=\begin{bmatrix}\lambda \mat{D^{-\frac14}
  \Big(\frac{\zeta^2+\zeta^{-2}}{2}\Big)^{-\frac 12}}{}{}
  {D^{\frac 14}\Big(\frac{\zeta^2+\zeta^{-2}}{2}\Big)^{\frac 12}}&\\
  &\mat{D^{\frac 14}\Big(\frac{\zeta^2+\zeta^{-2}}{2}\Big)^{\frac 12}}{}{}{D^{-\frac14}
  \Big(\frac{\zeta^2+\zeta^{-2}}{2}\Big)^{-\frac12}}
 \end{bmatrix} \\
 &\qquad\times\begin{bmatrix}1&-i\zeta^2&&\\0&1&&\\&&1&0\\&&-i\zeta^2&1\end{bmatrix}
 \mat{k_0}{0}{0}{k_0},
\end{align*}
where $k_0 \in \SU(2) = \{g \in \SL_2(\C):\: {}^{t}\bar{g} g = I_2\}$.
Hence, using (\ref{realWsharpproperty1eq}) and (\ref{realWsharpproperty2eq}), we get
\begin{align}\label{whittakerfirstformulaeq}
 &W^\#\Big(\eta \begin{bmatrix}\lambda t_0\mat{\zeta}{}{}{\zeta^{-1}}&\\
  &t_0^{-1}\mat{\zeta^{-1}}{}{}{\zeta}\end{bmatrix},s\Big) \nonumber \\
 &\qquad=\Big|\lambda D^{-\frac12}
  \big(\frac{\zeta^2+\zeta^{-2}}{2}\big)^{-1}\Big|^{3(s+\frac12)}
  \omega_\tau(\lambda)^{-1}
  W^{(0)}(\mat{\lambda D^{\frac12}\big(\frac{\zeta^2+\zeta^{-2}}{2}\big)}{0}{0}{1}).
\end{align}
Let $q\in\C$ be such that $\omega_\tau(y)=y^q$ for $y>0$.
It follows from (\ref{W1Whittakerformulaeq}), (\ref{archbesselformula3eq}) and
(\ref{whittakerfirstformulaeq}) that
\begin{align}\label{realI1eq}
 Z_{\infty}(s)&=a^+\pi D^{-\frac{3s}2-\frac 34+\frac q4}(4\pi)^{\frac q2}
  \int\limits_0^{\infty}\int\limits_1^{\infty}
  \lambda^{3s+\frac 32+l-\frac q2}
  \Big(\frac{\zeta^2+\zeta^{-2}}2\Big)^{-3s-\frac 32+\frac q2}
  W_{\frac l2, \frac{ir}2}\big(4 \pi \lambda
  D^{1/2}\frac{\zeta^2+\zeta^{-2}}2\big) \nonumber\\
 &\hspace{30ex}e^{-2\pi \lambda
  D^{1/2}\frac{\zeta^2+\zeta^{-2}}2}(\zeta-\zeta^{-3})\lambda^{-4}\,d\zeta\,d\lambda.
\end{align}
Substituting $u = (\zeta^2+\zeta^{-2})/2$ we get
$$
 Z_{\infty}(s) =a^+ \pi D^{-\frac{3s}2-\frac 34+\frac q4}(4\pi)^{\frac q2}
 \int\limits_1^{\infty}\int\limits_0^{\infty} \lambda^{3s-\frac 32+l-\frac q2}
 u^{-3s-\frac 32+\frac q2}
 W_{\frac l2, \frac{ir}2}(4 \pi \lambda
 D^{1/2}u) e^{-2\pi \lambda D^{1/2}u}\,\frac{d\lambda}{\lambda}\,du.
$$
We will first compute the integral with respect to $\lambda$. For
a fixed $u$ substitute $x = 4 \pi \lambda D^{1/2}u$ to get
$$
 Z_{\infty}(s) =a^+\pi D^{-3s-\frac l2+\frac q2}(4 \pi)^{-3s+\frac 32 -l+q}
 \int\limits_1^{\infty}u^{-6s-l+q}\int\limits_0^{\infty}W_{\frac l2,
 \frac{ir}2}(x)e^{-\frac x2}x^{3s-\frac 32 + l-\frac q2}\,\frac{dx}x\,du.
$$
Using the integral formula for the Whittaker function from \cite[p. 316]{MOS}, we get
\begin{align}\label{I1integraleq}
 Z_{\infty}(s) &=a^+\pi D^{-3s-\frac l2+\frac q2}(4\pi)^{-3s+\frac 32 -l+q}\,
  \frac{\Gamma(3s+l-1+\frac{ir}2-\frac q2)
  \Gamma(3s+l-1-\frac{ir}2-\frac q2)}{\Gamma(3s+\frac l2 - \frac 12-\frac q2)}
  \int\limits_1^{\infty}u^{-6s-l+q} du \nonumber \\
 &=a^+\pi D^{-3s-\frac l2+\frac q2}\,\frac{(4\pi)^{-3s+\frac 32 -l+q}}{6s+l-q-1}\,
  \frac{\Gamma(3s+l-1+\frac{ir}2-\frac q2)
  \Gamma(3s+l-1-\frac{ir}2-\frac q2)}{\Gamma(3s+\frac l2 - \frac 12-\frac q2)}\nonumber\\
 &=\frac{a^+}2\pi D^{-3s-\frac l2+\frac q2}\,(4\pi)^{-3s+\frac 32 -l+q}\,
  \frac{\Gamma(3s+l-1+\frac{ir}2-\frac q2)
  \Gamma(3s+l-1-\frac{ir}2-\frac q2)}{\Gamma(3s+\frac{l+1-q}2)}.
\end{align}
Here, for the calculation of the $u$-integral, we have assumed that ${\rm Re}(6s+l-q)>0$. ---
Now assume that \underline{$D \equiv 3 \pmod{4}$}. In this case we choose
$$\renewcommand{\arraystretch}{1.2}
 S(-D) = \mat{\frac{1+D}4}{\frac 12}{\frac 12}{1} =
 \mat{1}{\frac 12}{0}{1} \mat{\frac D4}{0}{0}{1} \mat{1}{0}{\frac 12}{1}.
$$
Let $T(\R)$, $R(\R)$, $\eta$, $B$ be the objects defined with this
$\mat{a}{b/2}{b/2}{c}=\mat{\frac{1+D}4}{\frac 12}{\frac 12}{1}$, and let
$\tilde T(\R)$, $\tilde R(\R)$, $\tilde\eta$, $\tilde B$
be the corresponding objects defined with
$\mat{\tilde a}{\tilde b/2}{\tilde b/2}{\tilde c}=\mat{\frac D4}{}{}{1}$. Let
$$
 h_0=\begin{bmatrix}1\\-\frac12&1\\&&1&\frac12\\&&&1\end{bmatrix}.
$$
Then
$$
 T^1(\R)=h_0\tilde T^1(\R)h_0^{-1},\qquad
 T(\R)=h_0\tilde T(\R)h_0^{-1},\qquad R(\R)=h_0\tilde R(\R)h_0^{-1}.
$$
Furthermore, $\eta=\tilde\eta h_0^{-1}$. The integral (\ref{realZinftyeq}) becomes
\begin{align*}
 Z_{\infty}(s)&= \int\limits_{R(\R)\backslash H(\R)}W^\#(\eta h, s) B(h)\,dh\\
 &= \int\limits_{h_0\tilde R(\R)h_0^{-1}\backslash H(\R)}
  W^\#(\tilde\eta h_0^{-1}h, s) B(h_0h_0^{-1}h)\,dh\\
 &= \int\limits_{h_0\tilde R(\R)h_0^{-1}\backslash H(\R)}
  W^\#(\tilde\eta h_0^{-1}hh_0, s) B(h_0h_0^{-1}hh_0)\,dh\\
 &= \int\limits_{\tilde R(\R)\backslash H(\R)}
  W^\#(\tilde\eta h, s) B(h_0h)\,dh\\
 &= \int\limits_{\tilde R(\R)\backslash H(\R)}
  W^\#(\tilde\eta h, s) \tilde B(h)\,dh.
\end{align*}
This integral can be computed just like the one in the case
$D\equiv0$ mod $4$, and we get the exactly same answer as in
(\ref{I1integraleq}).
\begin{theorem}\label{archmaintheorem}
 Let $l$ and $D$ be positive integers such that $D\equiv0,3$ mod $4$.
 Let $S(-D)=\mat{D/4}{}{}{1}$ if $D\equiv0$ mod $4$ and $S(-D)=\mat{(1+D)/4}{1/2}{1/2}{1}$
 if $D\equiv3$ mod $4$. Let $B:\:\GSp_4(\R)\rightarrow\C$ be the function defined
 in (\ref{archBesselformula2eq}), and let $W^\#(\,\cdot\,,s)$ be the function
 defined in (\ref{realWsharpdefeq}). Then, for ${\rm Re}(6s+l-q)>0$,
 \begin{align}\label{archmaintheoremeq1}
  \nonumber Z_{\infty}(s)&:= \int\limits_{R(\R)\backslash H(\R)}W^\#(\eta h, s) B(h)\,dh\\
  &=\frac{a^+}2\pi D^{-3s-\frac l2+\frac q2}\,(4\pi)^{-3s+\frac 32 -l+q}\,
  \frac{\Gamma(3s+l-1+\frac{ir}2-\frac q2)
  \Gamma(3s+l-1-\frac{ir}2-\frac q2)}{\Gamma(3s+\frac{l+1-q}2)}.
 \end{align}
 Here, $q\in\C$ is related to the central character of $\tau$ via
 $\omega_\tau(y)=y^q$ for $y>0$. The number $r\in\C$ is such
 that (\ref{W1Deltapropertyeq}) holds.
\end{theorem}
We will state two special cases of formula (\ref{archmaintheoremeq1}). First
assume that $\tau=\chi_1\times\chi_2$, an irreducible principal series representation of
$\GL(2,\R)$, where $\chi_1$ and $\chi_2$ are characters of $\R^\times$. Let
$\varepsilon_i\in\{0,1\}$ and $s_i\in\C$ be such that
$\chi_i(x)={\rm sgn}(x)^{\varepsilon_i}|x|^{s_i}$, for $i=1,2$. Then $\Delta$ acts
on $\tau$ by multiplication with $-\frac14(1-(s_1-s_2)^2)$.
Comparing with (\ref{W1Deltapropertyeq}), we get $(s_1-s_2)^2=-r^2$, so that
$ir=\pm(s_1-s_2)$. Furthermore, $q=s_1+s_2$. Therefore,
\begin{equation}\label{archmaintheoremeq2}
 Z_\infty(s)=\frac{a^+}2\pi D^{-3s-\frac l2+\frac{s_1+s_2}2}\,
  (4\pi)^{-3s+\frac 32 -l+s_1+s_2}\,
  \frac{\Gamma(3s+l-1-s_1)
  \Gamma(3s+l-1-s_2)}{\Gamma(3s+\frac{l+1-s_1-s_2}2)}.
\end{equation}
Now assume that $l_1$ is a positive integer, that $q\in\C$, and
that $\tau=\mathcal{D}_q(l_1)$, the discrete series (or limit of
discrete series) representation of $\GL(2,\R)$ with a lowest
weight vector of weight $l_1$ for which the central element
$Z=\mat{1}{}{}{1}$ acts by multiplication with $q$. Then
$ir=\pm(l_1-1)$, so that, from (\ref{archmaintheoremeq1}),
\begin{equation}\label{archmaintheoremeq3}
  Z_{\infty}(s)=\frac{a^+}2\pi D^{-3s-\frac l2+\frac q2}\,(4\pi)^{-3s+\frac 32 -l+q}\,
  \frac{\Gamma(3s+l-1+\frac{l_1-1}2-\frac q2)
  \Gamma(3s+l-1-\frac{l_1-1}2-\frac q2)}{\Gamma(3s+\frac{l+1-q}2)}.
\end{equation}
\section{Modular Forms}\label{modular forms sec}
Let $\A$ be the ring of adeles of $\Q$.
In this section we will consider a cuspidal, automorphic
representation $\pi$ of $\GSp_4(\A)$, obtained from a Siegel cusp
form, and a cuspidal, automorphic representation
$\tau$ of $\GL_2(\A)$, obtained from a Maa{\ss} form. We
want to obtain an integral formula for the $L$-function $L(s,\pi
\times \tau)$. We will use the local calculations from the previous
two sections to achieve this.

Given a quadratic field extension $L/\Q$, we define the groups
$G=GU(2,2)$, $H=\GSp_4$, $P=MN$ and $R=TU$ as in Sect.\ \ref{unitarygroupsec}
and \ref{besselsubgroupsec}, but now considered as algebraic groups over $\Q$.
\subsection{Siegel modular forms and Bessel models}\label{siegel modular form}
Let $\Gamma_2 = \SSp_4(\Z)$. For a positive integer $l$ denote by
$S_l(\Gamma_2)$ the space of Siegel cusp forms of degree $2$ and
weight $l$ with respect to $\Gamma_2$. If $\Phi \in S_l(\Gamma_2)$
then $\Phi$ satisfies
$$
 \Phi(\gamma\langle Z \rangle) = \det(J(\gamma,Z))^l \Phi(Z),
 \qquad \gamma \in \Gamma_2,\;Z \in \SH_2.
$$
Let us assume that $\Phi \in S_l(\Gamma_2)$ is a Hecke eigenform. It has a Fourier expansion
$$
 \Phi(Z) = \sum\limits_{S > 0}a(S,\Phi)e^{2 \pi i \tr(SZ)},
$$
where $S$ runs through all symmetric semi-integral positive
definite matrices of size two. Let us make the following two
assumptions about the function $\Phi$.
\begin{description}
 \item[Assumption $1$:] $a(S,\Phi) \neq 0$ for some $S =
  \mat{a}{b/2}{b/2}{c}$ such that $b^2-4ac = -D < 0$ where $-D$
  is the discriminant of the imaginary quadratic field $\Q(\sqrt{-D})$.
 \item[Assumption $2$:] The weight $l$ is a multiple of $w(-D)$, the
  number of roots of unity in $\Q(\sqrt{-D})$. We have
  $$
   w(-D) = \left\{\begin{array}{ll}
    4 & \hbox{ if } D=4, \\
    6 & \hbox{ if } D=3, \\
    2 & \hbox{ otherwise. }\end{array}\right.
  $$
\end{description}
Let us define a function $\phi= \phi_{\Phi}$ on $H(\A) = \GSp_4(\A)$ by
\begin{equation}\label{lift of siegel modular form to group}
 \phi(\gamma h_{\infty} k_0) =
 \mu_2(h_{\infty})^l\det(J(h_{\infty},I))^{-l} \Phi(h_\infty \langle I \rangle),
\end{equation}
where $\gamma \in H(\Q)$, $h_{\infty} \in H^{+}(\R)$, $k_0 \in
\prod\limits_{p < \infty}H(\Z_p)$. Here $I = \mat{i}{}{}{i}$.
Note that $\phi$ is invariant under the center $Z_H(\A)$ of $H(\A)$.
It can be shown (see \cite[p.\ 186]{AS}) that the function $\phi_{\Phi}$
is a cuspidal automorphic form.
Let $V_{\Phi}$ be the automorphic representation generated by $\phi_{\Phi}$.
This representation may not be irreducible, but
decomposes into a direct sum of finitely many irreducible,
cuspidal, automorphic representations of $H(\A)$.
Let $\pi_{\Phi}$ be one of these irreducible components, and write $\pi_{\Phi}$ as
a restricted tensor product $\pi_{\Phi}\cong \otimes'_p \pi_{p}$, where $\pi_{p}$ is an
irreducible, admissible, unitarizable representation of $H(\Q_p)$. Since
$\phi_{\Phi}$ is $H(\Z_p)$-invariant for all finite primes $p$, the representation $\pi_{p}$ has a non-zero, essentially unique $H(\Z_p)$-invariant vector.
The same calculations as in \cite{AS} show that the equivalence
class of $\pi_p$ depends only on $\Phi$ and not on the chosen global
irreducible component $\pi_{\Phi}$.

Let $\psi = \prod\limits_{p}\psi_{p}$ be a character of $\Q
\backslash \A$ which is unramified at every finite prime and such that
$\psi_{\infty}(x) = e^{-2\pi ix}$ for $x \in \R$. Let
\begin{equation}\label{special four coeff matrix}\renewcommand{\arraystretch}{1.2}
 S(-D) = \left\{\begin{array}{l@{\qquad}l}
    \mat{\frac D4}{0}{0}{1}& \hbox{ if } D \equiv 0 \pmod{4}, \\[3ex]
 \mat{\frac{1+D}{4}}{\frac 12}{\frac 12}{1}& \hbox{ if } D \equiv 3 \pmod{4}.
 \end{array}\right.
\end{equation}
Our quadratic extension is $L=\Q(\sqrt{-D})$. We have $T(\Q)
\simeq \Q(\sqrt{-D})^{\times}$. Let $\Lambda$ be an ideal class
character of $\Q(\sqrt{-D})$, i.e., a character of
$$T(\A)/T(\Q)T(\R) \prod\limits_{p<\infty}(T(\Q_p)\cap \GL_2(\Z_p)),$$
to be chosen below. We define the global
Bessel function of type $(\Lambda, \theta)$ associated to $\bar{\phi}$ by
\begin{equation}\label{global Bessel model defn}
B_{\bar{\phi}}(h) = \int\limits_{Z_H(\A)R(\Q)\backslash
R(\A)}(\Lambda \otimes \theta)(r)^{-1}\bar{\phi}(rh)dr,
\end{equation}
where $\theta(\mat{1}{X}{}{1})=\psi(\tr(S(-D)X))$ and
$\bar{\phi}(h) = \overline{\phi(h)}$. For a finite prime $p$, the
function $B_p(h_p) := B_{\bar{\phi}}(h_p)$, with $h_p \in
H(\Q_p)$, is in the Bessel model for $\pi_p$ with respect to the
character $\Lambda_p \otimes \theta_p$ of $R(\Q_p)$. The
uniqueness of the Bessel model for $\GSp_4$ (see \cite{N-PS})
gives us
\begin{equation}\label{global-local bessel model}
 B_{\bar{\phi}}(h) = B_{\bar{\phi}}(h_{\infty}) \prod\limits_{p < \infty} B_p(h_p),
\end{equation}
where $h = \otimes h_p$. From \cite[(1-17), (1-19), (1-26)]{Su},
we have, for $h_{\infty} \in H^+(\R)$,
\begin{equation}\label{Bessel model arch formula}
 B_{\bar{\phi}}(h_{\infty}) = |\mu_2(h_{\infty})|^l\,
 \overline{\det(J(h_{\infty},I))^{-l}}\,e^{-2\pi i\,
 \tr(S(-D)\overline{h_{\infty}\langle I\rangle})}
 \sum\limits_{j=1}^{h(-D)} \Lambda(t_j)^{-1}\overline{a(S_j,\Phi)},
\end{equation}
and $B_{\bar{\phi}}(h_{\infty}) = 0$ for $h_{\infty} \not\in H^+(\R)$.
Here, $h(-D)$ is the class number of $\Q(\sqrt{-D})$, the elements
$t_j$, $j = 1,\ldots,h(-D)$, are
representatives of the idele classes of $\Q(\sqrt{-D})$ and $S_j$,
$j = 1,\cdots,h(-D)$, are the representatives of $\SL_2(\Z)$
equivalent classes of primitive semi-integral positive definite
matrices of discriminant $-D$ corresponding to $t_j$. Thus, by
Assumption 1, there exists a $\Lambda$ such that
$B_{\bar{\phi}}(I_4) \neq 0$. We fix such a $\Lambda$. Note that
$B_{\bar{\phi}}(h_{\infty})$ is a non-zero constant multiple of
(\ref{archBesselformula2eq}). Let us abbreviate $a(\Lambda) =
\sum\limits_{j=1}^{h(-D)} \Lambda(t_j)a(S_j,\Phi)$.
\subsection{Maa{\ss} forms and Eisenstein series}\label{maass form sec}
Let $\SH_1 = \{z=x+iy \in \C : y >0\}$ be the complex upper half
plane. Fix a square-free integer $N$. Let $\Gamma_0(N) =
\{\mat{a}{b}{c}{d} \in \SL_2(\Z) : N|c \}$. A smooth function $f :
\SH_1 \rightarrow \C$ is called a Maa\ss \, cusp  form of weight
$l_1$ with respect to $\Gamma_0(N)$ if
\begin{enumerate}
\item For every $\mat{a}{b}{c}{d} \in \Gamma_0(N)$ and $z \in
\SH_1$ we have
$$f\Big(\frac{a z + b}{c z + d}\Big) = \Big(\frac{c z +
d}{|c z + d|}\Big)^{l_1} f(z).$$

\item $f$ is an eigenfunction of $\Delta_{l_1}$, where
$$\Delta_{l_1} = y^2 \Big(\frac{\partial^2}{\partial x^2} + \frac{\partial^2}{\partial
y^2}\Big) - i l_1 y \frac{\partial}{\partial x}.$$

\item $f$ vanishes at the cusps of $\Gamma_0(N)$.
\end{enumerate}
Denote the space of Maa{\ss} cusp forms of weight $l_1$ with
respect to $\Gamma_0(N)$ by $S_{l_1}^M(N)$. A function $f \in
S_{l_1}^M(N)$ has the Fourier expansion
\begin{equation}\label{Maass form four exp}
 f(x+iy) = \sum\limits_{n \neq 0} a_n W_{{\rm sgn}(n)\frac{l_1}2,
\frac{ir}2}(4 \pi|n|y)e^{2 \pi i nx},
\end{equation}
where $W_{\nu,\mu}$ is a classical Whittaker function (the same function
as in (\ref{W1Whittakerformulaeq})) and $(\Delta_{l_1}+\lambda) f=0$
with $\lambda = 1/4 + (r/2)^2$. Let $f \in S_{l_1}^M(N)$ be a  Hecke eigenform.

If $ir/2 = (l_2-1)/2$ for some integer $l_2 > 0$, then the cuspidal, automorphic
representation of $\GL_2(\A)$ constructed below is holomorphic at infinity of
lowest weight $l_2$. In this case we make the additional assumptions that
$l_2\leq l$ and $l_2\leq l_1$, where $l$ is the weight of the Siegel
cusp form $\Phi$ from the previous section.

Starting from $f$, we obtain another Maa{\ss} form $f_l \in
S_{l}^M(N)$ by applying the raising and lowering operators. The
raising operator $R_{\ast}$ maps $S_{\ast}^M(N)$ to
$S_{\ast+2}^M(N)$ and the lowering operator $L_{\ast}$ maps
$S_{\ast}^M(N)$ to $S_{\ast-2}^M(N)$; for more details on these
operators, see \cite[p.\ 3925]{P}. In particular, we have
\begin{equation}\label{raise-lower-maass form}
 f_l = \left\{\begin{array}{ll}
 R_{l-2}R_{l-4}\cdots R_{l_1+2}R_{l_1}f& \hbox{ if } l_1 < l,\\
 f& \hbox{ if } l_1 = l,\\
 L_{l+2}L_{l+4}\cdots L_{l_1-2}L_{l_1}f& \hbox{ if } l_1 > l.\end{array}\right.
\end{equation}
Note that, by Assumption 2 on the Siegel cusp form $\Phi$, the weight $l$
is always even. Also, $S_{l_1}^M(N)$ is empty if $l_1$ is odd.
Hence, (\ref{raise-lower-maass form}) makes sense.
If $ir/2 = (l_2-1)/2$, then the assumption $l_2\leq l$ guarantees that $f_l\neq0$.
Suppose $\{c(n) \}$ are the Fourier coefficients of $f_l$. In
later calculations we will need $c(1)$. By \cite[Lemma 2.5]{P},
\begin{equation}\label{1st coefficient}
 c(1) = \left\{\begin{array}{ll}
  a_1& \hbox{ if } l_1 \leq l, \\
  \displaystyle\prod\limits_{\substack{t=l+2 \\ t \equiv l \pmod{2}}}^{l_1}
  \Big(\frac{ir}2 + \frac 12 - \frac t2 \Big)
  \Big(\frac{ir}2 - \frac 12 + \frac t2 \Big)a_1& \hbox{ if } l_1 > l.
 \end{array}\right.
\end{equation}
Define a function $\hat{f}$ on $\GL_2(\A)$ by
\begin{equation}\label{maass form lift to group}
 \hat{f}(\gamma_0 m k_0) = \Big(\frac{\gamma i +\delta}{|\gamma i +\delta|}\Big)^{-l}
 f_l\Big(\frac{\alpha i + \beta}{\gamma i + \delta}\Big),
\end{equation}
where $\gamma_0 \in \GL_2(\Q)$, $m =
\mat{\alpha}{\beta}{\gamma}{\delta} \in \GL_2^{+}(\R)$, $k_0 \in
\prod\limits_{p | N}K^{(1)}(\p) \prod\limits_{p \nmid
N}\GL_2(\Z_p)$. Here, for $p | N$ we have
$K^{(1)}(\p) = \GL_2(\Q_p) \cap\mat{\Z_p^{\times}}{\Z_p}{\p}{1+\p}$
with $\p = p\Z_p$, as in (\ref{K1defeq}). Then $\hat f$ satisfies
\begin{equation}\label{hatfrighttransformationpropertyeq}
 \hat f(gr(\theta))=e^{il\theta}\hat f(g),\qquad g\in\GL_2(\A),\;
 r(\theta)=\mat{\cos(\theta)}{\sin(\theta)}{-\sin(\theta)}{\cos(\theta)}.
\end{equation}
Let $(\tau_f,V_f)$ be the cuspidal, automorphic representation of
$\GL_2(\A)$ generated by $\hat f$. By strong multiplicity one,
$\tau_f$ is irreducible. Note that $\tau_f$ has trivial central
character. Write $\tau_f$ as a restricted tensor product $\tau_f =
\otimes'_p \tau_p$. If $p \nmid N$ is a finite prime, then
$\tau_p$ is an irreducible, admissible, unramified representation
of $\GL_2(\Q_p)$. If $p | N$, then $\tau_p$ is an irreducible,
admissible representation of $\GL_2(\Q_p)$ with conductor $\p =
p\Z_p$. Since $\tau$ has trivial central character, $\tau_p$ is a
twisted Steinberg representation given by
$\tau_p=\Omega_p\St_{\GL_2(\Q_p)}$, where $\Omega_p$ is an
unramified, quadratic character of $\Q_p^\times$. Let
$$
 W^{(0)}(g) := \int\limits_{\Q \backslash \A}\hat{f}(\mat{1}{x}{}{1}g)\psi(x)dx,
$$
where $\psi$ is the additive character fixed in the previous section.
Then $W^{(0)}$ is in the Whittaker model of $\tau_f$ with respect
to the character $\psi^{-1}$. By (\ref{hatfrighttransformationpropertyeq}),
\begin{equation}\label{W0righttransformationpropertyeq}
 W^{(0)}(gr(\theta))=e^{il\theta}W^{(0)}(g),\qquad g\in\GL_2(\A),\;
 r(\theta)=\mat{\cos(\theta)}{\sin(\theta)}{-\sin(\theta)}{\cos(\theta)}.
\end{equation}
For any finite prime $p$, the function $W_p(g_p) := W^{(0)}(g_p)$, for
$g_p \in \GL_2(\Q_p)$, is in the Whittaker model of $\tau_p$. By
the uniqueness of Whittaker models for $\GL_2$, we get
$$W^{(0)}(g) = W^{(0)}(g_{\infty})\prod\limits_{p<\infty}W_p(g_p)$$
for $g = \otimes g_p$. Using the definition (\ref{maass form lift to group}) for $\hat{f}$ we get, for $t \in \R^{\times}$,
\begin{equation}\label{W0t1transformationpropertyeq}
 W^{(0)}(\mat{t}{}{}{1}) = \renewcommand{\arraystretch}{1.2}
  \left\{\begin{array}{ll}
 c(1)W_{\frac l2, \frac{ir}2}(4\pi t) & \hbox{ if } t > 0, \\
 c(-1)W_{-\frac l2, \frac{ir}2}(-4\pi t)& \hbox{ if } t <0 .\end{array}\right.
\end{equation}
We want to extend $\hat{f}$ to a function on $\GU(1,1;L)(\A)$. For
this, we need to construct a suitable character $\chi_0$ on
$L^{\times} \backslash \A_L^{\times}$.

\begin{lemma}\label{divisiblegrouplemma}
 Let $S$ be a divisible group, i.e., a group with the property that
 $S=\{s^n:\:s\in S\}$ for all positive integers $n$. Let $A$ and $B$ be abelian
 groups, and assume that $B$ is finite. Then every exact sequence
 $$
  1\longrightarrow S\longrightarrow A\longrightarrow B\longrightarrow1
 $$
 splits.
\end{lemma}
{\bf Proof.} Write $B$ as a product of cyclic groups $\langle b_i\rangle$.
Choose pre-images $a_i$ of $b_i$ in $A$. Modifying $a_i$ by suitable elements
of $S$, we may assume that $a_i$ has the same order as $b_i$. Then the group
generated by all $a_i$ is isomorphic to $B$.\qed

\vspace{1ex}
\begin{lemma}\label{globalchi0lemma}
 Let $L=\Q(\sqrt{-D})$ with $D>0$ be an imaginary quadratic number field. Let $\A^\times_L$
 be the group of ideles of $L$. Let $K_f$ be the subgroup given by
 $\prod_{v<\infty}\OF_{L,v}^\times$, where $v$ runs over all finite places of $L$,
 and $\OF_{L,v}$ is the ring of integers in the completion of $L$ at $v$.
 The archimedeam component of $\A^\times_L$ is isomorphic to
 $\C^\times=\R_{>0}\times S^1$, where $S^1$ is the unit circle.
 Let $l\in\Z$ be a multiple of $w(-D)$, the number of roots of unity in $L$.
 Then there exists a character $\chi_0$ of $\A^\times_L$ with the properties
 \begin{enumerate}
  \item $\chi_0$ is trivial on $\A_\Q^\times K_f L^\times$; and
  \item $\chi_0(\zeta)=\zeta^{-l}$ for all $\zeta\in S^1$.
 \end{enumerate}
\end{lemma}
{\bf Proof.} First note that $\A_\Q^\times K_f L^\times=\R_{>0}K_f L^\times$.
There is an exact sequence
$$
 1\longrightarrow W\backslash S^1\longrightarrow \R_{>0}K_fL^\times\backslash\A^\times_L
 \longrightarrow\C^\times K_fL^\times\backslash\A^\times_L\longrightarrow1,
$$
where $W$ is the group of roots of unity in $L$. The group on the right
is the ideal class group of $L$. By Lemma \ref{divisiblegrouplemma},
$$
 \R_{>0}K_fL^\times\backslash\A^\times_L\cong(W\backslash S^1)\times
 (\C^\times K_fL^\times\backslash\A^\times_L).
$$
By hypothesis, the map $S^1\ni\zeta\mapsto\zeta^l$ factors through
$W\backslash S^1$. The assertion follows.\qed

\vspace{3ex}
Let $\chi_0$ be a character of $\A^\times_L$ as in Lemma \ref{globalchi0lemma}
(observe our Assumption 2 above).
We extend $\hat{f}$ to $\GU(1,1;L)(\A)$ by
\begin{equation}\label{hatfextensioneq}
 \hat{f}(\zeta g) = \chi_0(\zeta)\hat{f}(g)
 \qquad\text{for }\zeta \in\A^\times_L,\;g \in \GL_2(\A).
\end{equation}
Since $l$ is even, this is well-defined; see (\ref{realKintersect3eq}) and
(\ref{hatfrighttransformationpropertyeq}).
Let $\chi$ be the character of $L^\times\backslash\A^\times_L$ given by
$\chi(\zeta)=\Lambda(\bar\zeta)^{-1}\chi_0(\bar\zeta)^{-1}$.
Let $K^\#_G(N)$ be the compact subgroup $\prod\limits_{p | N}K^\#(\P)
\prod\limits_{p \nmid N}K^\#(\P^0)$ of $\GU(2,2;L)(\A)$, where
$K^\#(\P^n)$ is as defined in (\ref{K2defeq}).
For a complex variable $s$, let us define a function $f_{\Lambda}(\,\cdot\,,s)$ on
$\GU(2,2;L)(\A)$ by
\begin{enumerate}
 \item $f_{\Lambda}(g,s) = 0$ if $g \not\in M(\A)N(\A)K_{\infty}K^\#_G(N)$.
 \item If $m = m_1m_2$, $m_i \in M^{(i)}(\A)$, $n \in N(\A)$, $k =
  k_0k_{\infty}$, $k_0 \in K_G^\#(N)$, $k_{\infty} \in K_{\infty}$, then
  \begin{equation}\label{section definition}
   f_{\Lambda}(mnk,s) = \delta_P^{\frac 12 + s}(m)\chi(m_1)
   \hat{f}(m_2) \det(J(k_{\infty},I))^{-l}.
  \end{equation}
  Recall from (\ref{deltaPformulaeq})
  that $\delta_P(m_1m_2) =|N_{L/\Q}(m_1)\mu_1(m_2)^{-1}|^3$.
\end{enumerate}
Here, $M^{(1)}(\A)$, $M^{(2)}(\A)$, $N(\A)$ are the adelic points of the
algebraic groups defined by (\ref{M1defn}), (\ref{M2defn}) and (\ref{Ndefn}).
The groups $K^\#(\P)$ are as defined in (\ref{K2defeq}), and $K_{\infty}$ is as defined
in Sect.\ \ref{real groups}.
In fact, $f_{\Lambda}$ is a section in the representation $I(s,\chi,\chi_0,\tau)$
of $\GU(2,2;L)(\A)$ obtained by parabolic induction from $P$;
see Sect.\ \ref{parabolicinductionsec}.

Let us define the Eisenstein series on $\GU(2,2;L)(\A)$ by
\begin{equation}\label{eisenstein series definition}
 E_{\Lambda}(g,s) = \sum\limits_{\gamma \in P(\Q) \backslash
 G(\Q)}f_{\Lambda}(\gamma g,s).
\end{equation}
This series is absolutely convergent for ${\rm Re}(s) > 1/2$, uniformly
convergent in compact subdomains and has a meromorphic
continuation to the whole complex plane; see \cite{L}.

{\bf Remark:} Note that our definition (\ref{section definition}) differs from the
formula for $f_\Lambda$ given on p.\ 209 of \cite{Fu}. In fact, the function
$f_\Lambda$ in \cite{Fu} is not well-defined, since there is a non-trivial
overlap between $M^{(2)}(\R)$ and $K_\infty$. It is necessary to extend the
function $\hat f$ to $\GU(1,1;L)(\A)$ using the character $\chi_0$ as in
(\ref{hatfextensioneq}), not the trivial character.
\subsection{Global integral and $L$-functions}
Let $\phi$ be as in (\ref{lift of siegel modular form to group}).
Let $f_\Lambda(\,\cdot\,,s)$ and $E_\Lambda(\,\cdot\,,s)$ be as in the previous section.
We shall evaluate the global integral
\begin{equation}\label{global integral calculation}
 Z(s,\Lambda) = \int\limits_{Z_H(\A)H(\Q)\backslash H(\A)}E_{\Lambda}(h,s)\bar{\phi}(h)dh.
\end{equation}
In Theorem 2.4 of \cite{Fu}, the following basic identity has been proved.
\begin{equation}\label{basic identity}
Z(s,\Lambda) = \int\limits_{R(\A) \backslash H(\A)}W_{\Lambda}(\eta h,
s)B_{\bar{\phi}}(h) dh,
\end{equation}
where
\begin{equation}\label{W_f defn}
W_{\Lambda}(g,s) = \int\limits_{\Q \backslash \A}
f_{\Lambda}\Big(\begin{bmatrix}1&&&\\&1&&x\\&&1&\\&&&1\end{bmatrix}g,s\Big)\psi(x)dx,
\end{equation}
$$
 \eta =\begin{bmatrix}1&0&&\\\alpha&1&&\\&&1&-\bar{\alpha}\\&&0&1\end{bmatrix},
 \qquad \alpha = \frac{b+\sqrt{-D}}{2},
$$
and $B_{\bar{\phi}}$ is as defined in (\ref{global Bessel model defn}). Note that
the value of $b$ above depends on the choice of $S(-D)$ in
(\ref{special four coeff matrix}). For the choice of $f_{\Lambda}$ in the
previous section, we get
$$
 W_{\Lambda}(g,s) = W_\infty(g_{\infty},s) \prod\limits_{p<\infty}W_p(g_p,s),
$$
where $W_p$ is the function $W^\#$ defined in Sect.\ \ref{W sharp function}.
For $g_{\infty} \in G(\R)$, the function $W_\infty(g_{\infty},s)$ is exactly the function
$W^\#$ from (\ref{realWsharpdefeq}). Note that, in this case, the values of $a^+, a^-$ in
(\ref{W1Whittakerformulaeq}) are given by $a^+ = c(1)$ and $a^- = c(-1)$.
From the basic identity (\ref{basic identity}) we therefore have
$$
 Z(s,\Lambda) = \prod_{p\leq\infty}Z_p(s),\qquad Z_p(s)=
 \int\limits_{R(\Q_p)\backslash H(\Q_p)}W_p(\eta h_p,s)B_p(h_p) dh_p.
$$
Here, $B_\infty$ is the function given in (\ref{Bessel model arch
formula}). If $p$ is a finite prime such that $p \nmid N$, then
all the local data satisfies the hypothesis of Theorem 3.7 from
\cite{Fu}, where the corresponding local integral is computed. For
$p | N$, we apply Theorem \ref{localintegralmaintheorem}, and for
the archimedean integral we apply Theorem \ref{archmaintheorem}.
We obtain the following integral representation.
\begin{thm}\label{main-global-thm}
 Let $\Phi \in S_l(\Gamma_2)$ be a cuspidal Siegel eigenform of degree $2$ and even
 weight $l$ satisfying the two assumptions from Sect.\ \ref{siegel modular form}.
 Let $L=\Q(\sqrt{-D})$, where $D$ is as in Assumption 1.
 Let $N$ be a square-free, positive integer. Let $f$ be a Maa{\ss}
 Hecke eigenform of weight $l_1 \in \Z$ with respect to $\Gamma_0(N)$. If $f$ lies in a
 holomorphic discrete series with lowest weight $l_2$, then assume that $l_2 \leq l$.
 Then the integral (\ref{global integral calculation}) is given by
 \begin{equation}\label{integral-l-fn-formula}
  Z(s,\Lambda) = \kappa_{\infty} \kappa_N \frac{L(3s+\frac 12, \pi_{\Phi} \times
  \tau_f)}{\zeta(6s+1) L(3s+1, \tau_f\times \AI(\Lambda))},
 \end{equation}
 where
 \begin{align*}
  \kappa_{\infty} &=\frac12\overline{a(\Lambda)}c(1)\pi D^{-3s-\frac l2}\,
   (4 \pi)^{-3s+\frac 32-l}\,\frac{\Gamma(3s+l-1+\frac{ir}2) \Gamma(3s+l-1-\frac{ir}2)}
   {\Gamma(3s+\frac{l+1}2)}, \\
   \kappa_N &=\prod\limits_{p | N} \frac{p(p-1)}{(p+1)(p^4-1)}
  (1-\Big(\frac Lp\Big) p^{-1})(1-p^{-6s-1})^{-1}.
 \end{align*}
 Here, the non-zero constant $c(1)$ is given by (\ref{raise-lower-maass form}),
 the non-zero constant $a(\Lambda)$ is defined at the end of
 Sect.\ \ref{siegel modular form}, and
 $$
  \Big(\frac Lp\Big)=\left\{\begin{array}{l@{\qquad\text{if }p}l}
  -1&\text{ is inert in }L,\\
  0&\text{ ramifies in }L,\\
  1&\text{ splits in }L.
  \end{array}\right.
 $$
 The quantity $\frac{ir}2$ is as in (\ref{Maass form four exp}).
\end{thm}
\subsection{The special value}\label{special values sec}
In this section, we will apply Theorem \ref{main-global-thm} to a
special case -- when $f$, from the previous section, is a
holomorphic cusp form of the same weight $l$ as the Siegel cusp
form $\Phi$ -- to obtain a special $L$-value result. This result fits
into the general conjecture of Deligne on special values of $L$-functions.

Let $\Psi(z) = \sum\limits_{n>0}b_ne^{2\pi inz}$ be a holomorphic
cuspidal eigenform of weight $l$ with respect to $\Gamma_0(N)$.
Here, $l$ is the same as the weight of the Siegel modular form
$\Phi$ from Sect.\ \ref{siegel modular form} and $N$ is a
square-free, positive integer. Let us normalize $\Psi$ so that
$b_1=1$. The function $f_{\Psi}$ defined by
$f_{\Psi}(z)=y^{l/2}\Psi(z)$ is a Maa{\ss} form in $S_l^M(N)$. Let
$\{c(n)\}$ be its Fourier coefficients; see (\ref{Maass form four
exp}). It follows from the formula
$W_{\mu+1/2,\mu}(z)=e^{-z/2}z^{\mu+1/2}$ for the Whittaker
function that
\begin{equation}\label{holom-maass-four-reln}
 c(n) = \left\{\begin{array}{ll}
           b_n(4\pi n)^{-l/2} & \hbox{ if } n>0, \\
           0& \hbox{ if } n<0.
         \end{array}\right.
\end{equation}
From (\ref{maass form lift to group}), we have
$$
 \hat f_{\Psi}(\gamma_0 m k_0) = \Big(\frac{\gamma i +\delta}
 {|\gamma i + \delta|}\Big)^{-l}f_{\Psi}\Big(\frac{\alpha i
 + \beta}{\gamma i + \delta}\Big) = \frac{\det(m)^{l/2}}{(\gamma i
 + \delta)^l}\Psi\Big(\frac{\alpha i + \beta}{\gamma i +\delta}\Big),
$$
where $\gamma_0 \in \GL_2(\Q)$, $m =
\mat{\alpha}{\beta}{\gamma}{\delta} \in \GL_2^{+}(\R)$, $k_0 \in
\prod\limits_{p | N}K^{(1)}(\p) \prod\limits_{p \nmid N}\GL_2(\Z_p)$.
Let us denote $z_{22}$ by $Z^\ast$ for $Z =
\mat{\ast}{\ast}{\ast}{z_{22}} \in \HH_2$. Let us set $\hat{Z} =
\frac i2(\,^t\!\bar{Z}-Z)$ for $Z \in \HH_2$. Let ${\rm Im}(z)$ denote the
imaginary part of a complex number $z$. Let $f_{\Lambda}$ be as
defined in (\ref{section definition}) and $I=\mat{i}{}{}{i} \in\HH_2$.
\begin{lem}
 For $g \in G^+(\R)$, we have
 \begin{equation}\label{alternative f-lambda formula}
  f_{\Lambda}(g,s) = \mu_2(g)^l\det(J(g,I))^{-l}\Big(\frac{\det
  \widehat{g\langle I \rangle}}{{\rm Im}(g \langle I
  \rangle)^\ast}\Big)^{3s+\frac 32 - \frac l2}\Psi((g\langle I\rangle)^\ast).
 \end{equation}
\end{lem}
{\bf Proof.} For $g \in G^+(\R)$ and $Z \in \HH_2$ we have
$\widehat{g\langle Z \rangle} = \mu_2(g)\,
^t\overline{J(g,Z)}^{-1} \hat{Z} J(g,Z)^{-1}$. This implies that
$\det(\widehat{ g\langle I \rangle}) = \mu_2(g)^2
|\det(J(g,I))|^{-2} \det(\hat{I}) = \mu_2(g)^2
|\det(J(g,I))|^{-2}$. It follows from (\ref{realKintersect1eq}) that we can
write the element $g \in G^+(\R)$ as
$$
 g = \begin{bmatrix}\zeta\\&1\\&&\zeta^{-1}\\&&&1\end{bmatrix}
 \begin{bmatrix}1\\&\alpha&&\beta\\&&\mu\\&\gamma&&\delta\end{bmatrix}
 \begin{bmatrix}1&x&x\bar{y}+w&y\\ &1&\bar{y}\\ &&1\\
 &&-\bar{x}&1\end{bmatrix} k,
$$
where $\zeta \in \R^{\times}$, $\mat{\alpha}{\beta}{\gamma}{\delta} \in \GL_2^+(\R)$,
$x,y \in \C$, $w\in \R$ and $k \in K_{\infty}$. Then we have
$$
 \det(J(g,I)) = \zeta^{-1}\mu (\gamma i + \delta)
 \det(J(k,I)) \mbox{ and } (g \langle I \rangle)^\ast =
 \frac{\alpha i + \beta}{\gamma i + \delta}.
$$
Hence, the right hand side of (\ref{alternative f-lambda formula}) is equal to
$$
 \mu^l \big(\zeta^{-1} \mu (\gamma i + \delta) \det(J(k,I))\big)^{-l}
 \Big(\frac{\mu^2|\zeta^{-1}\mu (\gamma i +
 \delta)\det(J(k,I))|^{-2}}{\mu |\gamma i +
 \delta|^{-2}}\Big)^{3s+\frac 32 - \frac l2} \Psi\Big(\frac{\alpha
 i + \beta}{\gamma i + \delta}\Big).
$$
Using the fact that $|\det(J(k,I))|^{-2} = \det(\widehat{k\langle I \rangle}) = 1$, we
get the lemma. \qed

\vspace{3ex}
{\bf Remark:} Eq.\ $(4.4.2)$ of \cite{Fu} claims that, for $g \in G^+(\R)$, the
function $f_{\Lambda}(g,s)$ satisfies a formula different
from (\ref{alternative f-lambda formula}). In this formula,
the term $\det({\rm Im}(g\langle I \rangle))$ replaces the term
$\det(\widehat{g\langle I \rangle})$ from
(\ref{alternative f-lambda formula}). Note that ${\rm Im}(Z) = \frac i2(\bar{Z}-Z)$ for
$Z \in\HH_2$. One easily checks that the resulting function is not
invariant under $N(\R)$ and hence cannot equal $f_{\Lambda}$, as
defined in (\ref{section definition}). If one replaces Eq.\
$(4.4.2)$ in \cite{Fu} by (\ref{alternative f-lambda formula}),
the subsequent arguments in \cite{Fu} remain valid.

Let $E_{\Lambda}$ be the Eisenstein series defined in
(\ref{eisenstein series definition}). From the above lemma, we see
that, for $g \in G^+(\R)$,
$\mu_2(g)^{-l}\det(J(g,I))^lE_{\Lambda}(g,s)$ only depends on $Z =
g \langle I \rangle$. Hence, we can define a function $\mathcal
E_{\Lambda}$ on $\HH_2$ by the formula
\begin{equation}\label{eis ser on siegel half plane}
 \mathcal E_{\Lambda}(Z,s) = \mu_2(g)^{-l}\det(J(g,I))^l\,
 E_{\Lambda}\big(g,\frac s3 + \frac l6 - \frac 12\big),
\end{equation}
where $g \in G^+(\R)$ is such that $g\langle I \rangle = Z$. The
series that defines $\mathcal E_{\Lambda}(Z,s)$ is absolutely
convergent for ${\rm Re}(s) > 3 - l/2$ (see \cite{Kl}).  Since $l
\geq 12$, we can set $s=0$ and obtain a holomorphic Eisenstein
series $\mathcal E_{\Lambda}(Z,0)$ on $\HH_2$.
For a finite place $p$ of $\Q$ recall the local congruence subgroups
$K^\#(\P^n)\subset G(\Z_p)$ and $K^\#(\p^n)=K^\#(\P^n)\cap H(\Z_p)$
defined in (\ref{K2defeq}) resp.\ (\ref{K2capHeq}). We let
$$
 \Gamma^\#_G(N) = G(\Q)\cap G(\R)^+K_G^\#(N),\qquad
 K_G^\#(N)=\prod_{p|N}K^\#(\P)\prod_{p\nmid N}K^\#(\P^0),
$$
and
$$
 \Gamma^\#(N) = H(\Q)\cap H(\R)^+K^\#(N),\qquad
 K^\#(N)=\prod_{p|N}K^\#(\p)\prod_{p\nmid N}K^\#(\p^0).
$$
Since the functions $f_\Lambda$ and $E_\Lambda$ are also invariant under the center, we let
$$
 \tilde\Gamma^\#(N) = H(\Q)\cap H(\R)^+\tilde K^\#(N),\qquad
 \tilde K^\#(N)=\prod_{p|N}(Z_H(\Z_p)K^\#(\p))\prod_{p\nmid N}(Z_H(\Z_p)K^\#(\p^0)).
$$
Explicitly,
$$
 \tilde\Gamma^\#(N)=\{h=(h_{ij})\in{\rm Sp}(4,\Z)\cap
 \begin{bmatrix}\Z&N\Z&\Z&\Z\\\Z&\Z&\Z&\Z\\N\Z&N\Z&\Z&\Z\\N\Z&N\Z&N\Z&\Z\end{bmatrix}:
 \:h_{11}\equiv h_{44}\;{\rm mod}\;N,\;h_{22}\equiv h_{33}\;{\rm mod}\;N\}.
$$
Then $\mathcal E_{\Lambda}(Z,0)$ is a modular form  of weight $l$
with respect to $\Gamma^\#_G(N)$. Its restriction to $\mathfrak{h}_2$ is a modular form
of weight $l$ with respect to $\tilde\Gamma^\#(N)$.
We see that $\mathcal E_{\Lambda}(Z,0)$ has a Fourier expansion
$$
 \mathcal E_{\Lambda}(Z,0) = \sum\limits_{\mathcal S \geq 0} b(\mathcal S, \mathcal
 E_{\Lambda})e^{2 \pi i \,\tr(\mathcal SZ)},
$$
where $\mathcal{S}$ runs through all Hermitian half-integral (i.e., $\mathcal S =
\mat{t_1}{\bar{t}_2}{t_2}{t_3},\: t_1, t_2 \in \Z, \:\sqrt{-D}\,t_2 \in
\mathcal{O}_{\Q(\sqrt{-D})}$) positive semi-definite matrices of size $2$. By \cite{Ha},
\begin{equation}\label{four coeff algeb}
 b(\mathcal S, \mathcal E_{\Lambda}) \in \bar{\Q} \qquad\mbox{ for any } \mathcal S.
\end{equation}
Here $\bar{\Q}$ denotes the algebraic closure of $\Q$ in $\C$. The
relation between the global integral $Z(s,\Lambda)$ defined in
(\ref{global integral calculation}) and the Eisenstein series
$\mathcal E_{\Lambda}$ is given in the following lemma.
\begin{lem}\label{integral matching formula}
 We have
 $$
  Z\big(\frac l6-\frac 12,\Lambda\big)=\frac12V_N
  \int\limits_{\tilde\Gamma^\#(N)\backslash \SH_2}
  \mathcal E_{\Lambda}(Z,0) \bar{\Phi}(Z) (\det(Y))^{l-3}\,dX\,dY,
 $$
 where $V_N = \prod\limits_{p|N}\frac{1}{(p^2-1)(p^4-1)}$ and $Z =X + iY$.
\end{lem}
{\bf Proof.} By definition,
$$
 Z\big(\frac l6 - \frac 12, \Lambda\big)=\int\limits_{Z_H(\A)H(\Q) \backslash H(\A)}
 E_{\Lambda}\big(h, \frac l6 - \frac 12\big) \bar{\phi}_{\Phi}(h)\,dh.
$$
Note that the integrand is right invariant under $K^H_\infty K^\#(N)$. Since
${\rm vol}(K_{\infty}^HK^\#(N)) =\prod\limits_{p|N}\frac{1}{(p^2-1)(p^4-1)} = V_N$,
it follows that
$$
 Z\big(\frac l6 - \frac 12, \Lambda\big)=V_N
 \int\limits_{Z_H(\A)H(\Q) \backslash H(\A)/K^H_\infty K^\#(N)}
 E_{\Lambda}\big(h, \frac l6 - \frac 12\big) \bar{\phi}_{\Phi}(h)\,dh.
$$
Note that
\begin{equation}\label{Hupperhalfplaneeq}
 Z_H(\A)H(\Q) \backslash H(\A) /K^H_\infty K^\#(N)
 = Z_H(\R)\tilde\Gamma^\#(N) \backslash H(\R)^+ /K^H_{\infty}
 = \tilde\Gamma^\#(N) \backslash \SH_2.
\end{equation}
The $H(\R)^+$-invariant measure on $\SH_2$ is given by $\frac 12\det(Y)^{-3}dXdY$. From
(\ref{lift of siegel modular form to group}) and (\ref{eis ser on siegel half plane})
we get, for $h \in H(\R)^+$,
\begin{align*}
 E_{\Lambda}\big(h,\frac l6-\frac 12\big)\bar{\phi}_{\Phi}(h) &=\mu_2(h)^l \det(J(h,I))^{-l}
 \,\mathcal E_{\Lambda}(h\langle I \rangle,0)\,\mu_2(h)^{l}\,\overline{\det(J(h,I))}^{-l}\,
 \bar{\Phi}(h\langle I \rangle) \\
 &=\det(Y)^l\,\mathcal E_{\Lambda}(Z,0)\,\bar{\Phi}(Z),
\end{align*}
where $Z = h\langle I \rangle=X+iY$. We get the last equality because, for $Z \in \SH_2$
and $h \in H(\R)^+$,
$$
 {\rm Im}(h \langle Z \rangle)
 =\mu_2(h)\,{}^t\!J(h,Z)^{-1}\,{\rm Im}(Z)\,\overline{J(h,Z)}^{-1}.
$$
This completes the proof of the lemma. \qed

\begin{lem}\label{integral-is-alg-lem}
With notations as above, we have
\begin{equation}\label{integral-is-alg}
 \frac{Z(\frac l6 - \frac 12, \Lambda)}{(\Phi,\Phi)_2} \in \bar{\Q},
\end{equation}
where
$$
 (\Phi,\Phi)_2 = \int\limits_{\Gamma_2 \backslash \SH_2}|\Phi(Z)|^2\det(Y)^{l-3}\,dX\,dY.
$$
\end{lem}
{\bf Proof.} Let $\Gamma^{(2)}(N) := \{g \in \SSp_4(\Z) : g \equiv
1 \pmod{N} \}$. Since $\Gamma^{(2)}(N) \subset \tilde\Gamma^\#(N)$
we know that $\mathcal E_{\Lambda}|_{\SH_2}$ is a holomorphic
Siegel modular form of weight $l$ with respect to
$\Gamma^{(2)}(N)$. Let us denote the space of all holomorphic
Siegel modular forms of weight $l$ with respect to
$\Gamma^{(2)}(N)$ by $M_l(\Gamma^{(2)}(N))$ and its subspace of
cusp forms by $S_l(\Gamma^{(2)}(N))$. For $\Phi_1, \Phi_2$ in
$M_l(\Gamma^{(2)}(N))$ with one of the $\Phi_i$ a cusp form, one
can define the Petersson inner product $\langle \Phi_1, \Phi_2
\rangle_N$ in the usual way. Let $V$ be the orthogonal complement
of $S_l(\Gamma^{(2)}(N))$ in $M_l(\Gamma^{(2)}(N))$ with respect
to the Petersson inner product. In Corollary 2.4.6 of \cite{Ha},
it is shown, using the Siegel operator, that $V$ is generated by
Eisenstein series. By Theorem 3.2.1 of \cite{Ha}, one can choose a
basis $\{E_j\}$ such that all the Fourier coefficients of each
$E_j$ are algebraic. From \cite[p. 460]{Ga1}, we can find an
orthogonal basis $\{\Phi_i\}$ of $S_l(\Gamma^{(2)}(N))$ such that
$\Phi = \Phi_1$ and all the Fourier coefficients of the $\Phi_i$
are algebraic. Let us write
\begin{equation}\label{E-lambda-development}
\mathcal E_{\Lambda}|_{\SH_2} = \sum\limits_i \alpha_i \Phi_i +
\sum\limits_j \beta_j E_j.
\end{equation}
 Given a $F \in M_l(\Gamma^{(2)}(N))$
and $\sigma \in Aut(\C/\bar{\Q})$, let $F^{\sigma} \in
M_l(\Gamma^{(2)}(N))$ be defined by applying the automorphism
$\sigma$ to the Fourier coefficients of $F$. Applying $\sigma$ to
(\ref{E-lambda-development}) we get
\begin{equation}\label{E-lambda-development1}
\mathcal E_{\Lambda}|_{\SH_2} = \sum\limits_i \sigma(\alpha_i)
\Phi_i + \sum\limits_j \sigma(\beta_j) E_j.
\end{equation}
 This follows from the
construction of the bases $\{\Phi_i\}, \{E_j\}$ and the property
(\ref{four coeff algeb}). From (\ref{E-lambda-development}) and
(\ref{E-lambda-development1}) we now get
$$\sigma\Big(\frac{\langle \mathcal E_{\Lambda}|_{\SH_2}, \Phi_1
\rangle}{\langle \Phi_1, \Phi_1 \rangle}\Big) = \frac{\langle
\mathcal E_{\Lambda}|_{\SH_2}, \Phi_1 \rangle}{\langle \Phi_1,
\Phi_1 \rangle} \mbox{ for all } \sigma \in Aut(\C/\bar{\Q})
\qquad \Rightarrow \qquad \frac{\langle \mathcal
E_{\Lambda}|_{\SH_2}, \Phi_1 \rangle}{\langle \Phi_1, \Phi_1
\rangle} \in \overline{\Q}.$$ Now, using Lemma \ref{integral
matching formula}, we get the result. \qed

Let $(\Psi,\Psi)_1 = (\SL_2(\Z) : \Gamma_1(N))^{-1}
\int\limits_{\Gamma_1(N)\backslash \SH_1}|\Psi(z)|^2y^{l-2}dxdy$,
where $\Gamma_1(N) := \{\mat{a}{b}{c}{d} \in \Gamma_0(N) : a, d
\equiv 1 \pmod{N} \}$. We have the following generalization of Theorem 4.8.3 of \cite{Fu}.
\begin{thm}\label{special values thm}
 Let $\Phi$ be a cuspidal Siegel eigenform of weight $l$ with respect to
 $\Gamma_2$ satisfying the two assumptions from Section \ref{siegel modular form}.
 Let $\Psi$ be a normalized, holomorphic,
 cuspidal eigenform of weight $l$ with respect to $\Gamma_0(N)$, with $N$ a square-free,
 positive integer. Then
 \begin{equation}\label{special-values-eqn}
  \frac{L(\frac l2 - 1, \pi_{\Phi} \times \tau_{\Psi})}{\pi^{5l-8}
  (\Phi,\Phi)_2(\Psi,\Psi)_1} \in \bar{\Q}
 \end{equation}
\end{thm}
{\bf Proof.} By Theorem \ref{main-global-thm}, we have
$$Z(\frac l6 - \frac 12, \Lambda) = C \pi^{4-2l}\frac{L(\frac l2 - 1, \pi_{\Phi} \times \tau_{\Psi})}{\zeta(l-2) L(\frac{l-1}2, \tau_{\Psi} \times \AI(\Lambda))},$$
where
$$
 C = \overline{a(\Lambda)}D^{-l+\frac 32}2^{-4l+6}(2l-5)!\prod\limits_{p | N}
 \frac{p(p-1)}{(p+1)(p^4-1)} (1-\Big(\frac{\Q(\sqrt{-D})}{p}\Big)
 p^{-1})(1-p^{-l+2})^{-1} \in \bar{\Q};
$$
observe that $\frac{ir}2=\frac{l-1}2$, and that $c(1)=(4\pi)^{-l/2}$ by
(\ref{holom-maass-four-reln}). It is well known that
$\zeta(l-2)\pi^{2-l} \in \Q$. Using \cite{Sh}, by the same
argument as in the proof of Theorem 4.8.3 in \cite{Fu}, we get
$$\frac{L(\frac{l-1}2, \tau_{\Psi} \times \AI(\Lambda))}{\pi^{2l-2}(\Psi,\Psi)_1} \in \bar{\Q}.$$
Together with (\ref{integral-is-alg}), this implies the theorem. \qed

We remark that it would be interesting to know the behavior of
the quantity $\frac{L(\frac l2 - 1, \pi_{\Phi} \times \tau_{\Psi})}{\pi^{5l-8}
(\Phi,\Phi)_2(\Psi,\Psi)_1}$ under the action of
${\rm Aut}(\C)$. This subject will be considered in a future work.

\end{document}